\documentclass[reqno]{amsart}
\usepackage{tikz-cd}
\usepackage{amsmath, bbm, color}
\usepackage{amssymb}
\usepackage{amsthm}
\usepackage{mathtools}
\usepackage{tikz}
\usepackage[final]{showkeys}
\usepackage[top=1in,bottom=1in,left=1in,right=1in]{geometry}
\usepackage{hyperref}
\usepackage{enumerate}
\usepackage{color}
\usetikzlibrary{patterns,positioning,arrows,shapes,trees}
\usetikzlibrary{calc}

\newtheorem{thm}{Theorem}[section]
 
\newtheorem{corollary}[thm]{Corollary}

\newtheorem{lemma}[thm]{Lemma}
\newtheorem{exm}[thm]{Example}
\newtheorem{prop}[thm]{Proposition}

\newtheorem{definition}[thm]{Definition}
\newtheorem{rem}[thm]{Remark}
\allowdisplaybreaks

\numberwithin{equation}{section}

\newcounter{counterConstant}
\newcommand{\const}[1]{
    \addtocounter{counterConstant}{1}
    \edef#1{\arabic{counterConstant}}
}

\def\R{\mathbb{R}}
\def\P{\mathbb{P}}
\def\E{\mathbb{E}}

\def\sL{\mathcal{L}}
\def\sC{\mathcal{C}}

\def\1{{\mathbbm 1}}

\begin{document}

\title{Dirichlet heat kernel estimates for rectilinear stable processes}
\date{}

\author[Z.-Q. Chen]{Zhen-Qing Chen}
\address{Zhen-Qing Chen: Department of Mathematics, University of Washington, Seattle, WA 98195, USA}
\email{zqchen@uw.edu}

\author[E. Hu]{Eryan Hu}
\address{Eryan Hu: Center for Applied Mathematics and KL-AAGDM, Tianjin University, Tianjin 300072, P.R. China}
\email{eryan.hu@tju.edu.cn}

\author[G. Zhao]{Guohuan Zhao}
\address{Guohuan Zhao: Applied Mathematics, Chinese Academy of Science, Beijing 100081, P.R. China}
\email{gzhao@amss.ac.cn}

\thanks{The research of Z.-Q. Chen is partially supported by a Simons Foundation Grant. The research of E. Hu is supported by National Key R\&D Program of China (No. 2022YFA1006000) and by the National Natural Science Foundation of China (No. 12171354). The research of G. Zhao is supported by the National Natural Science Foundation of China (Nos. 12288201, 12271352, 12201611).}

\begin{abstract}
  Let $d \ge 2$, $\alpha \in (0,2)$, and $X$ be the rectilinear $\alpha$-stable process on $\R^d$. We first present a geometric characterization of open subset $D\subset \R^d$ so that the part process $X^D$ of $X$ in $D$ is irreducible. We then study the properties of the transition density functions of $X^D$, including the strict positivity property as well as their sharp two-sided bounds in $C^{1,1}$ domains in $\mathbb{R}^d$. Our bounds are shown to be sharp for a class of $C^{1,1}$ domains.
\end{abstract}

\maketitle
\tableofcontents

\section{Introduction and main results}\label{S:1}

  Dirichlet heat kernels for non-local operators are a fundamental subject  both in analysis and in probability theory. Sharp two-sided Dirichlet heat kernel estimates for fractional Laplacian $\Delta^{\alpha/2}:= -(-\Delta)^{\alpha/2}$ in $C^{1,1}$ open subsets of $\R^d$ with $\alpha \in (0, 2)$ have first been obtained in \cite{ChenKimSong.2010.JEMSJ1307}. Since then, there are many works in extending it to certain classes of symmetric Markov processes and their lower order perturbations as well as to more general open sets; see, for example, \cite{BaeKim.2020.PA661, BogdanGrzywnyRyznar.2010.AP1901, BogdanGrzywnyRyznar.2014.SPA3612, ChenKimWang.2022.MA373, ChenKimSong.2011.JLMS258, ChenKimSong.2012.AP2483, ChenKimSong.2014.PLMS390, ChenKimSong.2016.JRAM111, ChenTokle.2011.PTRF373, ChoKangKim.2021.JLMS2823, GrzywnyKimKim.2020.SPA431, KimMimica.2018.EJP45} and the references therein. In all these works, the jump measures of the Markov processes are absolutely continuous with respect to the Lebesgue measure.

  Let $d \ge 2$ and $\alpha \in (0,2)$. The purpose of this paper is to study the Dirichlet heat kernels for
  \begin{equation*}
    \mathcal{L} := - \sum_{k=1}^d \left(- \frac{\partial^2}{(\partial x^{(k)})^2}\right)^{\alpha/2}
  \end{equation*}
  in open subsets of $\R^d$. Here $x^{(k)}$ is the $k^{th}$-coordinate of a point $x = (x^{(1)},x^{(2)},\cdots, x^{(d)}) \in \mathbb{R}^d$. We call $\sL$ a rectilinear fractional Laplace operator, which is more singular than the usual isotropic  fractional Laplacian $\Delta^{\alpha/2}$. The rectilinear fractional Laplacian  $\sL$ is the infinitesimal generator of the rectilinear $\alpha$-stable process
  $$
    X= \left\{ X_t = \big(X^{(1)}_t,X^{(2)}_t,\cdots,X^{(d)}_t \big); t\geq 0 \right\}
  $$
  on $\R^d$, where $X^{(1)}$, $X^{(2)}, \cdots, X^{(d)}$ are independent one-dimensional symmetric $\alpha$-stable processes. The process $X$ is a L\'evy process on $\R^d$ whose L\'evy exponent is $\Psi(\xi)=\sum_{j=1}^d |\xi^{(j)}|^\alpha$ for $\xi = (\xi^{(1)},\xi^{(2)}, \cdots, \xi^{(d)}) \in \R^d$; that is,
  $$
	\E e^{i\xi \cdot (X_t-X_0)} = e^{-t \sum_{j=1}^d |\xi^{(j)}|^\alpha}\quad\hbox{for } t>0 \hbox{ and } \xi \in \R^d.
  $$
  The L\'evy measure of $X$ is singular with respect to the Lebesgue measure on $\R^d$; see \eqref{e:1.5}.

\medskip

  Unlike the isotropic (or, rotationally symmetric) $\alpha$-stable process $Z$ on $\R^d$, the distribution of the increments of the rectilinear $\alpha$-stable process $X$ is not rotationally invariant. The isotropic $\alpha$-stable process $Z$ is a L\'evy process on $\R^d$ having infinitesimal generator $\Delta^{\alpha/2}$ and L\'evy exponent $|\xi|^\alpha :=\big(\sum_{j=1}^d |\xi^{(j)}|^2\big)^{\alpha/2}$. For $f\in C_c^2 (\R^d)$,
  \begin{equation} \label{e:1.0}
    \Delta^{\alpha/2} f(x)= \int_{\R^d \setminus \{0\}} \left( f(x+z) - f(x) - \nabla f(x) \cdot z \1_{\{|z|\leq 1\}} \right) \frac{\sC_{d,\alpha}}{|z|^{d+\alpha}} dz,
  \end{equation}
  while
  \begin{equation} \label{e:1.1}
    \sL f(x)= \sum_{j=1}^d \int_{\R \setminus \{0\}} \left( f(x+we_j) - f(x) - w\1_{\{|w|\leq 1\}} \frac{\partial f(x)}{\partial x^{(j)}} \right) \frac{\sC_{1,\alpha}}{|w|^{1+\alpha}} dw,
  \end{equation}
  where $e_j$ is the unit vector in the positive $x^{(j)}$-direction and
  \begin{equation} \label{e:1.2}
    \sC_{d,\alpha} = \frac{\alpha 2^{\alpha-1} \Gamma((d+\alpha)/2)}{\pi^{d/2} \Gamma(1-\alpha/2)}.
  \end{equation}
  Here $\Gamma$ is the usual Gamma function defined by $\Gamma(\lambda) := \int_0^\infty t^{\lambda-1}e^{-t}dt$ for $\lambda > 0$. Note that the formula \eqref{e:1.1} can be rewritten as:
  \begin{align}
    \sL f(x) =&~ \sum_{j=1}^d \int_{\R} \left( f(x+we_j) + f(x-we_j) - 2f(x) \right) \frac{\sC_{1,\alpha}}{2|w|^{1+\alpha}} dw, \label{e:1.1-2}\\
    =&~ \sum_{j=1}^d \lim_{\varepsilon \rightarrow 0+} \int_{|w|>\varepsilon} \left( f(x+we_j) - f(x)\right) \frac{\sC_{1,\alpha}}{|w|^{1+\alpha}} dw. \label{e:1.1-3}
  \end{align}
  Another way to see the differences, the operator $\sL$ has Fourier symbol $\Psi(\xi)=\sum_{j=1}^d |\xi^{(j)}|^\alpha$, while $\Delta^{\alpha/2}$ has Fourier symbol $|\xi|^\alpha$, which is smoother than $\Psi(\xi)$.

 \medskip

  It is well known (see \cite[Theorem 2.1]{BlumenthalGetoor.1960.TAMS263} via stable scaling) that the isotropic $\alpha$-stable process $Z$ on $\R^d$ has a smooth density function $p^{(d,\alpha)}(t,x,y)$ with respect to the Lebesgue measure on $\R^d$ and there are positive constants $c_2>c_1>0$ that depend only on $d$ and $\alpha$ so that
  \begin{equation}\label{e:1.3}
    c_1 \left( t^{-d/{\alpha}}\wedge \frac{t}{|x-y|^{d+\alpha}} \right) \leq  p^{(d, \alpha)}(t,x ,y )\leq c_2 \left( t^{-{d}/{\alpha}}\wedge \frac{t}{|x-y|^{d+\alpha}} \right) \quad \hbox{for } t>0,~x,  y \in \mathbb{R}^d.
  \end{equation}
  In this paper, we will use := as a way of definition. For $a, b\in \R$, $a\wedge b:= \min \{a, b\}$ and $a\vee b:= \max \{a, b\}$. Since for $a>0$ and $b>0$,
  $$
    \frac{ab}{a+b} \leq   a\wedge b \leq \frac{2ab}{a+b},
  $$
  we can rewrite the estimates in \eqref{e:1.3} by
  \begin{equation}\label{e:1.4}
    \frac{c_3 \, t}{(t^{1/\alpha} + |x-y|)^{d+\alpha}} \leq p^{(d,\alpha)}(t,x,y) \leq \frac{c_4 \, t}{(t^{1/\alpha} + |x-y|)^{d+\alpha}} \quad \hbox{for } t>0 \hbox{ and } x,  y \in \mathbb{R}^d,
  \end{equation}
  where constants $c_4>c_3>0$ depend only on $d$ and $\alpha$.

  \const{\CTstable}
  By the independence between its coordinate processes, the rectilinear $\alpha$-stable process $X$ on $\R^d$ has a smooth transition density function
  \begin{equation}\label{eq:ppalpha}
    p(t,x,y) = \prod_{k=1}^d p^{(1, \alpha)}(t,x^{(k)},y^{(k)})  \quad \hbox{for } t>0 \hbox{ and } x=(x^{(k)}) ,  \, y=(y^{(k)}) \in \R^d,
  \end{equation}
  with respect to the Lebesgue measure on $\R^d$. By \eqref{e:1.3}, there is a constant $C_{\CTstable} = C_{\CTstable}(d,\alpha) > 1$ so that
  \begin{equation}\label{eq:pEst}
    C_{\CTstable}^{-1} \, \prod_{k=1}^d \left(t^{-1/\alpha}\wedge\frac{t}{|x^{(k)}-y^{(k)}|^{1+\alpha}}\right)  \le p(t,x,y) \le C_{\CTstable} \,   \prod_{k=1}^d \left(t^{-1/\alpha}\wedge\frac{t}{|x^{(k)}-y^{(k)}|^{1+\alpha}}\right)
   \end{equation}
  for all $t>0$, and $x=(x^{(k)}), \, y=(y^{(k)}) \in \mathbb{R}^d$. This is clearly quite different from the estimates for the heat kernel $p^{(d,\alpha)}(t,x,y)$ of the isotropic or rotationally symmetric $\alpha$-stable process $Z$ on $\R^d$.

\medskip

  Though the rectilinear $\alpha$-stable process $X$ and the isotropic  $\alpha$-stable process $Z$ are both L\'evy processes that are invariant under the $\alpha$-stable scaling, they have many fundamentally different properties, which will be further revealed in this paper. For instance, it is shown in \cite{BassChen.2010.MZ489} that Harnack inequality fails for rectilinear stable processes, while scale invariant Harnack inequality holds for isotropic stable processes. The root of these differences lies in the fact that the isotropic $\alpha$-stable process $Z$ can jump in any direction uniformly, while the rectilinear $\alpha$-stable process $X$ can only jump along the directions of coordinate axes, one at a time, and thus is much singular. The L\'evy measure of $X$ is
  \begin{equation}\label{e:1.5}
    \mu(dz) =\sum_{j=1}^d \frac{\sC_{1, \alpha} }{|z^{(j)}|^{1+\alpha}} dz^{(j)} \otimes \prod_{\genfrac{}{}{0pt}{4}{k=1}{k\not = j}}^d \delta_{\{0\}} (dz^{(k)}),
  \end{equation}
  where $\delta_{\{0\}}$ denotes the Dirac measure concentrated at $0$ and $z=(z^{(1)},z^{(2)},\cdots,z^{(d)})\in \R^d$. The L\'evy measure $\mu$ describes how the rectilinear $\alpha$-stable process $X$ jumps. For any  non-negative measurable function $f$ on $\mathbb{R}_+\times \mathbb{R}^d\times \mathbb{R}^d$ with $f(s,x,x)=0$ for any $s\geq 0$ and $x\in \R^d$ and for any stopping time $S$ with respect to the minimum augmented filtration generated by $X$, we have
  \begin{equation}\label{eq:XLS-0}
    \mathbb{E}_x \left[\sum_{s\le S} f(s,X_{s-},X_s)\right] = \mathbb{E}_x \left[\int_0^S \int_{\R^d} f(s,X_s,X_s+z) \mu(dz) ds\right].
  \end{equation}
  See, e.g., \cite[proof of Lemma 4.7]{ChenKumagai.2003.SPA27} and \cite[Appendix A]{ChenKumagai.2008.PTRF277}. We mention that recently it is shown in \cite{KassmannKimKumagai.2022.JMPA91} that the transition density functions of the symmetric pure jump processes whose jumping measure $J(dx,dy)$ is comparable to $dx \mu(dy-x)$, where $\mu$ is the L\'evy measure of \eqref{e:1.5}, have the two-sided estimates \eqref{eq:pEst}. This result has been further extended to more general rectilinear L\'evy processes in \cite{KimWang.2022.SPA165}.

\medskip
	
  For any non-empty open subset $D\subset \mathbb{R}^d$, let $\tau_D(\omega) = \inf\{t>0,X_t(\omega) \notin D\}$ denote the first exit time from $D$ by $X$. Taking $f(x,y)=\1_D(x) \1_{D^c}(y) \varphi(y)$ and $S=\tau_D$ in \eqref{eq:XLS-0}, where $\varphi$ is a bounded function defined on $D^c$, yields
  \begin{equation} \label{e:1.10}
    \mathbb{E}_x \left[\varphi(X_{\tau_D}); X_{\tau_D-} \not= X_{\tau_D} \right] = \mathbb{E}_x \left[\int_0^{\tau_D}   \int_{D^c} \varphi(z) \mu (dz -X_s) ds\right].
  \end{equation}
  The subprocess $X^D$ of $X$ killed upon leaving $D$ is defined as
  $$
    X^{D}_t(\omega) =
    \begin{cases}
      X_t(\omega) \quad &\hbox{for } t< \tau_D(\omega) \cr
      \partial &\hbox{for } t \ge \tau_D(\omega) \cr
    \end{cases},
  $$
  where $\partial$ is a cemetery state. The subprocesses of other Markov processes in an open set can be defined in a similar way. Denote by $\sL^D$ the infinitesimal generator of $X^D$, which is the non-local operator $\sL$ in $D$ satisfying the zero exterior condition.

  Let $\{P_t; t\geq 0\}$ be the transition semigroup of the rectilinear $\alpha$-stable process $X$; that is, for $t>0$, $ x\in \R^d$ and $f\geq 0$ on $\R^d$,
  $$
    P_t f(x):= \E_x \left[ f(X_t) \right] = \int_{\R^d} p(t,x,y) f(y) dy .
  $$
  By \eqref{eq:pEst}, $\{P_t; t\geq 0\}$ is a strongly continuous semigroup in the Banach space $C_\infty (\R^d)$ of continuous functions that vanish at infinity equipped with the uniform norm $\|f\|_\infty:=\sup_{x\in \R^d} |f(x)|$. Moreover, since $p(t,x,y)$ is jointly continuous and has estimates \eqref{eq:pEst}, $\{P_t; t\geq 0\}$ has strong Feller property in the sense that for every $t>0$ and any bounded function $f$ on $\R^d$, $P_t f$ is a bounded continuous function on $\R^d$. Thus the L\'evy process $X$ is a Feller process having strong Feller property. By the proof of \cite[Theorem on p.68]{Chung.1986.63}, the semigroup $\{P^D_t; t\geq 0\}$ of $X^D$ has strong Feller property for any non-empty open subset $D\subset \mathbb{R}^d$. (Observe that the proof of the strong Feller property of  $\{P^D_t; t\geq 0\}$ in \cite{Chung.1986.63} does not need regular assumption on $D$.) In this paper, we will show that $X^{D}$ has a jointly locally H\"older continuous transition density $p_D(t,x,y)$ with respect to the Lebesgue measure. Furthermore, we will investigate the strict positivity property and the two-sided estimates of $p_D(t,x,y)$ for a class of open subsets $D\subset \R^d$.

 \medskip

  \begin{thm}\label{thm:pD-exist}
    For any non-empty open set $D\subset \R^d$, the subprocess $X^D$ has a jointly (locally) H\"older continuous density function $p_D(t,x,y)$ on $(0, \infty) \times D \times D$; that is, for any $x\in D$ and any non-negative Borel measurable function $\varphi$ on $D$,
    \begin{equation}\label{e:1.11}
      \E_x \left[ \varphi(X^D_t) \right] = \int_D p_D(t,x,y) \varphi(y) dy.
    \end{equation}
  \end{thm}

  Throughout this paper, we use the convention that any function $\varphi$ defined on $D$ is extended to $\partial$ by setting $\varphi(\partial) = 0$. We also call $p_D(t,x,y)$ the heat kernel of $X^D$ (or, equivalently,  of $\sL^D$), or the Dirichlet heat kernel of $X$ (or, equivalently,  of $\sL$) in $D$.

\medskip

  Unlike the rotationally symmetric $\alpha$-stable process $Z$, the behaviors of $X^D$ and $p_D(t,x,y)$  are strongly dependent on the shape of the domain $D$ due to the special structure of the L\'evy measure of the rectilinear $\alpha$-stable process $X$. For example, $X^D$ can be reducible for some smooth bounded open sets $D$.

  \begin{definition}
    We say a Markov process $\{Y,\mathbb{P}_x\}$ on  a topological state  space $E$ is irreducible if for any non-empty open subset $U \subset E$,
    \begin{equation*}
      \mathbb{P}_x( \sigma_U < \infty) > 0 \quad \hbox{for every } x  \in E,
    \end{equation*}
    where $\sigma_U:= \inf\{t>0: Y_t \in U\}$. Otherwise, we say the process $\{Y,\mathbb{P}_x\}$ is reducible.
  \end{definition}

\medskip

  The next result gives a geometric criterion on $D$ for the irreducibility of the rectilinear $\alpha$-stable subprocess $X^D$ in $D$.

  \begin{thm}\label{thm:irreducible}
    Let $D\subset \mathbb{R}^d$ be a non-empty open set. The subprocess $X^D$ is irreducible if and only if
    \begin{equation}  \label{e:1.13}
      \begin{split}
        &\emph{for every $x,y \in D$, there is $N\ge 1$ and some $\{x_i\}_{i=0}^N\subset D$ with $x_0 = x$ and $x_N = y$ so} \\
        &\emph{that each consecutive pair $(x_{i-1}, x_i)$, $1\leq i\leq N$, has only one different coordinate.}
      \end{split}
    \end{equation}
    Moreover, $X^D$ is irreducible if and only if  $p_D(t,x,y)>0$ for every $t>0$ and $x, y\in D$.
  \end{thm}

  Theorem \ref{thm:irreducible} together with Theorem  \ref{thm:pD-exist} in particular implies that for any connected open set $D$, $X^D$ is irreducible and has a strictly positive continuous transition density function $p_D(t,x,y)$.

  \begin{corollary}\label{C:1.4}
    Suppose that $D\subset \mathbb{R}^d$ is  a non-empty open set, and $D_1$ and $D_2$ are two  disjoint connected components of $D$. Then
    \begin{enumerate} [\rm (i)]
      \item $p_D(t,x,y) >0$ for every $t>0$ and $x, y\in D_1$.

      \item Either $p_D(t,x,y)>0$ for every $(t,x,y) \in  (0, \infty) \times D_1 \times D_2$ or $p_D(t,x,y)=0$ for every $(t,x,y) \in  (0, \infty) \times  D_1 \times D_2$. The former happens if and only if  there exists a finite sequence $\{x_i\}_{i=0}^N \subset D $ with $x_0 \in D_1$ and $x_N \in D_2$ so that each consecutive pair $(x_{i-1}, x_i)$, $1\leq i\leq N$,  has only one different coordinate.
    \end{enumerate}
  \end{corollary}

\smallskip

  See Theorems \ref{T:5.2}-\ref{T:5.3} for further information on the positivity of $p_D(t,x,y)$. To obtain more precise bounds on $p_D(t,x,y)$, we need certain smoothness of $D$ and some additional geometric condition on $D$ beyond \eqref{e:1.13} (or equivalently, the irreducibility of $X^D$).

 \medskip

  Recall that an open set $D \subset \mathbb{R}^d$ is said to be $C^{1,1}$ with characteristics $(R,\Lambda)$ for some $R,\Lambda>0$, if for every $z \in \partial D$, there is a $C^{1,1}$-function $\phi = \phi_z: \mathbb{R}^{d-1} \rightarrow \mathbb{R}$ satisfying $\phi(0) = 0$, $\nabla \phi(0) = 0$, $|\nabla\phi(\tilde{x}) - \nabla \phi(\tilde{y})| \le \Lambda |\tilde{x}-\tilde{y}|$, $\tilde{x},\tilde{y} \in \mathbb{R}^{d-1}$,  and an orthogonal coordinate system $CS_z: y = (y^{(1)},\cdots,y^{(d-1)},y^{(d)}) =:(\tilde{y},y^{(d)})$ with its origin at $z$ such that
  \begin{equation*}
    B(z,R) \cap D = \left\{y=(\tilde{y},y^{(d)}) \in B(0,R)\text{ in } CS_z: y^{(d)} > \phi(\tilde{y}) \right\}.
  \end{equation*}
  The pair $(R,\Lambda)$ is called the characteristics of the $C^{1,1}$ open set $D$. Note that the $C^{1,1}$ open set $D$ may be disconnected and may have infinite number of components. However, the distances between any two distinct connected components of $D$ are bounded from below uniformly by a positive constant.

\medskip

  For an open set $D \subset \mathbb{R}^d$ and $x\in D$, let $\delta_D(x)$ be the Euclidean distance between $x$ and $D^c$. We say $D$ satisfies the uniform interior ball condition with radius $R_1 > 0$, if, for every $x\in D$ with $\delta_D(x) \le R_1$, there is $z_x \in \partial D$ such that $|x-z_x| = \delta_D(x)$ and $B(x_0,R_1) \subset D$ for $x_0 := z_x + R_1(x-z_x)/|x-z_x|$; see \cite{ChenKimSong.2011.JLMS258,Cho.2006.PA387}. Similarly, we can define the uniform exterior ball condition.

\medskip

  It is well known that $D$ being a $C^{1,1}$ open set is equivalent to that $D$ satisfies both the uniform interior and exterior ball conditions. Thus without loss of generality, in this paper, for a $C^{1,1}$ open set $D$, we always assume its $C^{1,1}$ characteristics $(R,\Lambda)$ have the property that $R\le 1$, $\Lambda \ge 1$ and it satisfies the uniform interior and exterior ball conditions with radius $R$.

\medskip

  For $u=(u^{(1)}, u^{(2)},\cdots,u^{(d)}) \in \mathbb{R}^d$, $a\in\mathbb{R}$ and $1\le i\le d$, let
  \begin{equation*}
    [u]_a^i := (u^{(1)},\cdots,u^{(i-1)},a,u^{(i+1)},\cdots,u^{(d)});
  \end{equation*}
  that is, $[u]_a^i$ is the point in $\R^d$ by changing its $i^{th}$-coordinate to $a$. For $x,y \in \R^d$ and a permutation $\{i_1,i_2,\cdots,i_d\}$ of $\{1,2,\cdots,d\}$, let
  \begin{equation*}
    \overline{xy}_{1} := [x]_{y^{(i_1)}}^{i_1}, \  \   \overline{xy}_{2} := [\overline{xy}_{1}]_{y^{(i_2)}}^{i_2}, \  \ \overline{xy}_{3} := [\overline{xy}_{2}]_{y^{(i_3)}}^{i_3},  \   \  \cdots,  \    \     \overline{xy}_{d} := [\overline{xy}_{d-1}]_{y^{(i_d)}}^{i_d}=y.
  \end{equation*}
  That is, $\overline{xy}_{j}$ is the point obtained by swapping consecutively the $i_k^{th}$-coordinate of $x$ with that of $y$ for $k=1,2, \cdots, j$.

\medskip

  We consider the following geometric condition on an open set $D$. Let $\gamma \in (0,1]$.

\medskip

  \begin{enumerate}
    \item[$\bf (H_{\gamma})$:] \label{eq:conH}
        An open set  $D\subset \R^d$ is said to satisfy condition \hyperref[eq:conH]{$\bf (H_{\gamma})$} if for any $x,y \in D $ with $\delta_D(x)\wedge\delta_D (y) \ge r >0$, there exists a permutation $\{i_1,i_2,\cdots,i_d\}$ of $\{1,2,\cdots,d\}$ so that
        \begin{equation*}
          B(\overline{xy}_{k},\gamma r) \subset D,\qquad k=1,2,\cdots,d.
        \end{equation*}
  \end{enumerate}

\smallskip

\begin{figure}[htp]
  \begin{tikzpicture}[
      scale = 1.5
    ]
    \draw [->] (-2.5,0) -- (2.5,0);
    \draw [->] (0,-0.4) -- (0,2.2);
    \node [below right] {$0$};
    \draw (-1.3,1.1) circle (1);
    \fill (-1.3,1.1) circle (1pt);
    \draw [->] (-1.3,1.1) -- (-1.3+0.866,1.1-0.5)  node[pos =0.5,sloped,above] {$r=1$};
    \node [above left] at (-1.3, 1.1) {\Large $X_0$};
    \draw (1.3,1.1) circle (1);
    \fill (1.3,1.1) circle (1pt);
    \draw [->] (1.3,1.1) -- (1.3+0.866,1.1-0.5) node[pos =0.5,sloped,above] {$r=1$};
    \node [above left] at (1.3, 1.1) {\Large $Y_0$};
  \end{tikzpicture}\qquad
  \begin{tikzpicture}[
      scale = 1.5
    ]
    \draw [->] (-2,0) -- (2.5,0);
    \draw [->] (-1.7,-0.4) -- (-1.7,2.2);
    \draw [rounded corners = 10pt] (-1.3,0.2) -- (2.2,0.2) -- (2.2,2.2) -- (-1.3,2.2) -- cycle;
    \node at (1,1) {\Large$D$};
    \node [below right] at (-1.7,0) {$0$};
  \end{tikzpicture}
  \caption{The set $D := B(X_0,r) \cup B(Y_0,r)$ with $r = 1$, and the set $D:=$ the cubes with round corners in $\mathbb{R}^2$}
  \label{fig:parallel}
\end{figure}
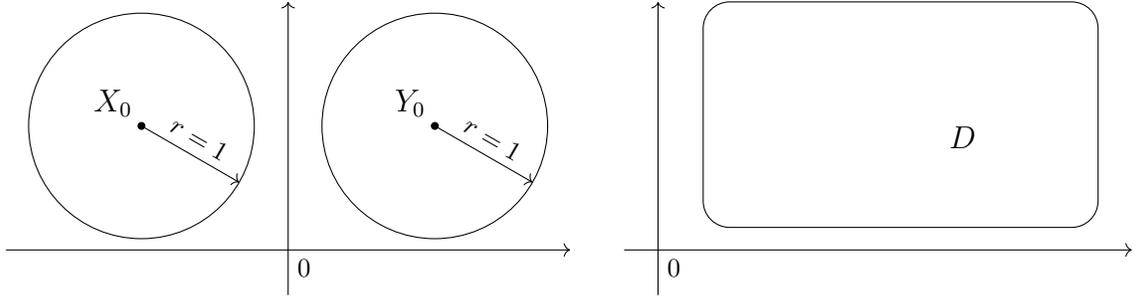

  Clearly, for any $0<\gamma_1<\gamma_2 \leq 1$, condition $\bf (H_{\gamma_2})$ implies $\bf (H_{\gamma_1})$, and any \hyperref[eq:conH]{$\bf (H_{\gamma})$} with $\gamma \in (0,1]$ implies the irreducibility condition \eqref{e:1.13}. Many open sets in $\mathbb{R}^d$ satisfy condition \hyperref[eq:conH]{$\bf (H_{\gamma})$}. For example, all balls, complements of closed balls, and the open sets shown in Figure \ref{fig:parallel} satisfy \hyperref[eq:conH]{$\bf (H_{1})$}. But there are also many open sets which do not satisfy condition \hyperref[eq:conH]{$\bf (H_{\gamma})$}; see Section \ref{S:6} for some examples.

\smallskip

  Recall that $p(t,x,y)$ is the transition density function (also called the heat kernel) of the rectilinear stable process $X$, and, for any open set $D\subset \R^d$, $p_D(t,x,y)$ is the heat kernel of $X^D$.  Recall also that $\sL^D$ is the infinitesimal generator of $X^D$.

\const{\CTmainUp}
\const{\CTmainLow}
\const{\CTeigen}
\const{\CTmainL}
  \begin{thm}\label{T:1.1}
    Let $D\subset \mathbb{R}^d$ be a $C^{1,1}$ open set with characteristics $(R,\Lambda)$.
    \begin{enumerate}
      \item[\rm (i)] For any $T>0$, there exists $C_{\CTmainUp} = C_{\CTmainUp}(d,\alpha,R,\Lambda,T) > 0$ such that for all $t \in (0,T]$ and $x,y \in D$,
          \begin{equation}\label{eq:pDEst-small-1}
            p_D(t,x,y) \le C_{\CTmainUp}\left(1\wedge \frac{\delta_D(x)^{\alpha/2}}{\sqrt{t}}\right)\left(1\wedge \frac{\delta_D(y)^{\alpha/2}}{\sqrt{t}}\right)p(t,x,y).
          \end{equation}

      \item[\rm (ii)] Assume in addition that $D$ satisfies \hyperref[eq:conH]{$\bf (H_{\gamma})$} for some $\gamma \in (0,1]$. Then, for any $T>0$, there exists $C_{\CTmainLow} = C_{\CTmainLow}(d,\alpha,R,\Lambda,\gamma,T) > 0$ such that for all $t \in (0,T]$ and $x,y \in D$,
          \begin{equation}\label{eq:pDEst-small-2}
            p_D(t,x,y) \ge C_{\CTmainLow}\left(1\wedge \frac{\delta_D(x)^{\alpha/2}}{\sqrt{t}}\right)\left(1\wedge \frac{\delta_D(y)^{\alpha/2}}{\sqrt{t}}\right)p(t,x,y).
          \end{equation}

      \item[\rm (iii)] Assume in addition that $D$ is bounded and satisfies \hyperref[eq:conH]{$\bf (H_{\gamma})$} for some $\gamma \in (0,1]$.  Denote by $\lambda_1(D)$ the first eigenvalue of $-\mathcal{L^D}$.     Then, there exists $C_{\CTeigen} = C_{\CTeigen}(d,\alpha,R,\Lambda,\gamma,\mathrm{diam}(D)) > 1$ such that
          \begin{equation}\label{eq:eigenBnd}
             C_{\CTeigen}^{-1}\le \lambda_1 (D)  \le C_{\CTeigen},
          \end{equation}
          and, for any $T>0$, there exists $C_{\CTmainL} = C_{\CTmainL}(d,\alpha,R,\Lambda,\gamma,T,\mathrm{diam}(D)) > 0$ such that for all $t \in[T,\infty)$ and $x,y \in D$,
          \begin{equation}\label{eq:pDEst-large}
            C_{\CTmainL}^{-1} e^{-\lambda_1 (D) t}\delta_D(x)^{\alpha/2}\delta_D(y)^{\alpha/2} \le p_D(t,x,y) \le C_{\CTmainL} e^{-\lambda_1 (D) t}\delta_D(x)^{\alpha/2}\delta_D(y)^{\alpha/2}.
          \end{equation}
    \end{enumerate}
  \end{thm}

\medskip

  For the comparable lower bound estimate \eqref{eq:pDEst-small-2} to hold, certain geometric condition beyond smoothness of the bounded open set $D$ is needed. We show by Examples  \ref{em:1} and \ref{em:2} in Section \ref{S:6} that there are smooth connected bounded domains, even some smooth convex domains, that do not satisfy condition \hyperref[eq:conH]{$\bf (H_{\gamma})$} for any $\gamma \in (0, 1]$ and the lower bound of Dirichlet heat kernel estimate \eqref{eq:pDEst-small-2} fails. In Example \ref{em:3}, a bounded smooth open set $D$ is given that does not satisfy the condition \eqref{e:1.13} and thus $X^D$ is not irreducible. These facts are in strong contrast with that for the rotationally symmetric stable processes in $\R^d$, whose subprocesses in open sets are always irreducible and the comparable two-sided Dirichlet heat kernel estimates, obtained in  \cite{ChenKimSong.2010.JEMSJ1307}, are known to hold for any  $C^{1,1}$ smooth open sets. We next present a bounded smooth open set $D$ that satisfies the irreducibility condition \eqref{e:1.13} but does not satisfy condition \hyperref[eq:conH]{$\bf (H_{\gamma})$} for any $\gamma \in (0,1]$ and the lower bound \eqref{eq:pDEst-small-2} fails for $p_D(t,x,y)$, nevertheless for which we can still derive comparable upper and lower bound of Dirichlet heat kernel estimates.

 \medskip

\begin{figure}[htp]
  \begin{tikzpicture}[
    declare function = {
        r = 2;
        },
      scale = 0.55
    ]
    \draw [->] (-3-r-0.5,0) -- (9+r+0.5,0);
    \draw [->] (0,-3-r-0.2) -- (0,3+r+0.2);
    \node [below right] {$0$};

    \draw [rounded corners = 10pt] (-3-r,-3-r) -- (-3+r,-3-r) -- (-3+r,-3+r) -- (-3-r,-3+r) -- cycle;
    \fill (-3,-3) circle (1.5pt);
    \draw [->] (-3,-3) -- (-3+r,-3)  node[pos =0.5,sloped,above] {$r=1$};
    \node [above left] at (-3, -3) {\Large $O_1$};
    \node [below] at (-3, -3.5) {\Large $A_1$};

    \draw [rounded corners = 10pt] (3-r,-3-r) -- (3+r,-3-r) -- (3+r,-3+r) -- (3-r,-3+r) -- cycle;
    \fill (3,-3) circle (1.5pt);
    \node [above left] at (3, -3) {\Large $O_2$};
    \node [below] at (3, -3.5) {\Large $A_2$};

    \draw [rounded corners = 10pt] (3-r,3-r) -- (3+r,3-r) -- (3+r,3+r) -- (3-r,3+r) -- cycle;
    \fill (3,3) circle (1.5pt);
    \node [above left] at (3, 3) {\Large $O_3$};
    \node [below] at (3, 2.5) {\Large $A_3$};

    \draw [rounded corners = 10pt] (9-r,3-r) -- (9+r,3-r) -- (9+r,3+r) -- (9-r,3+r) -- cycle;
    \fill (9,3) circle (1.5pt);
    \node [above left] at (9, 3) {\Large $O_4$};
    \node [below] at (9, 2.5) {\Large $A_4$};

     \end{tikzpicture}
  \caption{The set $D$ is the union of four squares with round corners in $\mathbb{R}^2$}
  \label{fig:4squares}
\end{figure}
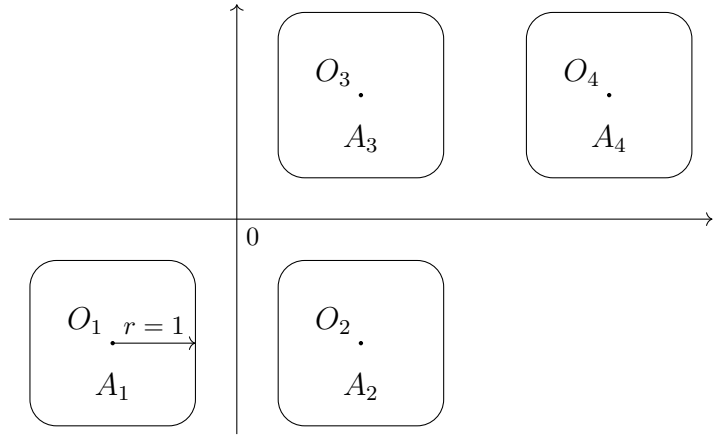

  Let $O_i\in \mathbb{R}^2$, $i=1,\cdots,4$ be four points such that, the line through $O_1$ and $O_2$ is paralleled with $x$-axis, the line through $O_2$ and $O_3$ is paralleled with $y$-axis, the line through $O_3$ and $O_4$ is paralleled with $x$-axis, and $|O_i - O_{i+1}| = 3$ for $i = 1,2,3$; see Figure \ref{fig:4squares}. Let $A_i$, $i=1,\cdots,4$ be four squares with round corners and with edge-length 2 centered at $O_i$ respectively. Consider the open set
  \begin{equation*}
    D := \bigcup_{i=1}^4 A_i.
  \end{equation*}

  Note that $D\subset \mathbb{R}^2$ is a bounded smooth open set that satisfies the irreducibility condition \eqref{e:1.13} but does not satisfy condition \hyperref[eq:conH]{$\bf (H_{\gamma})$} for any $\gamma \in (0,1]$ as for any $x \in A_1$ and $y\in A_4$, and swapping any coordinate of $x$ by that of $y$ results a point falling outside $D$.

%\medskip

  \const{\CTexmDa}
  \const{\CTexmDb}
  \const{\CTexmDc}
  \begin{thm} \label{T:1.6}
    Let $D\subset \mathbb{R}^2$ be the above smooth open set as shown in Figure \ref{fig:4squares} and $T >0$.
    \begin{enumerate}
      \item[\rm (i)] There exists $C_{\CTexmDa} = C_{\CTexmDa}(\alpha,T) > 0$ such that for all $t \in (0,T]$, $x \in A_i$ and  $y \in A_j$ with $|i-j|\le 2$,
          \begin{equation}\label{eq:pDEst-exm-square-1}
            p_D(t,x,y) \stackrel{C_{\CTexmDa}}{\asymp} \left(1\wedge \frac{\delta_D(x)^{\alpha/2}}{\sqrt{t}}\right) \left(1\wedge\frac{\delta_D(y)^{\alpha/2}}{\sqrt{t}}\right) p(t,x,y).
          \end{equation}

      \item[\rm (ii)] There exist  $C_{\CTexmDb} = C_{\CTexmDb}(\alpha,T) > 0$ and $C_{\CTexmDc}=C_{\CTexmDc}(\alpha,T)>0$ such that for all $t \in (0,T]$,  $x \in A_1$ and $y \in A_4$,
          \setlength{\arraycolsep}{1.5pt}
          \begin{eqnarray}
            p_D(t,x,y) & \stackrel{C_{\CTexmDb}}{\asymp}& t^3 \left(1\wedge \frac{\delta_D(x)^{\alpha/2}}{\sqrt{t}}\right) \left(1\wedge\frac{\delta_D(y)^{\alpha/2}}{\sqrt{t}}\right) \label{eq:pDEst-exm-square-2}\\
            &\stackrel{C_{\CTexmDc}}{\asymp}& t \left(1\wedge\frac{\delta_D(x)^{\alpha/2}}{\sqrt{t}}\right) \left(1\wedge \frac{\delta_D(y)^{\alpha/2}}{\sqrt{t}}\right) p(t,x,y).  \label{e:6.12}
        \end{eqnarray}
    \end{enumerate}
  \end{thm}

 \smallskip

  Here and in the sequel, for two functions $f,~g$ and a positive constant $C$, the notation $f\stackrel{C}{\asymp}g$ means that $C^{-1} f \le g\le Cf$ holds true on their common domains.  Theorem \ref{T:1.6}(ii) shows that, for the smooth open set $D$ in Figure \ref{fig:4squares}, the lower bound \eqref{eq:pDEst-small-2} fails for $p_D(t,x,y)$.

 \medskip

  For an open set $D\subset \R^d$,  we call
  $$
    G_D(x, y) := \int_0^\infty p_D(t,x,y) dt
  $$
  the Green function of $X$ in $D$. It follows from \eqref{e:1.11} that  for any $x\in D$ and any non-negative Borel measurable function $\varphi$ on $D$,
  $$
   \E_x \int_0^{\tau_D} \varphi(X_s) ds = \int_D G_D(x, y) \varphi(y) dy.
  $$
  From the Dirichlet heat kernel estimates in Theorem \ref{T:1.1} for $p_D(t,x,y)$, one can clearly derive the Green function estimates for $G_D(x,y)$.

\medskip

    Finally, we mention that  boundary regularity of solutions to the Dirichlet problem for the generator of isotropic stable processes
 is   studied in \cite{Ros-OtonSerra.2014.JMPA}.
  These regularity results  were later extended to more general stable operators in \cite{Ros-OtonSerra.2016.JDE, Ros-OtonValdinoci.2016.Adv}, and further to   nonlinear nonlocal equations in \cite{Ros-OtonSerra.2016.Duke}.
  However, for Dirichlet heat kernel estimates, the singular nature of rectilinear stable processes poses significant challenges.
   For instance, as mentioned above, Harnack inequality fails for rectilinear stable processes. Thus we can not use the approach developed in \cite{ChenKimSong.2010.JEMSJ1307} for the study of Dirichlet heat kernel estimates for rotationally symmetric stable processes directly. To see this through a concrete case, we invite the interested reader to pause for a few minutes and think about possible ways to establish the joint local H\"older continuity of the transition density $p_D(t,x,y)$ of $X^D$  in any open subset $D\subset \R^d$ before reading the proof of Theorem \ref{thm:qHolder}. Some new ideas and methods are needed for the study of rectilinear stable processes. We employ a combination of the probabilistic and analytic methods in our investigation.

\medskip

  The rest of the paper is organized as follows. In Section \ref{S:2}, we show that the part process $X^D$ in any open subset $D\subset \R^d$ has a locally H\"older continuous transition density function. Boundary properties of the harmonic measures of $\sL$, or equivalently, the exit distributions of $X$, in $C^{1,1}$ open sets are investigated  in Section \ref{S:3}, using testing function methods developed in \cite{BogdanBurdzyChen.2003.PTRF89,  ChenKimSongEtAl.2012.TAMS4169}. Various Dirichlet heat kernel estimates are obtained  in Section \ref{S:4},  and  the proof  of Theorem \ref{T:1.1} is given in Subsection \ref{Subs:4.3}. For the upper bound estimates of $p_D(t,x,y)$, we use the exit time estimates, strong Markov property and the L\'evy system of the rectilinear stable process $X$. For the lower bound estimates of  $p_D(t,x,y)$, we first obtain its near diagonal interior estimate in Lemma \ref{lem:pUnle-int}, and then the interior estimates under the condition \hyperref[eq:conH]{$\bf (H_{\gamma})$} for some $\gamma \in (0,1]$ in Lemma \ref{lem:pUle-int} using the Chapman-Kolmogorov equation, a chaining argument and a delicate probability lower bound estimate for $X^U_t$ taking values in suitable cubes. The sharp lower bound estimates for $p_D(t,x,y)$ over some bounded time interval $(0,t_*]$ in any $C^{1,1}$ open set $D$ satisfying the condition \hyperref[eq:conH]{$\bf (H_{\gamma})$} for some $\gamma \in (0,1]$ is established in Lemmas \ref{lem:pDTle}-\ref{lem:pDtle} though a careful probabilistic argument that boils down to the exit time estimates for $X$. The proof of Theorem \ref{T:1.1} is given Subsection \ref{Subs:4.3}, where in particular the lower bound estimate in Lemma \ref{lem:pDtle} over {\it some} bounded time interval is shown to hold over {\it any} bounded time interval through a chaining argument. For any two fixed distinct points $x, y\in D$, a geometric criteria for the positivity of $p_D(t,x,y)$ is given in Theorems \ref{T:5.2} and \ref{T:5.3}, whose proof uses some of the lower bound estimates derived in Section \ref{S:4}. From these, the proof of Theorem \ref{thm:irreducible} is given in Section \ref{S:irred}. In addition to the proof of Theorem \ref{T:1.6}, three additional examples of bounded smooth open sets are given in Section \ref{S:6}, two of them  are connected open sets,  that do not satisfy the condition \hyperref[eq:conH]{$\bf (H_{\gamma})$} for any $\gamma \in (0,1]$, for which the lower bound estimate \eqref{eq:pDEst-small-2} is shown to fail.

\smallskip

  In this paper, for $x=(x^{(1)},x^{(2)},\cdots,x^{(d)}) \in \mathbb{R}^d$ and $r > 0$,  we will use  $Q(x,r)$ to denote the cube centered at $x$ with edge-length $2r$, that is,
  \begin{equation*}
    Q(x,r) := \left\{y=(y^{(1)},y^{(2)},\cdots,y^{(d)}) \in \mathbb{R}^d: |x^{(i)} - y^{(i)}| < r,~i=1,2,\cdots,d \right\}.
  \end{equation*}
  For an open set $U\subset \R^d$ and $\lambda >0$, unless otherwise stated, we define
  \begin{equation*}
    \lambda U := \{\lambda y: y \in U\}.
  \end{equation*}
  For a measurable set $A\subset\mathbb{R}^d$, we use $|A|$ to denote its Lebesgue measure.

\medskip

\section{H\"older regularity of Dirichlet heat kernel}\label{S:2}

In this section, we fix a non-empty open set $D\subset \mathbb{R}^d$. Recall that $\tau_D(\omega) = \inf\{t>0,X_t(\omega) \notin D\}$ is the first exit time from $D$ by the rectilinear $\alpha$-stable process $X$. Since $X_t$ has a continuous transition density function with respect to the Lebesgue measure, we have the following property by the same proof as that in \cite[Proposition 1.20]{ChungZhao.1995.287}.

\begin{prop}\label{P:2.1}
  For every $t>0$ and $x\in \R^d$, $\P_x (\tau_D =t)=0$.
\end{prop}

For $t > 0$, $x,y\in \mathbb{R}^d$, set
\begin{equation*}
  r_D(t,x,y) := \mathbb{E}_x\left[p(t-\tau_D,X_{\tau_D},y);\tau_D < t\right],
\end{equation*}
and
\begin{equation}\label{eq:pU-def}
  p_D(t,x,y) := p(t,x,y) - r_D(t,x,y).
\end{equation}
Note that the function $p_D(t,x,y)$ is pointwise well-defined.  One can follow the proof of \cite[Theorem 2.4, p. 33]{ChungZhao.1995.287} to prove the following lemma. Let $\mathcal{B}(\mathbb{R}^d)$ denote the collection of all Borel measurable sets in $\mathbb{R}^d$.

\begin{lemma}\label{lem:dhk-def-1}
  For any $t > 0$, $x\in\mathbb{R}^d$ and $A \in \mathcal{B}(\mathbb{R}^d)$,
  \begin{equation}\label{eq:pUprobability}
    \mathbb{P}_x\left(X_t \in A, \  t<\tau_D\right) = \int_A p_D(t,x,y) dy.
  \end{equation}
  The function $p_D(t,x,y)$ is almost surely symmetric on $\mathbb{R}^d\times \mathbb{R}^d$: for all $t > 0$,
  \begin{equation}\label{eq:pUSymm-a.e.}
    p_D(t,x,y) = p_D(t,y,x) \quad \text{for a.a. } (x,y) \in D\times D.
  \end{equation}
  Moreover, for any $s,t > 0$ and $x\in \mathbb{R}^d$, we have
  \begin{equation}\label{eq:pUC-K-a.e.}
    p_D(t+s,x,y) = \int_{\mathbb{R}^d} p_D(t,x,z) p_D(s,z,y) dz \quad \text{for a.a. } y \in D.
  \end{equation}
\end{lemma}

Unfortunately, we can not use the approach in \cite[Theorem 2.4 on p.33]{ChungZhao.1995.287} to establish the joint continuity of $p_D(t,\cdot,\cdot)$ on $(\R^d \setminus \partial D) \times (\R^d \setminus \partial D)$, and hence to improve the identities in \eqref{eq:pUSymm-a.e.} and \eqref{eq:pUC-K-a.e.} from almost every point to every point. The main issue is that unlike Brownian motion or rotationally symmetric stable processes case, in our setting, the function $\1_{\{\tau_D < t\}} p(t-\tau_D,X_{\tau_D},y)$ is unbounded. However, we can apply the result in \cite{BassChen.2010.MZ489} and the ideas in \cite[Proposition 2.5, p.1603]{ChenCroydonKumagai.2015.AP1594} to establish the joint  H\"older continuity of $p_D(t, \cdot, \cdot)$ on $D\times D$ in Theorem \ref{thm:qHolder} (see also \cite{BendikovGrigoryanHuEtAl.2021.ASNSPCS5399} or \cite{GrigoryanHuHu.2018.AM433} for another approach). In order to make the proof as self-contained as possible, we show all the details.

\smallskip

The rectilinear $\alpha$-stable process $X$ has the following scaling property: for any $\lambda > 0$, the processes $\{\lambda  X_{\lambda^{-\alpha}t} ; t\geq 0\}$ conditioned on  $X_0= x$ has the same distribution as $\{X_t; t\geq 0\}$ conditioned on $X_0=\lambda x$. Consequently, since the heat kernel $p(t,x,y)$ is continuous, it has the following scaling property: for any $\lambda>0$,
\begin{equation}\label{eq:HKscale}
  p(t,x,y) = \lambda^{-d} p(\lambda^{-\alpha}t,\lambda^{-1}x,\lambda^{-1}y),\qquad t> 0, ~ x,y \in \mathbb{R}^d.
\end{equation}
Moreover, $\{\lambda X_{\lambda^{-\alpha}t}^D; t\geq 0\}$ conditioned on $X_0=x\in D$ has the same distribution as $\{X_t^{\lambda D}; t\geq 0\}$ conditioned on $X_0=\lambda x$. It follows that for every $\lambda > 0$, $t>0$ and $x\in D$,
\begin{equation}\label{eq:DHKscale}
  p_D(t,x,y) = \lambda^d p_{\lambda D}(\lambda^\alpha t,\lambda x,\lambda y)  \quad \hbox{for a.a. }  y \in D.
\end{equation}

The following lemma gives the exit time estimates for $X$ from balls.

\begin{lemma}\label{lem:meanExit}
  There exist positive constants $c_i:=c_i(d,\alpha)>0$, $i=1,2$, such that for any $x_0\in \R^d$ and $r>0$,
  \begin{enumerate}[{\rm (i)}]
    \item $\mathbb{E}_x\left[\tau_{B(x_0,r)} \right] \le c_1 r^\alpha$ for all $x \in B(x_0,r)$;
    \item $\mathbb{E}_x\left[\tau_{B(x_0,r)}\right] \ge c_2 r^\alpha$ for all $x \in B(x_0,r/2)$.
  \end{enumerate}
\end{lemma}

\begin{proof}
  By the L\'evy property of $X$ we may assume that $x_0=0$. When $r = 1$, this lemma follows directly from \cite[Proposition 2.1, p. 492]{BassChen.2010.MZ489} with the matrix $A$ being the $d\times d$  identity matrix. For general $r>0$, the desired property follows from the scaling property of $X^{B(0,r)}$ that $\mathbb{E}_x \left[\tau_{B(0,r)} \right]  = r^\alpha \mathbb{E}_{x/r} \left[\tau_{B(0,1)} \right] $.
\end{proof}

\begin{definition}
  A bounded function $h$ on $\mathbb{R}^d$ is said to be harmonic (with respect to $X$) in a ball $B\subset \mathbb{R}^d$ if
  \begin{equation*}
    h(x) = \mathbb{E}_x\left[h(X_{\tau_B})\right]\quad \text{ for all } x \in B.
  \end{equation*}
\end{definition}

For a non-negative function $f$ on $D$, let
\begin{equation}\label{e:2.6}
  P_t^Df(x) := \int_D p_D(t,x,z) f(z) dz = \mathbb{E}_x  [ f(X_t^D)]   ,\quad t > 0,~ x \in \mathbb{R}^d,
\end{equation}
where the second equality is due to  \eqref{eq:pUprobability}. Thus $\{P_t^D, t> 0\}$ is the transition semigroup of $X^D$. It follows from \eqref{eq:pUSymm-a.e.} that $\{P_t^D, t> 0\}$ is a strongly continuous symmetric contractive semigroup on $L^2(D; dx)$. Moreover, by \eqref{e:2.6} and the Markov property of $X^D$, we have for any $s,t>0$
\begin{equation}\label{eq:smgp-pointwise}
  P_{t+s}^Df(x) = P_t^DP_s^Df(x),\quad x \in \mathbb{R}^d \setminus \partial D.
\end{equation}
We further define the $1$-potential
\begin{equation*}
  G^D_1f(x) := \mathbb{E}_x\left[\int_0^\infty e^{-t}f(X_t^D) dt\right] = \int_0^\infty e^{-t} P_t^D f(x) dt,\quad x \in \mathbb{R}^d \setminus \partial D.
\end{equation*}
Let $\{\theta_t; t\geq 0\}$ be the time-shifting operators on the canonical probability sample space $\Omega$ for the L\'evy process $X$; that is, $X_r(\theta_t \omega) = X_{r+t}(\omega)$ for $\omega \in \Omega$ and $t, r\geq 0$.

\begin{thm}\label{thm:qHolder}
  There is a jointly (locally) H\"older continuous function $q(t,x,y)$ on $D\times D$ so that
  \begin{enumerate}[{\rm (i)}]
    \item for any $t > 0$ and $x \in D$,
        \begin{equation*}
           q(t,x,z) =p_D(t,x,z)   \quad \text{for a.a. } z \in D;
        \end{equation*}
    \item for any $t > 0$, $q(t,x,y)$ is symmetric on $D\times D$;
    \item for any $s,t>0$ and $x,y \in D$,
        \begin{equation*}
          q(t+s,x,y) = \int_D q(t,x,z) q(s,z,y) dz.
        \end{equation*}
  \end{enumerate}
\end{thm}

\begin{proof}
  The proof uses a result from \cite{BassChen.2010.MZ489} and the ideas from \cite[Proposition 2.5, p.1603]{ChenCroydonKumagai.2015.AP1594}. For the reader's convenience, we spell out all the details. We divide the proof into five steps.

\medskip

  \textbf{Step 1.} Let $r \in (0,1]$ and $B:=B(x_0,r) \subset D$. Suppose that $h$ is a bounded function on $\mathbb{R}^d$ and is harmonic with respect to $X$ in $B$. By the H\"older regularity obtained in \cite[Theorem 2.9, p. 499]{BassChen.2010.MZ489}, taking the matrix $A$ there the $d\times d$ identity matrix, there exist positive constants $c_1$ and $\beta$ depending only on $d$ and $\alpha$ such that
  \begin{equation}\label{eq:H-Holder}
    |h(x) - h(y)| \le c_1\left(\frac{|x-y|}{r}\right)^\beta \sup_{z\in\mathbb{R}^d} |h(z)| \ \ \text{ for all } x,y \in B(x_0,r/2).
  \end{equation}

\medskip

  \textbf{Step 2.} Let $f \in L^\infty (D) \cap L^2(D)$. By the strong Markov property, we obtain that for any $x \in B$,
  \begin{align*}
    G^D_1f(x) =&~ \mathbb{E}_x\left[\int_0^{\tau_B} e^{-t}f(X_t^D) dt\right] + \mathbb{E}_x\left[\int_{\tau_B}^\infty e^{-t}f(X_t^D) dt\right]\\
    =&~ \mathbb{E}_x\left[\int_0^{\tau_B} e^{-t}f(X_t^D) dt\right] + \mathbb{E}_x\left[e^{-\tau_B}  \left(\int_0^\infty e^{-t}f(X_t^D) dt \right)  \circ \theta_{\tau_B}  \right] \\
    =&~ \mathbb{E}_x\left[\int_0^{\tau_B} e^{-t}f(X_t^D) dt\right] + \mathbb{E}_x\left[e^{-\tau_B}G^D_1f(X^D_{\tau_B})\right]\\
    =&~ \mathbb{E}_x\left[\int_0^{\tau_B} e^{-t}f(X_t^D) dt\right] + \mathbb{E}_x\left[(e^{-\tau_B}-1)G^D_1f(X^D_{\tau_B})\right] +\mathbb{E}_x\left[G^D_1f(X^D_{\tau_B})\right]\\
    =&~ I_1(x) + I_2(x) + I_3(x).
  \end{align*}
   By Lemma \ref{lem:meanExit}(i) and the elementary inequality that  $1-e^{-a} \leq a$ for $a\geq 0$, we have
  \begin{equation*}
    |I_1(x)| \le \|f\|_{L^\infty(D)}\mathbb{E}_x\left[\tau_B\right] \le c_2 r^\alpha \|f\|_{L^\infty(D)},
  \end{equation*}
  and
  \begin{equation*}
    |I_2(x)| \le \|G^D_1f\|_{L^\infty(D)} \mathbb{E}_x\left[\tau_B\right] \le c_2 r^\alpha \|f\|_{L^\infty(D)}.
  \end{equation*}
  Since $z\mapsto I_3(z) = \mathbb{E}_z\left[G^D_1f(X^D_{\tau_B})\right]$ is bounded and harmonic in $B$, we have by (\ref{eq:H-Holder}) that for any $x,y \in B(x_0,r/2)\subset D$,
  \begin{equation*}
    |I_3(x) - I_3(y)| \le c_1\left(\frac{|x-y|}{r}\right)^\beta \sup_{z\in\mathbb{R}^d}|I_3(z)| \le c_1\left(\frac{|x-y|}{r}\right)^\beta \|G^D_1f\|_{L^\infty(D)}\le c_1\left(\frac{|x-y|}{r}\right)^\beta \|f\|_{L^\infty(D)}.
  \end{equation*}
  Combining the above four formulas, we obtain that for all $x,y \in B(x_0,r/2) \subset  D$,
  \begin{equation}\label{eq:potential-Holder}
    |G^D_1f(x) - G^D_1f(y)| \le (4c_2+c_1) \left(r^\alpha+\frac{|x-y|^\beta}{r^\beta}\right) \|f\|_{L^\infty(D)}.
  \end{equation}

\medskip

  \textbf{Step 3.} Recall that $\mathcal{L}^D$ is the generator of the heat semigroup $\{P_t^D, t> 0\}$ on $L^2(D)$. Note that $\mathcal{L}^D$ is {\em negative} definite self-adjoint. By general theory of heat semigroup, we have that for any $s,s' >0$ and $f \in L^2(D)$,
  \begin{equation}\label{eq:smgp-interchange}
    P_s^D \mathcal{L}^D P_{s'}^Df = P_{s'}^D\mathcal{L}^D P_s^Df   \quad \text{ a.e.}.
  \end{equation}
  For a fixed  $f \in L^\infty (D) \cap L^2(D)$, set
  \begin{equation*}
    h_t := P_t^Df - \mathcal{L}^DP_t^D f \in L^2(D),\quad t>0.
  \end{equation*}
  By spectral theory, there exists a spectral family $\{E_\lambda,\lambda \in \mathbb{R}\}$ such that
  \begin{equation}\label{eq:spectralRepre}
    \sL^D=-\int_0^\infty \lambda dE_\lambda,  \quad  f = \int_0^\infty dE_\lambda f
    \quad \hbox{and} \quad
    P_t^Df  \stackrel{\text{a.e.}}{=} \int_0^\infty e^{-\lambda t} dE_\lambda f.
  \end{equation}
  Consequently,
  \begin{equation}\label{e:2.13}
    (1- \sL^D) G^D_1 f \stackrel{\text{a.e.}}{=} f \quad \text{ and } \quad h_t \stackrel{\text{a.e.}}{=}  \int_0^\infty (1+\lambda)e^{-\lambda t} dE_\lambda f.
  \end{equation}
  For any $g \in L^1(D)$, by (\ref{eq:pEst}), we have for any $t > 0$
  \begin{equation*}
    \|P_t^Dg\|_{L^\infty(D)} = \left\|\int_D p_D(t,\cdot,z)g(z)dz\right\|_{L^\infty(D)} \le \left\|\int_D |p(t,\cdot,z)||g(z)|dz\right\|_{L^\infty(D)} \le C_{\CTstable} t^{-d/\alpha} \|g\|_{L^1(D)},
  \end{equation*}
  and then,
  \begin{equation}\label{eq:PtUgL1->L2}
    \|P_t^Dg\|_{L^2(D)} \le \sqrt{\|P_t^Dg\|_{L^1(D)}\|P_t^Dg\|_{L^\infty(D)}} \le \sqrt{\|g\|_{L^1(D)}C_{\CTstable} t^{-d/\alpha} \|g\|_{L^1(D)}} = \sqrt{C_{\CTstable}} t^{-\frac{d}{2\alpha}} \|g\|_{L^1(D)}.
  \end{equation}
  Using the above inequality, Cauchy-Schwarz and the facts that
  \begin{equation*}
    \sup_{\lambda >0} (1+\lambda)e^{-\lambda t} \le (t\wedge 1)^{-1} < \infty \quad \text{ and }\quad \sup_{\lambda >0} (1+\lambda)e^{-\lambda t/2} \le 2(t\wedge 1)^{-1} < \infty,
  \end{equation*}
  we obtain
  \begin{align*}
    (h_t,g)_{L^2(D)} =&~ \int_0^\infty (1+\lambda)e^{-\lambda t} d(E_\lambda f,g)_{L^2(D)}\\
    \le&~ \left(\int_0^\infty (1+\lambda)e^{-\lambda t} d(E_\lambda f,f)_{L^2(D)}\right)^{1/2} \left(\int_0^\infty (1+\lambda)e^{-\lambda t} d(E_\lambda g,g)_{L^2(D)}\right)^{1/2}\\
    \le&~ 2(t\wedge 1)^{-1} \left(\int_0^\infty d(E_\lambda f,f)_{L^2(D)}\right)^{1/2} \left(\int_0^\infty e^{-\lambda t/2} d(E_\lambda g,g)_{L^2(D)}\right)^{1/2}\\
    =&~ 2(t\wedge 1)^{-1} \|f\|_{L^2(D)}\|P_{t/2}^Dg\|_{L^2(D)}\\
    \le&~ 2\sqrt{C_{\CTstable}} (t\wedge 1)^{-1}t^{-\frac{d}{2\alpha}} \|f\|_{L^2(D)}\|g\|_{L^1(D)}.
  \end{align*}
  Since $g \in L^1(D)$ is arbitrary, we obtain
  \begin{equation*}
    \|h_t\|_{L^\infty} \le 2\sqrt{C_{\CTstable}} (t\wedge 1)^{-1}t^{-\frac{d}{2\alpha}} \|f\|_{L^2(D)}.
  \end{equation*}
  On the other hand, we have by \eqref{eq:smgp-interchange} and \eqref{e:2.13} that a.e. on $D$,
  $$
    G^D_1 h_t = G^D_1 P^D_t f - G^D_1 \sL^D P^D_t f = P^D_t \left(G^D_1  f -\sL^D G^D_1 f\right) =P^D_t f.
  $$
  As noted earlier, $P^D_t f $ and $P_t^D h_t$ are continuous functions on $D$ by the strong Feller property of $X^D$. By the dominated convergence theorem, $G^D_1 h_t(x)= \int_0^\infty e^{-s}P^D_sh_t(x)ds$  is a bounded continuous function on $D$. Hence we have
  $$
   P^D_t f(x) = G^D_1 h_t (x)   \quad \hbox{for every } x\in D.
  $$
  This together with (\ref{eq:potential-Holder})  yields that for all $t >0$ and $x,y \in B(x_0,r/2)\subset D$,
  \begin{equation}\label{eq:PtU-OSC-1}
  \begin{split}
    |P_t^Df(x) - P_t^Df(y)| =&~|G^D_1 h_t (x) - G^D_1  h_t  (y)| \\
    \le&~ (4c_2+c_1) \left(r^\alpha+\frac{|x-y|^\beta}{r^\beta}\right) \| h_t \|_{L^\infty(D)}\\
    \le&~ 2\sqrt{C_{\CTstable}}(4c_2+c_1) \left(r^\alpha+\frac{|x-y|^\beta}{r^\beta}\right) (t\wedge 1)^{-1}t^{-\frac{d}{2\alpha}}\|f\|_{L^2(D)}.
  \end{split}
  \end{equation}

\medskip

  \textbf{Step 4.} For any fixed compact set $K \subset D$ and $x,y \in K$, let $\delta_K = \frac{1}{4}(\mathrm{dist}(K,\partial D) \wedge \mathrm{dist}(K,\partial D)^2 \wedge 1)$ and $x_0 = x$.

  Case 1: $|x-y| < \delta_K$. Setting $r := |x-y|^{1/2}$, we have
  \begin{equation*}
    \mathrm{dist}(K,\partial D) > \sqrt{\delta_K} > |x-y|^{1/2} = r > 2|x-y|.
  \end{equation*}
  Applying (\ref{eq:PtU-OSC-1}) for this $r$ and $x_0 = x$, we have
  \begin{align*}
    |P_t^Df(x) - P_t^Df(y)| \le&~ 2\sqrt{C_{\CTstable}}(4c_2+c_1) \left(r^\alpha+\frac{|x-y|^\beta}{r^\beta}\right) (t\wedge 1)^{-1}t^{-\frac{d}{2\alpha}}\|f\|_{L^2(D)}\\
    \le&~ 4\sqrt{C_{\CTstable}}(4c_2+c_1) |x-y|^{(\alpha\wedge\beta)/2} (t\wedge 1)^{-1}t^{-\frac{d}{2\alpha}}\|f\|_{L^2(D)}.
  \end{align*}

  Case 2: $|x-y| \ge \delta_K$. By the definition of $P_t^Df$ and Cauchy-Schwarz inequality, we have for all $t > 0$ and $z \in D$,
  \begin{equation*}
    |P_t^Df(z)| \le \sqrt{\int_D p_D(t,z,y)^2 dy} \, \|f\|_{L^2(D)} \le \sqrt{p(2t,z,z)}\|f\|_{L^2(D)}\le \sqrt{C_{\CTstable}} (2t)^{-d/(2\alpha)}\|f\|_{L^2(D)}.
  \end{equation*}
  This implies that
  \begin{align*}
    |P_t^Df(x) - P_t^Df(y)| \le &~ |P_t^Df(x)| + |P_t^Df(y)| \le 2\sqrt{C_{\CTstable}} t^{-\frac{d}{2\alpha}}\|f\|_{L^2(D)}\\
    \le&~ 2\sqrt{C_{\CTstable}} \left(\frac{|x-y|}{\delta_K}\right)^{(\alpha\wedge\beta)/2} t^{-\frac{d}{2\alpha}} \|f\|_{L^2(D)}.
  \end{align*}
  Combining the above two cases, we obtain for any $x, y\in K$,
  \begin{equation}\label{eq:PtU-OSC-2}
    |P_t^Df(x) - P_t^Df(y)| \le c_3(d,\alpha,\mathrm{dist}(K,\partial D)) |x-y|^{(\alpha\wedge\beta)/2} (t\wedge 1)^{-1}t^{-\frac{d}{2\alpha}}\|f\|_{L^2(D)}.
  \end{equation}

\medskip

  \textbf{Step 5.} For any $0<s<t$, we define
  \begin{equation*}
    q_s(t,x,y) : = \int_D p_D(t-s,x,z) p_D(s,y,z)dz = P_{t-s}^Dp_D(s,y,\cdot)(x)= P_s^Dp_D(t-s,x,\cdot)(y),\quad x,y \in D.
  \end{equation*}
  By (\ref{eq:pUSymm-a.e.}) and (\ref{eq:pUC-K-a.e.}), we obtain for any $0<s<t$ and $x \in D$,
  \begin{equation}\label{eq:qs=pU-a.e.}
    q_s(t,x,y) = p_D(t,x,y),\quad \text{a.a. } y \in D.
  \end{equation}

  On the other hand, let $K$ be any compact subset of $D$ as in Step 4, $t > 0$ and $x,y,x',y' \in K$. Replacing $t$ by $t-s$ and $f$ by $p_D(s,y,\cdot)$ in (\ref{eq:PtU-OSC-2}) and using (\ref{eq:pEst}), we obtain
  \begin{align*}
    |q_s(t,x,y)-q_s(t,x',y)| =&~ |P_{t-s}^Df(x) - P_{t-s}^Df(x')|\\
    \le&~ c_3(d,\alpha,\mathrm{dist}(K,\partial D)) |x-x'|^{(\alpha\wedge\beta)/2} ((t-s)\wedge 1)^{-1}(t-s)^{-\frac{d}{2\alpha}} \|p_D(s,y,\cdot)\|_{L^2(D)}\\
    \le&~ c_3 |x-x'|^{(\alpha\wedge\beta)/2} ((t-s)\wedge 1)^{-1}(t-s)^{-\frac{d}{2\alpha}} \sqrt{p(2s,y,y)}\\
    \le&~ c_3\sqrt{C_{\CTstable}} |x-x'|^{(\alpha\wedge\beta)/2} ((t-s)\wedge 1)^{-1}(t-s)^{-\frac{d}{2\alpha}} (2s)^{-\frac{d}{2\alpha}}.
  \end{align*}
  Similarly, replacing $t$ by $s$ and $f$ by $p_D(t-s,x',\cdot)$ in (\ref{eq:PtU-OSC-2}) and using (\ref{eq:pEst}), we obtain
  \begin{align*}
    |q_s(t,x',y)-q_s(t,x',y')| =&~ |P_{t-s}^Df(y) - P_{t-s}^Df(y')|\\
    \le&~ c_3\sqrt{C_{\CTstable}} |y-y'|^{(\alpha\wedge\beta)/2} (s\wedge 1)^{-1}s^{-\frac{d}{2\alpha}} (2(t-s))^{-\frac{d}{2\alpha}}.
  \end{align*}
  Adding up the above two inequalities, we have that $q_s(t,x,y)$ is jointly H\"older continuous on $K\times K$, and hence, is jointly (locally) H\"older continuous on $D\times D$. Moreover, the (locally) H\"older continuity of $q_s(t,x,y)$ on $D\times D$ and (\ref{eq:qs=pU-a.e.}) imply that $q_s$ does not depend on the choice of $s$, that is, for any $0<s,s'<t$,
  \begin{equation*}
    q_s(t,x,y) = q_{s'}(t,x,y),\quad \text{ for all } x,y \in D.
  \end{equation*}
  Hence, we can define
  \begin{equation*}
    q(t,x,y) = q_s(t,x,y),\quad t>0,~x,y \in D.
  \end{equation*}
  This together with (\ref{eq:qs=pU-a.e.}) yields (i). By the definition of $q_s$ and (\ref{eq:qs=pU-a.e.}), we obtain for all $t >0$ and $x,y \in D$,
  \begin{equation*}
    q(t,x,y) = q_{t/2}(t,x,y) = \int_D q(t/2,x,z)q(t/2,y,z)dz= \int_D q(t/2,y,z)q(t/2,x,z)dz = q_{t/2}(t,y,x) = q(t,y,x),
  \end{equation*}
  which is (ii). By the definition of $q_s$, (\ref{eq:qs=pU-a.e.}) and the symmetry of $q$, we obtain for all $0<s<t$ and $x,y \in D$,
  \begin{align*}
    q(t,x,y) =&~ q_s(t,x,y) = \int_D p_D(t-s,x,z) p_D(s,y,z)dz \\
    =&~ \int_D q(t-s,x,z) q(s,y,z)dz\\
    =&~ \int_D q(t-s,x,z) q(s,z,y)dz,
  \end{align*}
  which is (iii). The proof is complete.
\end{proof}

\medskip

\begin{proof}[Proof of Theorem \ref{thm:pD-exist}]
  Theorem \ref{thm:pD-exist} follows directly from Lemma \ref{lem:dhk-def-1} and Theorem \ref{thm:qHolder}. Indeed, we need only to rename the heat kernel $q(t,x,y)$ from Theorem \ref{thm:qHolder} to $p_D(t,x,y)$.
\end{proof}

In the sequel, for any open set $D \subset \mathbb{R}^d$, we always use $p_D(t,x,y)$ to denote the (locally) H\"older continuous version of the heat kernel obtained in Theorem \ref{thm:pD-exist}.

\section{Harmonic measures}\label{S:3}

The estimates in Lemmas \ref{lem:alphaTest1} and  \ref{lem:alphaTest2} can be found   in \cite[Lemma 2.3]{Ros-OtonSerra.2016.Duke}.   For the reader's convenience,  we provide a  proof here that  adopts the approach from \cite[Lemmas 4.1 and 5.1]{BogdanBurdzyChen.2003.PTRF89}  for censored stable processes in upper half space and \cite[Lemma 2.1]{ChenKimSongEtAl.2012.TAMS4169} for the process that is the independent sum of Brownian motion and isotropic stable processes.

\const{\CTtest}
\begin{lemma}\label{lem:alphaTest1}
  Let $\Delta^{\alpha/2}$ be the fractional Laplacian on $\mathbb{R}$. Let $p > 0$ and define the function $w_p(x) = (x\vee 0)^p,~ x \in \mathbb{R}$. Then, there is a constant $C_{\CTtest} = C_{\CTtest}(p,\alpha)$ such that
  \begin{equation*}
  \left(\Delta^{\alpha/2} w_p\right)(x) = C_{\CTtest} x^{p-\alpha},\quad x >0,
  \end{equation*}
  with
  \begin{equation*}
    C_{\CTtest}
    \begin{cases}
      > 0,&\quad p\in(\frac{\alpha}{2},\alpha),\\
      = 0,&\quad p = \frac{\alpha}{2},\\
      < 0,&\quad p \in (0,\frac{\alpha}{2}).
    \end{cases}
  \end{equation*}
\end{lemma}
\begin{proof}
  Let $p \in (0,\alpha)$. Recall $\sC_{1,\alpha}$ is the constant defined in \eqref{e:1.2}. By the definition of $\Delta^{\alpha/2}$ (see \eqref{e:1.0} for $d=1$), we have for any $x > 0$,
  \begin{align*}
    \left(\Delta^{\alpha/2} w_p\right)(x) =&~ \frac{\sC_{1,\alpha}}{2} \int_{-\infty}^\infty \frac{w_p(x+z)+w_p(x-z) -2w_p(x)}{|z|^{1+\alpha}} dz\\
    =&~ \sC_{1,\alpha}\int_0^\infty \frac{w_p(x+z)+w_p(x-z) -2w_p(x)}{|z|^{1+\alpha}} dz\\
    =&~ \sC_{1,\alpha} \int_0^\infty \frac{w_p(x+xt)+w_p(x-xt) -2w_p(x)}{(xt)^{1+\alpha}} xdt\\
    =&~ \sC_{1,\alpha}  x^{p-\alpha}\int_0^\infty \frac{w_p(1+t)+w_p(1-t) -2w_p(1)}{t^{1+\alpha}} dt\\
    =&~ \left(\Delta^{\alpha/2} w_p\right)(1) \cdot x^{p-\alpha}\\
    =&\!\!:\, C_{\CTtest}(p,\alpha) x^{p-\alpha}.
  \end{align*}
  Note that
  \begin{equation*}
    (\sC_{1,\alpha})^{-1}C_{\CTtest} = \left(\int_0^1 + \int_1^\infty\right) \frac{w_p(1+t)+w_p(1-t) -2w_p(1)}{t^{1+\alpha}} dt =: I_1 + I_2.
  \end{equation*}
    By integration by parts,
  \begin{align*}
    I_2 =&~ \int_1^\infty\frac{(1+t)^p -2}{t^{1+\alpha}} dt= \int_1^\infty \left((1+t)^p -2\right) d(-\frac{1}{\alpha}t^{-\alpha})\\
    =&~ -\frac{1}{\alpha} t^{-\alpha}((1+t)^p -2)\Big|_1^\infty + \frac{1}{\alpha} \int_1^\infty t^{-\alpha} \cdot p (1+t)^{p-1} dt\\
    =&~ \frac{1}{\alpha}(2^p-2) + \frac{p}{\alpha} \int_1^\infty t^{-\alpha}(1+t)^{p-1} dt.
  \end{align*}
  By integration by parts again, we have
  \begin{align*}
    I_1 =&~ \int_0^1 \frac{(1+t)^p+(1-t)^p-2}{t^{1+\alpha}} dt= \int_0^1 ((1+t)^p +(1-t)^p-2) d(-\frac{1}{\alpha}t^{-\alpha})\\
    =&~ -\frac{1}{\alpha} t^{-\alpha}((1+t)^p+(1-t)^p-2)\Big|_0^1 + \frac{1}{\alpha} \int_0^1 t^{-\alpha} \cdot p ((1+t)^{p-1}-(1-t)^{p-1}) dt\\
    =&~ -\frac{1}{\alpha}(2^p-2) + \frac{p}{\alpha} \int_0^1 t^{-\alpha}((1+t)^{p-1}-(1-t)^{p-1})  dt.
  \end{align*}
  Combining the above three formulas, we obtain
  \begin{equation*}
    (\sC_{1,\alpha})^{-1} C_{\CTtest} = \frac{p}{\alpha} \int_0^\infty t^{-\alpha}(1+t)^{p-1} dt - \frac{p}{\alpha} \int_0^1 t^{-\alpha}(1-t)^{p-1}  dt.
  \end{equation*}
  By a change of variable $t= \frac{r}{1+r}$,
  \begin{align*}
    \int_0^1 t^{-\alpha}(1-t)^{p-1} dt =&~ \int_0^\infty \left(\frac{r}{1+r}\right)^{-\alpha} \left(1-\frac{r}{1+r}\right)^{p-1} d (\frac{r}{1+r})  \\
    =&~ \int_0^\infty r^{-\alpha}(1+r)^{\alpha-p-1} dr.
  \end{align*}
  Combining the above two formulas, we obtain
  \begin{equation*}
    \frac{\alpha}{p \sC_{1,\alpha}} C_{\CTtest} = \int_0^\infty t^{-\alpha}(1+t)^{p-1}(1-(1+t)^{\alpha-2p}) dt.
  \end{equation*}
  Note that $1+t > 1$ for all $t >0$, we obtain that
  \begin{equation*}
    \frac{\alpha}{p \sC_{1,\alpha}}C_{\CTtest}
    \begin{cases}
      > 0,&\quad p\in(\frac{\alpha}{2},\alpha),\\
      = 0,&\quad p = \frac{\alpha}{2},\\
      < 0,&\quad p \in (0,\frac{\alpha}{2}).
    \end{cases}
  \end{equation*}
  The proof is complete.
\end{proof}

Let $\Pi$ be a hyperplane described by the function
\begin{equation*}
  \Phi(x) = (a,x-x_0),\qquad x \in \mathbb{R}^d,
\end{equation*}
where $x_0 \in \mathbb{R}^d$, $0\neq a = (a^{(1)},a^{(2)},\cdots,a^{(d)}) \in \mathbb{R}^d$ and $(\cdot,\cdot)$ is the inner product in $\mathbb{R}^d$. That is, $\Pi = \{y: \Phi(y) = 0\}$. We define $\delta_{\Pi}(y) = (\Phi(y)\vee 0)|a|^{-1}$ to be the distance from $y$ to the lower half space separated by the hyperplane $\Pi$.  For simplicity, we denote $-(-\partial_{x^{(k)}x^{(k)}}^2)^{\alpha/2}$ by $\Delta_k^{\alpha/2}$ for $1\le k\le d$. That is, by \eqref{e:1.1-2} and \eqref{e:1.1-3},
\begin{align}
  \Delta^{\alpha/2}_k f(x) =&~ \frac{\sC_{1,\alpha}}{2}\int_{\mathbb{R}} \frac{(f(x+te_k) + f(x-te_k) - 2f(x)}{|t|^{1+\alpha}} dt,\notag\\
  =&~ \lim_{\varepsilon\rightarrow 0+}  \sC_{1,\alpha} \int_{|t|>\varepsilon} \frac{(f(x+te_k) - f(x)}{|t|^{1+\alpha}} dt.\label{eq:Delta-k}
\end{align}
Recall that $e_k$ is the unit vector in the positive $x^{(k)}$-direction.
%ZC deleted: {\blue  The following lemma can be also found in \cite[Lemma 2.3]{Ros-OtonSerra.2016.Duke}. We give the details of the proof for the readers' convenience.}

\const{\harU}
\const{\harL}
\begin{lemma}\label{lem:alphaTest2}
  Let $p > 0$ and suppose that $\Pi$ and $\delta_{\Pi}$ are defined as above. Then, there are two constants $C_i=C_i(d,\alpha,p) > 0,i=\harU,\harL$ such that for every $x \in \mathbb{R}^d$ with $\Phi(x) >0$,
  \begin{align}
    C_{\harL} \delta_{\Pi}(x)^{p-\alpha} \le &\left(\mathcal{L} \delta_{\Pi}^p\right)(x) \le C_{\harU} \delta_{\Pi}(x)^{p-\alpha}, && \text{ if }  ~p \in (\frac{\alpha}{2},\alpha),\label{eq:alpha+}\\
    &\left(\mathcal{L} \delta_{\Pi}^p\right)(x) = 0,&&  \text{ if }~p = \frac{\alpha}{2},\label{eq:alpha}\\
    -C_{\harU} \delta_{\Pi}(x)^{p-\alpha} \le &\left(\mathcal{L} \delta_{\Pi}^p\right)(x) \le -C_{\harL} \delta_{\Pi}(x)^{p-\alpha}, && \text{ if } ~p \in(0, \frac{\alpha}{2}).\label{eq:alpha-}
  \end{align}
\end{lemma}
\begin{proof}
  By the definitions of $\Delta^{\alpha/2}_k$ and $\delta_{\Pi}$, for fixed $x \in \mathbb{R}^d$ with $\Phi(x) >0$, we have
  \begin{align*}
   \left(\Delta^{\alpha/2}_k\delta_{\Pi}^p\right)(x) =& \frac{\sC_{1,\alpha}}{2|a|^p}\int_{\mathbb{R}} \frac{(\Phi(x+te_k)\vee0)^p + (\Phi(x-te_k)\vee0)^p - 2(\Phi(x)\vee0)^p}{|t|^{1+\alpha}} dt\\
    =& \frac{\sC_{1,\alpha}}{2|a|^p} \int_{\mathbb{R}} \frac{((\Phi(x)+a^{(k)}t)\vee0)^p + ((\Phi(x)-a^{(k)}t)\vee0)^p - 2(\Phi(x))^p}{|t|^{1+\alpha}} dt\\
    =& \frac{\sC_{1,\alpha}\Phi(x)^{p-\alpha}}{2|a|^p} \int_{\mathbb{R}} \frac{((1+a^{(k)}t)\vee0)^p + ((1-a^{(k)}t)\vee0)^p - 2}{|t|^{1+\alpha}} dt\\
    =& \frac{\sC_{1,\alpha}\Phi(x)^{p-\alpha}}{2|a|^p} |a^{(k)}|^\alpha \int_{\mathbb{R}} \frac{((1+t)\vee0)^p + ((1-t)\vee0)^p - 2}{|t|^{1+\alpha}} dt\\
    =& \frac{\Phi(x)^{p-\alpha}}{|a|^p} |a^{(k)}|^\alpha \cdot \left(\Delta^{\alpha/2} w_p\right)(1) = C_{\CTtest}\frac{\Phi(x)^{p-\alpha}}{|a|^p}|a^{(k)}|^\alpha,
  \end{align*}
  where $C_{\CTtest}$ is the constant in Lemma \ref{lem:alphaTest1}. Note that, the last equality still holds for $a^{(k)} = 0$. Combining the last equality and the fact that
  \begin{equation*}
    d^{-\alpha/2} |a|^\alpha = d^{-\alpha/2} \left(\sum_{k=1}^d |a^{(k)}|^2\right)^{\alpha/2} \le \sum_{k=1}^d |a^{(k)}|^\alpha \le d \left(\sum_{k=1}^d |a^{(k)}|^2\right)^{\alpha/2} = d |a|^{\alpha},
  \end{equation*}
  we have
  \begin{align*}
   \left(\mathcal{L}\delta_{\Pi}^p\right)(x) =&~  C_{\CTtest} \frac{\Phi(x)^{p-\alpha}}{|a|^p}\left(\sum_{k=1}^d |a^{(k)}|^\alpha\right) \stackrel{d}{\asymp} C_{\CTtest} \frac{\Phi(x)^{p-\alpha}}{|a|^p} \cdot |a|^\alpha = C_{\CTtest} \delta_{\Pi}^{p-\alpha}(x).
  \end{align*}
  Therefore, by Lemma \ref{lem:alphaTest1}, we finish the proof.
\end{proof}

\const{\CTLAlpha}
\const{\CTLAlphaU+}
\const{\CTLAlphaL+}
Let $D\subset \mathbb{R}^d$ be a $C^{1,1}$ open set with characteristics $(R,\Lambda)$. For $Q\in \partial D$ and the coordinate system $CS_Q$, we define
\begin{equation*}
  \rho_Q(y) := y^{(d)} - \phi_Q(\tilde{y}),\qquad y:=(\tilde{y},y^{(d)}) \in CS_Q.
\end{equation*}
Note that for every $Q\in \partial D$ and $y \in B(Q,R)\cap D$, we have
\begin{equation*}
  \frac{\rho_Q(y)}{\sqrt{1+\Lambda^2}} \le \delta_D(y) \le \rho_Q(y).
\end{equation*}
Set
\begin{equation}\label{eq:R0r0-def}
  R_0 :=R_0(R,\Lambda) = \frac{R}{\sqrt{1+\Lambda^2}}\quad \text{ and }\quad r_0 = r_0(R,\Lambda) = \frac{R}{4(1+\Lambda^2)}.
\end{equation}

\begin{lemma}\label{thm:harmonicFunEst}
  Let $Q\in \partial D$ and fix the coordinate system $CS_Q$ so that
  \begin{equation*}
    B(Q,R) \cap D = \{y = (\tilde{y},y^{(d)}) \in B(0,R) \text{ in } CS_Q: y^{(d)} > \phi_Q(\tilde{y})\}.
  \end{equation*}
  For $p \in [\alpha/2,\alpha)$, we define
  \begin{equation*}
    h_p(y) = \left(\rho_Q(y)\right)^p \1_{D\cap B(Q,R_0)}(y),\qquad y \in \mathbb{R}^d.
  \end{equation*}
  Then, there exist $C_i = C_i(d,\alpha,R,\Lambda,p) > 0$, $i=\CTLAlpha,\CTLAlphaU+,\CTLAlphaL+$ such that for all $x\in D$ with $\rho_Q(x) < r_0$ and $|\tilde{x}| < r_0$,
  \begin{enumerate}
    \item[\rm (1)] if $p = \alpha/2$, then, we have
        \begin{equation}\label{eq:Lhax}
          \left|\mathcal{L} h_p(x)\right| \le C_{\CTLAlpha}|\ln\rho_Q(x)|,
        \end{equation}
        \item[\rm (2)] if $\alpha/2<p<\alpha$, then, we have
        \begin{equation}\label{eq:Lhpx}
          C_{\CTLAlphaL+} \left(\rho_Q(x)\right)^{p-\alpha} \le \mathcal{L} h_p(x) \le C_{\CTLAlphaU+} \left(\rho_Q(x)\right)^{p-\alpha}.
        \end{equation}
  \end{enumerate}
\end{lemma}

\begin{proof}    Note that any $C^{1,1}$ open set is  locally very close to the upper half space. We will use this property and apply Lemma \ref{lem:alphaTest2} to prove this lemma.

  Fix $x =(\tilde{x},x^{(d)}) \in CS_Q$ with $\rho_Q(x) < r_0$ and $|\tilde{x}| < r_0$, and choose $x_0\in \partial D$ with $\tilde{x} = \tilde{x}_0$. See Figure \ref{fig:Qx} for a special case.

  \begin{figure}[htp]
  \begin{tikzpicture}[
    declare function = {
        thetaA = 215;
        thetaB = 325;
        thetax = 255;
        r = 3;
        range=0.4;
        thetas = 244;
        thetat = 296;
        },
      scale = 2,
    ]
    \coordinate (Q) at ({r*cos((thetaA+thetaB)/2)},{r*sin((thetaA+thetaB)/2)});

    \fill [color= gray!40] ({r*cos(thetas)},{r*sin(thetas)}) arc (thetas:thetat:r) -- ($({r*cos(thetat)},{r*sin(thetat)})+(0,{range*r})$) arc (thetat:thetas:r) -- cycle;

    \draw ({r*cos(thetaA)},{r*sin(thetaA)}) arc (thetaA:thetaB:r);

    \fill (Q) circle (1pt) node [above] {$Q$};
    \coordinate (Q1) at ($(Q)!0.3*r!90:(0,0)$);
    \coordinate (Q2) at ($(Q)!0.3*r!-90:(0,0)$);
    \draw[dashed] (Q1) -- (Q2);

    \coordinate (x0) at ({r*cos(thetax)},{r*sin(thetax)});
    \coordinate (L1) at ($(x0)!0.23*r!90:(0,0)$);
    \coordinate (L2) at ($(x0)!0.36*r!-90:(0,0)$);
    \draw (L1) -- (L2);
    \fill (x0) circle (0.7pt) node [above] {$x_0$};

    \coordinate (x) at ($(x0)+(0,0.25*r)$);
    \fill (x) circle (0.7pt) node [above] {$x$};

    \coordinate (ya) at ($(x0)+({-0.5*r*sin(thetax)},{0.5*r*cos(thetax)})$);
    \draw[rotate=thetax+90] (ya) rectangle ($(ya)+({0.03*r*cos(45)},{0.03*r*sin(45)})$);
    \coordinate (y) at ($(ya)+({-0.44*r*cos(thetax)},{-0.44*r*sin(thetax)})$);
    \fill (y) circle (0.7pt) node [above] {$y$};
    \draw (ya) -- (y);
    \coordinate (yaleft) at ($(ya)+({0.03*r*sin(thetax)},{-0.03*r*cos(thetax)})$);
    \coordinate (yleft) at ($(y)+({0.03*r*sin(thetax)},{-0.03*r*cos(thetax)})$);
    \draw [gray,thin,|<->|] (yaleft) -- (yleft) node [pos=0.5,sloped, above] {$\delta_\Pi(y)$};

    \coordinate (yb) at ($(ya)+({-0.44*r*tan(90-thetax)*sin(thetax)},{0.44*r*tan(90-thetax)*cos(thetax)})$);
    \draw (yb) -- (y);
    \fill (yb) circle (0.7pt);
    \coordinate (ybright) at ($(yb)+({0.2*r},0)$);
    \coordinate (yright) at ($(y)+({0.2*r},0)$);
    \draw [gray,very thin] (ybright) -- ($(yb)+({0.03*r},0)$);
    \draw [gray,very thin] (yright) -- ($(y)+({0.03*r},0)$);
    \draw [gray,thin,|<->|] (ybright) -- (yright) node [pos=0.5,sloped, below] {$\overline{h}(y)$};

    \coordinate (yc) at ($(y)+(0,{-0.27*r},0)$);
    \coordinate (ycright) at ($(yc)+({0.03*r},0)$);
    \fill (yc) circle (0.7pt);
    \draw [gray, thin,|<->|] (ycright) -- ($(y)+({0.03*r},0)$) node [pos=0.5,sloped, below] {$\rho_Q(y)$};

    \draw [->] ($({r*cos(thetaB)},{r*sin(thetaB)})+(0.2*r,-0.05*r)$) -- ($({r*cos(thetaB)},{r*sin(thetaB)})+(0.01*r,-0.02*r)$);
    \node [right] at ($({r*cos(thetaB)},{r*sin(thetaB)})+(0.2*r,-0.05*r)$) {$\partial D$};

    \draw [->] ($(x0)+({-1.08*r*sin(thetax)},{1.08*r*cos(thetax)})+(0.25*r,0.1)$) -- ($(x0)+({-1.08*r*sin(thetax)},{1.08*r*cos(thetax)})+(0.02*r,0)$);
    \node [right] at ($(x0)+({-1.08*r*sin(thetax)},{1.08*r*cos(thetax)})+(0.25*r,0.1)$) {$\Pi$};

    \coordinate (w1) at ($(Q)+(0.01*r,range*r+0.1*r)$);
    \coordinate (w2) at ($(Q)+(-0.02*r,(range*r-0.1*r)$);
    \draw [->] (w1) -- (w2);
    \node [above right] at (w1) {\!\!\!\!\!The set $\left\{z =(\tilde{z},z^{(d)}) \in CS_Q: \rho_Q(z) < r_0, ~|\tilde{z}| < r_0\right\}$};

  \end{tikzpicture}
  \caption{The points $Q$, $x$ and $x_0$, etc}
  \label{fig:Qx}
\end{figure}
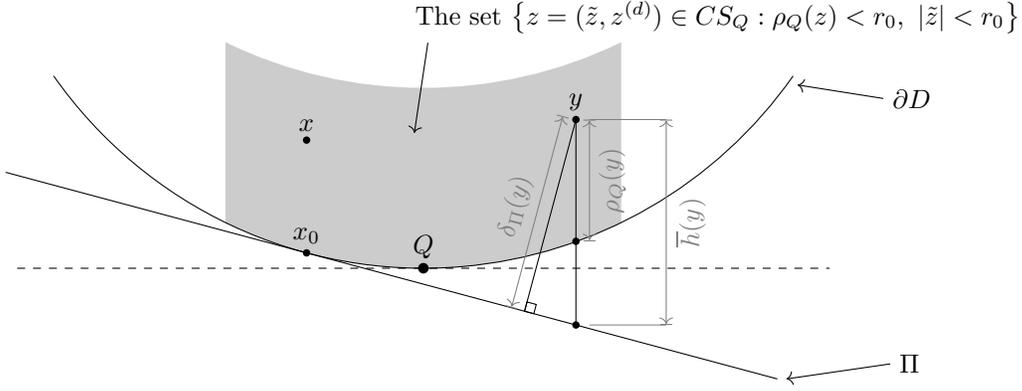

  Denote by $\Pi$ the hyperplane tangent to $\partial D$ at point $x_0$. Then, the function $\Gamma^*: \mathbb{R}^{d-1} \mapsto \mathbb{R}$ defined by $\Gamma^*(\tilde{y}) := \phi_Q(\tilde{x}_0) + \nabla\phi_Q(\tilde{x}_0)(\tilde{y}-\tilde{x})$ describes the plane $\Pi$. We use the following items:
  \begin{align*}
    \overline{h}(y) \coloneqq&~ \overline{h}_x(y) := (y^{(d)}-\Gamma^*(\tilde{y}))\vee 0,\\
    D_{\Gamma^*} =&~ \{y\in \mathbb{R}^d:y^{(d)}>\Gamma^*(\tilde{y})\},\\
    \delta_{\Pi}(y) =&~ {\rm dist}(y,\Pi)\1_{D_{\Gamma^*}}(y),\quad y \in \mathbb{R}^d,\\
    b_x \coloneqq &~\sqrt{1+|\nabla\phi_Q(\tilde{x})|^2},\\
    h_{x,p}(y) \coloneqq&~ (\overline{h}(y))^p, \ \ \text{ for }\alpha/2 \le p <\alpha.
  \end{align*}
  Note that $1\le b_x \le \sqrt{1+\Lambda^2}$ and $h_{x,p}(x) = h_p(x)$. Since $\overline{h}(y) = b_x\delta_{\Pi}(y)$, by   (\ref{eq:alpha+})
  and (\ref{eq:alpha}), we have for $y \in D_{\Gamma^*}$,
  \begin{equation}\label{eq:LhAalpha}
    \mathcal{L} h_{x,p}(y) = b_x^p \mathcal{L} \delta_{\Pi}^p(y) = 0,\quad p = \alpha/2,
  \end{equation}
  and
  \begin{equation*}
    C_{\harL} b_x^p(\delta_{\Pi}(y))^{p-\alpha} \le \mathcal{L} h_{x,p}(y) = b_x^p \mathcal{L} \delta_{\Pi}^p(y) \le C_{\harU} b_x^p (\delta_{\Pi}(y))^{p-\alpha},\quad \alpha/2 < p < \alpha.
  \end{equation*}
  Note that $b_x \delta_{\Pi}(x) = \rho_Q(x)$. By the last inequality, we have for $\alpha/2 < p < \alpha$,
  \begin{equation}\label{eq:Lhp}
  \begin{split}
    C_{\harL} (\rho_Q(x))^{p-\alpha} \le C_{\harL} b_x^\alpha (b_x \delta_{\Pi}(x))^{p-\alpha} &\le \mathcal{L} h_{x,p}(x) \\
    &\le C_{\harU} b_x^\alpha(b_x \delta_{\Pi}(x)))^{p-\alpha} \le C_{\harU} (1+\Lambda^2)^{\alpha/2} (\rho_Q(x))^{p-\alpha}.
  \end{split}
  \end{equation}
  We claim that,
  \begin{equation}\label{eq:hf-claim}
    \left|\mathcal{L}(h_p-h_{x,p})(x)\right| \le
    \begin{cases}
       c_1 < \infty, &p \in(\frac{\alpha}{2},\alpha),\\
       c_2|\ln \rho_Q(x)|, &p=\frac{\alpha}{2},
    \end{cases}
  \end{equation}
  for some constant $c_1=c_1(d,\alpha,p,R,\Lambda) > 0$ and $c_2=c_2(d,\alpha,R,\Lambda) > 0$, which together with (\ref{eq:LhAalpha}) and (\ref{eq:Lhp}) will establish this lemma.

  Let
  \begin{align*}
    A &\coloneqq \{y: \Gamma^*(\tilde{y}) < y^{(d)} \le \phi_Q(\tilde{y}) \text{ and } |\tilde{y}-\tilde{x}| < r_0\}\cup \{y: \Gamma^*(\tilde{y}) \ge y^{(d)} > \phi_Q(\tilde{y}) \text{ and } |\tilde{y}-\tilde{x}| < r_0\},\\
    E &\coloneqq \{y\in D\backslash A: |\tilde{y}-\tilde{x}| < r_0 \text{ and } \rho_Q(y) < r_0(2+\Lambda)\}.
  \end{align*}
  Note that, if $ y\in D\cap B(x,r_0)$, then
  \begin{equation*}
    \rho_Q(y) = y^{(d)} - \phi_Q(\widetilde{y})\le |y^{(d)}-x^{(d)}| + |x^{(d)}-\phi_Q(\tilde{x})|+|\phi_Q(\tilde{x}) - \phi_Q(\tilde{y})| < r_0(2+\Lambda).
  \end{equation*}
  If $|\tilde{y} - \tilde{x}| < r_0$ and $\rho_Q(y) < r_0(2+\Lambda)$, then
  \begin{align*}
    |y-Q|^2 =&~ |\tilde{y}|^2 + |y^{(d)}|^2 \\
    \le&~ (|\tilde{y}-\tilde{x}|+|\tilde{x}|)^2 + (|y^{(d)}-\phi_Q(\tilde{y})| + |\phi_Q(\tilde{y})|)^2 \\
    \le&~ (2r_0)^2 + (r_0(2+\Lambda) + |\phi_Q(\tilde{y})|)^2\\
    \le&~ (2r_0)^2 + (r_0(2+\Lambda) + \Lambda |\tilde{y}|)^2 < R_0^2.
  \end{align*}
  Consequently, we have
  \begin{equation*}
    D\cap B(x,r_0) \subset D\cap \{y:|\tilde{y}-\tilde{x}| < r_0 \text{ and } \rho_Q(y) < r_0(2+\Lambda)\} \subset D\cap B(Q,R_0),
  \end{equation*}
  and then,
  \begin{equation}\label{eq:AcsubsetE}
    A^c\cap D\cap B(x,r_0) \subset E\subset D\cap B(Q,R_0).
  \end{equation}
  We first consider the case for $p \in (\alpha/2,\alpha)$. Let $e_k = (0,\cdots,0,1,0,\cdots,0)\in\mathbb{R}^d$ be the unit vector along the $k^{th}$ axis for $1 \le k \le d$. Set
  \begin{align*}
    F :=&~ \{y := (\tilde{y},y^{(d)}) \in CS_Q : |\tilde{y}-\tilde{x}| < r_0\},\\
    H :=&~ \{y := (\tilde{y},y^{(d)}) \in CS_Q : y^{(d)} - \Gamma(\tilde{y}) > 0\}.
  \end{align*}
  Note that
  \begin{align*}
    &\{z : h_p(z) > 0\} = D\cap B(Q,R_0),\\
    &\{z : h_{x,p}(z) > 0\} = H,\\
    &\{z : h_p(z) - h_{x,p}(z) \neq 0 \}\subset (D\cap B(Q,R_0))\cup H \subset D\cup H,\\
    &A = (D^c \cap H \cap F) \cup (D\cap H^c\cap F),\\
    &A^c = (H\cup D)^c \cup (D\cap H) \cup F^c.
  \end{align*}
  Consequently, we have
  \begin{align*}
    A^c \cap (D\cup H) \subset&~  (A^c \cap D) \cup (A^c \cap D^c \cap H) = A^c \cap D,
  \end{align*}
  and then, by \eqref{eq:AcsubsetE} and the fact that $B(x,r_0) \subset F$,
  \begin{align*}
    B(x,r_0) \cap \{z : h_p(z) - h_{x,p}(z) \neq 0 \} \subset&~ B(x,r_0) \cap (D\cup H)\\
    =&~  B(x,r_0) \cap (A \cup (A^c \cap (D\cup H)))\\
    \subset&~ B(x,r_0) \cap (A \cup (A^c \cap D))\\
    \subset&~ A \cup (A^c \cap D \cap B(x,r_0)) \\
    \subset&~ A \cup E.
  \end{align*}
  Moreover, by the above formula and \eqref{eq:Delta-k}, we have for $1 \le k \le d$,
  \begin{align}\label{eq:decLhp}
    &(\sC_{1,\alpha})^{-1} \left|\Delta^{\alpha/2}_k (h_p-h_{x,p})(x)\right|\notag\\
    \le& \left|\int_{|t|\ge r_0} (h_p(x+te_k)-h_{x,p}(x+te_k)) \frac{dt}{|t|^{1+\alpha}} \right| + \lim_{\varepsilon\rightarrow 0} \int_{r_0>|t|>\varepsilon} |h_p(x+te_k)-h_{x,p}(x+te_k)| \frac{dt}{|t|^{1+\alpha}}  \notag\\
    \le&\int_{|t|\ge r_0} \left|h_p(x+te_k)-h_{x,p}(x+te_k)\right| \frac{dt}{|t|^{1+\alpha}} + \varliminf_{\varepsilon\rightarrow 0} \int_{\{r_0>|t|>\varepsilon\}\cap \{t: x+te_k \in A\}} \!\! |h_p(x+te_k)-h_{x,p}(x+te_k)| \frac{dt}{|t|^{1+\alpha}}\notag\\
    &+ \varliminf_{\varepsilon\rightarrow 0} \int_{\{r_0>|t|>\varepsilon\}\cap \{t:x+te_k\in E\}} |h_p(x+te_k)-h_{x,p}(x+te_k)| \frac{dt}{|t|^{1+\alpha}} \notag \quad \text{ (by (\ref{eq:AcsubsetE}))}\\
    =:& I_1 + I_2 + I_3.
  \end{align}
  We estimate $I_1,I_2,I_3$ separately. Note that for any $y\in\mathbb{R}^d$,
  \begin{equation*}
    0\le h_{x,p}(y) \le \left(|y^{(d)}-x^{(d)}|+|x^{(d)}-\phi_Q(\tilde{x})| + |\phi_Q(\tilde{x})-\Gamma^*(\tilde{y})|\right)^p \le (|x-y| + r_0 + \Lambda|x-y|)^p\le (r_0 + 2\Lambda|x-y|)^p.
  \end{equation*}
  Combining this and the facts that  $0\le h_p\le 1$ and $p < \alpha$, we have
  \begin{align*}
    I_1\le&~ \int_{|t|\ge r_0} \frac{1+(r_0+2\Lambda|t|)^p}{|t|^{1+\alpha}}dt\\
    \le&~ \int_{|t|\ge r_0} \frac{1+2^p(r_0^p+(2\Lambda)^p|t|^p)}{|t|^{1+\alpha}}dt\\
    =&~ 2\int_{r_0}^\infty \frac{1+(2r_0)^p + (4\Lambda)^p t^p}{t^{1+\alpha}} dt\\
    =&~ 2\left(\frac{1+(2r_0)^p}{\alpha}r_0^{-\alpha} + \frac{(4\Lambda)^p}{\alpha-p} r_0^{p-\alpha}\right) < \infty.
  \end{align*}
  For $y \in A$, we have
  \begin{align}\label{eq:BHP314}
    |h_p(y)| + |h_{x,p}(y)| \le&~ |y^{(d)} - \phi_Q(\tilde{y})|^p + |y^{(d)} -\Gamma^*(\tilde{y})|^p\notag \\
    \le&~ 2|\phi_Q(\tilde{y})-\Gamma^*(\tilde{y})|^p = 2|\phi_Q(\tilde{y}) - \phi_Q(\tilde{x}_0) - \nabla\phi_Q(\tilde{x}_0)(\tilde{y}-\tilde{x})|^p \notag\\
    \le&~  2\Lambda^p|\tilde{x}-\tilde{y}|^{2p} \le 2\Lambda^p|x-y|^{2p}
  \end{align}
  (see also \cite[(3.14)]{ChenKimSongEtAl.2012.TAMS4169}). Then, $|h_{x,p}(x+te_k)| + |h_p(x+te_k)| \le 2 \Lambda^p|t|^{2p}$. And since $p >\alpha/2$, we have
  \begin{equation*}
    I_2 \le \int_{|t|\le r_0} \frac{|h_{x,p}(x+te_k)| + |h_p(x+te_k)|}{|t|^{1+\alpha}} dt \le \int_0^{r_0}4\Lambda^{p}t^{2p-\alpha-1} dt = \frac{4\Lambda^p}{2p-\alpha} r_0^{2p-\alpha} < \infty.
  \end{equation*}

  For $y \in E\subset D\cap B(Q,R_0)$, we have $h_p(y) = \rho_Q(y)^p$. In view of this and the following two inequalities: for $a_1,a_2 > 0$,
  \begin{equation*}
    |a_1^p - a_2^p| \le
    \begin{cases}
      |a_1 - a_2|^p, &~~p\in (0,1),\\
      p(a_1\vee a_2)^{p-1}|a_1-a_2|, &~~ p\in [1,\infty),
    \end{cases}
  \end{equation*}
  we have for $y \in E$,
  \begin{equation}\label{eq:Difhp}
    |h_{x,p}(y) - h_p(y)| \le \left\{\begin{aligned}
                                 &|\overline{h}(y) - \rho_Q(y)|^p,\qquad  p\in (0,1),\\
                                 &p|\overline{h}(y) - \rho_Q(y)|,\qquad p\in [1,\infty).
                               \end{aligned}\right.
  \end{equation}
  On the other hand, by the definitions of $\overline{h}(y)$ and $\rho_Q(y)$, we have for $y \in E$,
  \begin{align}\label{eq:BHP37}
    |\overline{h}(y) - \rho_Q(y)| =&~ |\phi_Q(\tilde{y}) - \Gamma^*(\tilde{y})| = |\phi_Q(\tilde{y}) - \phi_Q(\tilde{x}_0) - \nabla\phi_Q(\tilde{x}_0)(\tilde{y}-\tilde{x})|  \notag \\
    \le&~ \Lambda |\tilde{x}-\tilde{y}|^2 \le \Lambda|x - y|^2,
  \end{align}
  (see also \cite[(3.7)]{ChenKimSongEtAl.2012.TAMS4169}). Using the last two inequalities, we have
  \setlength{\arraycolsep}{1.5pt}
  \begin{eqnarray*}
    I_3 & \le& \int_{\{r_0>|t|\}\cap \{t:x+te_k\in E\}} |h_p(x+te_k)-h_{x,p}(x+te_k)| \frac{dt}{|t|^{1+\alpha}}\\
    & \le& \begin{cases}
                   2\Lambda^p \int_0^{r_0} t^{2p-\alpha-1}dt =\frac{2\Lambda^p}{2p-\alpha} r_0^{2p-\alpha},
                      \quad  &\hbox{for } p\in (\alpha/2,1),  \medskip \\
                   2p\Lambda\int_0^{r_0} t^{1-\alpha}dt = \frac{2p\Lambda}{2-\alpha}r_0^{2-\alpha},   &\hbox{for  }  p\in [1,\infty).
                \end{cases}
  \end{eqnarray*}
  Combining (\ref{eq:decLhp}) and the estimates of $I_1,I_2$ and $I_3$,  and using the expression $\mathcal{L} = \sum_{k=1}^d\Delta^{\alpha/2}_k$,
  we can prove the first part of our claim  (\ref{eq:hf-claim}).

\medskip

  It remains to show the second part in (\ref{eq:hf-claim}) which is the case when $p = \alpha/2$. Similar to (\ref{eq:decLhp}), we have for $1\le k \le d$ and $p = \alpha/2$,
  \begin{align}\label{eq:decLhp2}
    &(\sC_{1,\alpha})^{-1} \left|\Delta^{\alpha/2}_k (h_p-h_{x,p})(x)\right|\notag\\
    \le& \left|\int_{|t|\ge r_0} (h_p(x+te_k)-h_{x,p}(x+te_k)) \frac{dt}{|t|^{1+\alpha}} \right| + \int_{r_0>|t|>\rho_Q(x)/(2\sqrt{1+\Lambda^2})} |h_p(x+te_k)-h_{x,p}(x+te_k)| \frac{dt}{|t|^{1+\alpha}} \notag\\
    &+  \varliminf_{\varepsilon\rightarrow 0} \int_{\varepsilon < |t| \le \rho_Q(x)/(2\sqrt{1+\Lambda^2})} |h_p(x+te_k)-h_{x,p}(x+te_k)| \frac{dt}{|t|^{1+\alpha}}\\
    \le&\int_{|t|\ge r_0} \left|h_p(x+te_k)-h_{x,p}(x+te_k)\right| \frac{dt}{|t|^{1+\alpha}} \notag\\
    &+ \int_{\{r_0>|t|>\rho_Q(x)/(2\sqrt{1+\Lambda^2})\}\cap \{t: x+te_k \in A\}} |h_p(x+te_k)-h_{x,p}(x+te_k)| \frac{dt}{|t|^{1+\alpha}}\notag\\
    &+ \int_{\{r_0>|t|>\rho_Q(x)/(2\sqrt{1+\Lambda^2})\}\cap \{t:x+te_k\in E\}} |h_p(x+te_k)-h_{x,p}(x+te_k)| \frac{dt}{|t|^{1+\alpha}} \notag\\
    &+  \varliminf_{\varepsilon\rightarrow 0} \int_{\varepsilon < |t| \le \rho_Q(x)/(2\sqrt{1+\Lambda^2})} |h_p(x+te_k)-h_{x,p}(x+te_k)| \frac{dt}{|t|^{1+\alpha}}\notag\\
    =:& J_1 + J_2 + J_3 + J_4.\notag
  \end{align}
  We can estimate $J_1$ similar to $I_1$, and have
  \begin{equation*}
    J_1 \le 2\left(\frac{1+(2r_0)^{\alpha/2}}{\alpha}r_0^{-\alpha} + \frac{(4\Lambda)^{\alpha/2}}{\alpha/2} r_0^{-\alpha/2}\right) < \infty.
  \end{equation*}
  Similar to $I_2$, by (\ref{eq:BHP314}) for $p = \alpha/2 \in (0,1)$, we have
  \begin{equation*}
    J_2 \le 2\int_{\rho_Q(x)/(2\sqrt{1+\Lambda^2})}^{r_0} 2   \Lambda^{\alpha/2}  |t|^{\alpha-\alpha-1} dt = (4\Lambda^{\alpha/2})(\ln(2r_0\sqrt{1+\Lambda^2}) - \ln \rho_Q(x)),
  \end{equation*}
  and similar to $I_3$, by (\ref{eq:Difhp}) for $p = \alpha/2$ and (\ref{eq:BHP37}),
  \begin{equation*}
    J_3 \le 2\int_{\rho_Q(x)/(2\sqrt{1+\Lambda^2})}^{r_0} \Lambda^{\alpha/2}|t|^{\alpha-\alpha-1} dt = (2\Lambda^{\alpha/2})(\ln(2r_0\sqrt{1+\Lambda^2}) - \ln \rho_Q(x)).
  \end{equation*}
  For $t \in (0,\rho_Q(x)/(2\sqrt{1+\Lambda^2})]$, we have
  \begin{equation*}
    \delta_D(x+te_k) \ge \delta_D(x) - |t| \ge \frac{\rho_Q(x)}{\sqrt{1+\Lambda^2}} - \frac{\rho_Q(x)}{2\sqrt{1+\Lambda^2}}= \frac{\rho_Q(x)}{2\sqrt{1+\Lambda^2}} > 0,
  \end{equation*}
  and by the definition of $\overline{h}$, for $y = x +te_k$,
   \begin{align}\label{eq:hy-lower}
    \overline{h}(y) =&~ y^{(d)} - \phi_Q(\tilde{x}_0) - \nabla\phi_Q(\tilde{x}_0)(\tilde{y}-\tilde{x}) \notag\\
    =&~ \overline{h}(x) + y^{(d)} - x^{(d)} - \nabla\phi_Q(\tilde{x}_0)(\tilde{y}-\tilde{x}) \notag\\
    \ge&~ \overline{h}(x) - \Lambda|t| \ge \rho_Q(x) - \frac{\Lambda\rho_Q(x)}{2\sqrt{1+\Lambda^2}}\ge \frac{\rho_Q(x)}{2} > 0.
  \end{align}
  The above two inequalities imply that for any $t \in (0,\rho_Q(x)/(2\sqrt{1+\Lambda^2})]$, $y:=x+te_k \in D$ (also, $y^{(d)} > \phi_Q(\tilde{y})$) and $y^{(d)} > \Gamma(\tilde{y})$ respectively. Hence, this together with \eqref{eq:AcsubsetE} implies that $x + te_k \in E$ for all $t \in (0,\rho_Q(x)/(2\sqrt{1+\Lambda^2})]$. Furthermore, combining (\ref{eq:BHP37}), \eqref{eq:hy-lower} and the following inequality:
  \begin{equation*}
    |a_1^{\alpha/2} - a_2^{\alpha/2}| \le a_1^{\alpha/2-1}|a_1-a_2|,\qquad a_1,a_2 > 0,
  \end{equation*}
  we have for $y= x+te_k$
  \begin{equation*}
    |h_{x,\alpha/2}(y) - h_{\alpha/2}(y)| \le (\overline{h}(y))^{\alpha/2-1}|\overline{h}(y) - \rho_Q(y)| \le \left(\frac{\rho_Q(x)}{2}\right)^{\alpha/2-1}\Lambda|x-y|^2  =\left(\frac{\rho_Q(x)}{2}\right)^{\alpha/2-1}\Lambda t^2.
  \end{equation*}
  Therefore,
  \begin{align*}
    J_4 \le&~ \int_0^{\rho_Q(x)/(2\sqrt{1+\Lambda^2})} 2\Lambda\left(\frac{\rho_Q(x)}{2}\right)^{\alpha/2-1}t^{1-\alpha} dt\\
    =&~ \frac{2\Lambda}{2-\alpha}(1+\Lambda^2)^{\alpha/2-1} \left(\frac{\rho_Q(x)}{2}\right)^{1-\alpha/2} \le\frac{2\Lambda r_0^{1-\alpha/2}}{2-\alpha}.
  \end{align*}
  Combining (\ref{eq:decLhp2}) and the estimates of $J_1,\cdots,J_4$, and using the expression $\mathcal{L} = \sum_{k=1}^d\Delta^{\alpha/2}_k$ yields the second part of our claim (\ref{eq:hf-claim}). In view of (\ref{eq:LhAalpha}), (\ref{eq:Lhp}) and our claim  (\ref{eq:hf-claim}), we get the desired results of this lemma.
\end{proof}

Recall that $\rho_Q(x):=x^{(d)} - \phi_Q(\tilde{x})$ for every $Q \in \partial D$ and
\begin{equation*}
  x \in B(Q,R) \cap D = \big\{y = (\tilde{y},y^{(d)}) \in B(0,R) \text{ in } CS_Q: y^{(d)} > \phi_Q(\tilde{y})\big\}.
\end{equation*}
We define for $r_1, r_2 > 0$
\begin{equation}\label{eq:Dr1r2-def}
  D(r_1,r_2) :=D_Q(r_1,r_2) :=\{y\in D: r_1 > \rho_Q(y) >0, ~|\tilde{y}| < r_2\}.
\end{equation}
Recall that the constants $R_0$ and $r_0$ are defined in (\ref{eq:R0r0-def}).

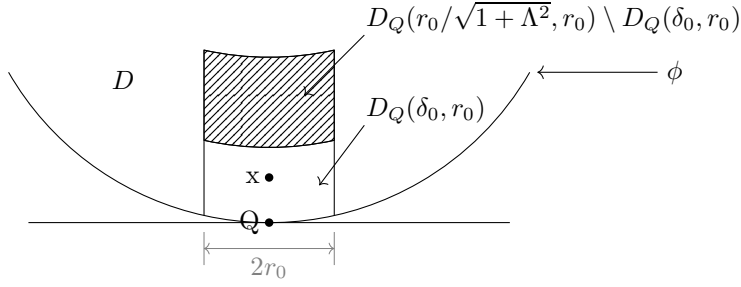
\begin{figure}[htp]
  \begin{tikzpicture}[
    declare function = {
        thetaS = 90;
        thetaA = 257.5;
        thetaB = 25;
        thetaC = 60;
        r = 2;
        sr1 = 0.5;
        sr2 = 1.1;
        },
      scale = 2,
      rotate = 0,
    ]

    \coordinate (a) at ({sr2*cos(thetaS) + r*cos(thetaA)},{sr2*sin(thetaS) + r*sin(thetaA)});
    \coordinate (b) at ({sr2*cos(thetaS) + r*cos(thetaA+thetaB)},{sr2*sin(thetaS) + r*sin(thetaA+thetaB)});
    \draw (a) arc (thetaA:thetaA+thetaB:r);

    \coordinate (c) at ({sr1*cos(thetaS) + r*cos(thetaA)},{sr1*sin(thetaS) + r*sin(thetaA)});
    \coordinate (d) at ({sr1*cos(thetaS) + r*cos(thetaA+thetaB)},{sr1*sin(thetaS) + r*sin(thetaA+thetaB)});
    \draw (c) arc (thetaA:thetaA+thetaB:r);

    \draw ({r*cos(270-thetaC)},{r*sin(270-thetaC)}) arc (270-thetaC:270+thetaC:r);

    \coordinate (e) at ({r*cos(thetaA)},{r*sin(thetaA)});
    \coordinate (f) at ({r*cos(thetaA+thetaB)},{r*sin(thetaA+thetaB)});
    \draw (a) -- (e);
    \draw (b) -- (f);

    \filldraw [pattern = north east lines] (a) arc (thetaA:thetaA+thetaB:r) -- (d) arc (thetaA+thetaB:thetaA:r) -- (c) -- cycle;

    \coordinate (Q) at ({r*cos(thetaA+thetaB/2)},{r*sin(thetaA+thetaB/2)});
    \coordinate (L1) at ($(Q)!0.8!90:(0,0)$);
    \coordinate (L2) at ($(Q)!0.8!-90:(0,0)$);
    \draw (L1) -- (L2);
    \fill (Q) circle (0.8pt) node [left] {Q};

    \coordinate (x) at ($(0,0)!0.85!(Q)$);
    \fill (x) circle (0.8pt) node [left] {x};

    \draw [->] ($(d)+(0.2,0.1)$) -- ($(d)+(-0.1,-0.3)$);
    \node [above right] at ($(d)+(0.15,0.05)$) {$D_Q(\delta_0,r_0)$};

    \draw [->] ($(b)+(0.2,0.1)$) -- ($(b)+(-0.2,-0.3)$);
    \node [above right] at ($(b)+(0.15,0.05)$) {$D_Q(r_0/\sqrt{1+\Lambda^2},r_0)\setminus D_Q(\delta_0,r_0)$};

    \coordinate (phi) at ({r*cos(270+thetaC)+0.05},{r*sin(270+thetaC)});
    \draw [->] ($(phi)+(0.8,0)$) -- (phi);
    \node [right] at ($(phi)+(0.8,0)$) {$\phi$};

    \coordinate (e1) at ($(a)!1.1!(e)$);
    \coordinate (e2) at ($(a)!1.3!(e)$);
    \draw [gray,very thin] (e1) -- (e2);
    \coordinate (f1) at ($(b)!1.1!(f)$);
    \coordinate (f2) at ($(b)!1.3!(f)$);
    \draw [gray,very thin] (f1) -- (f2);
    \coordinate (e3) at ($(a)!1.2!(e)$);
    \coordinate (f3) at ($(b)!1.2!(f)$);
    \draw [gray,very thin,<->] (e3) -- (f3) node [pos=0.5,sloped, below] {$2r_0$};

    \node [left] at ($(a)+(-0.4,-0.20)$) {$D$};
  \end{tikzpicture}
  \caption{The points $Q$ and $x$, and the set $D_Q(\delta_0,r_0)$, etc}
  \label{fig:DQdelta0r0}
\end{figure}

\const{\CTexitL}
\const{\CTexitU}
\begin{lemma}\label{lem:exitDistri}
  There are positive constants $\delta_0=\delta_0(d,\alpha,R,\Lambda)\in (0,r_0/(2\sqrt{1+\Lambda^2}))$ and $C_i=C_i(d,\alpha,R,\Lambda)$, $i = \CTexitL,\CTexitU$ such that for every $Q\in \partial D$ and $x \in D_Q(\delta_0,r_0)$ with $\tilde{x} = 0$,
  \begin{equation}\label{eq:exitLower}
    \mathbb{P}_x\left(X_{\tau_{D_Q(\delta_0,r_0)}} \in D_Q(r_0/\sqrt{1+\Lambda^2},r_0)\right) \ge C_{\CTexitL} \delta_D(x)^{\alpha/2},
  \end{equation}
  \begin{equation}\label{eq:exitUpper}
    \mathbb{P}_x\left(X_{\tau_{D_Q(\delta_0,r_0)}} \in D\right) \le C_{\CTexitU} \delta_D(x)^{\alpha/2},
  \end{equation}
  and
  \begin{equation}\label{eq:meanEupper}
    \mathbb{E}_x\left[\tau_{D_Q(\delta_0,r_0)} \right] \le C_{\CTexitU} \delta_D(x)^{\alpha/2}
  \end{equation}
  (cf. Figure \ref{fig:DQdelta0r0}).
\end{lemma}

\begin{proof}
  Recall the notation $w =(w^{(1)},w^{(2)},\cdots,w^{(d)}) :=(\tilde{w},w^{(d)}) \in \mathbb{R}^d$. Since $D$ is a $C^{1,1}$ open set with characteristics $(R,\Lambda)$, let $\phi :\mathbb{R}^{d-1}\mapsto \mathbb{R}$ be the $C^{1,1}$-function satisfying
  \begin{align*}
    &(1).\quad \phi( 0) =0,~~ \nabla\phi(0 ) =0;\\
    &(2).\quad \|\nabla\phi\|_\infty \le \Lambda;\\
    &(3).\quad |\nabla\phi(\tilde{y}) - \nabla\phi(\tilde{z})| \le \Lambda|\tilde{y} - \tilde{z}|.
  \end{align*}
  Let $CS_Q$ be the corresponding coordinate system such that
  \begin{equation*}
    B(Q,R) \cap D =  \left\{(\tilde{y},y^{(d)})\in B(0,R) \text{ in } CS_Q: y^{(d)} > \phi(\tilde{y}) \right\}.
  \end{equation*}
  Let $p \in (\alpha/2,\alpha)$ and define
  \begin{align*}
    \rho(y) :=&~ y^{(d)} - \phi(\tilde{y}),\\
    h(y) :=&~ \rho(y)^{\alpha/2} \1_{B(Q,R_0)\cap D}(y),\\
    h_p(y) :=&~ \rho(y)^p \1_{B(Q,R_0)\cap D}(y).
  \end{align*}
  Since $\rho(y) \le \sqrt{1+\Lambda^2}\delta_D(y)$, we have
  \begin{equation*}
    0\le h \le (\sqrt{1+\Lambda^2}R_0)^{\alpha/2} \le R^{\alpha/2} \le 1.
  \end{equation*}
  Let $\psi:\mathbb{R}^d \mapsto \mathbb{R}$ be a smooth positive function with bounded first and second order derivatives such that
  \begin{equation*}
    \psi(y) = \frac{2^{p+2}|\tilde{y}|^2}{r_0^2},\quad \text{ for } |y - Q| < r_0/(2\sqrt{1+\Lambda^2}),
  \end{equation*}
  and
  \begin{equation*}
    2^{p+2} \le \psi(y) \le 2^{p+3},\quad \text{ for } |y - Q| \ge r_0/\sqrt{1+\Lambda^2}.
  \end{equation*}
  Then, there exists $c_1 = c_1(d,\alpha,R,\Lambda,p)>0$ such that
  \begin{equation}\label{eq:Lpsi}
    \|\mathcal{L}\psi\|_\infty \le c_1.
  \end{equation}

\medskip

  \textbf{Step 1. Constructing suitable superharmonic and subharmonic functions with respect to $\mathcal{L}$.}

  We consider
  \begin{equation*}
    u_1(y):= h(y) + h_p(y),
  \end{equation*}
  and
  \begin{equation*}
    u_2(y):= h(y) + \psi(y) - h_p(y).
  \end{equation*}
  Since $p \in (\alpha/2,\alpha)$, by (\ref{eq:Lhax}) and (\ref{eq:Lhpx}), there exists $\delta_0 :=\delta_0(d,\alpha,R,\Lambda) \in (0,r_0/(2\sqrt{1+\Lambda^2}))$ such that for $y\in D_Q(\delta_0,r_0)$,
  \begin{equation}\label{eq:u1subharm}
    \mathcal{L}u_1(y) = \mathcal{L}h(y) + \mathcal{L}h_p(y)\ge -C_{\CTLAlpha} |\ln \rho(y)| + C_{\CTLAlphaL+} \rho(y)^{p-\alpha} \ge 0,
  \end{equation}
  and by (\ref{eq:Lpsi}),
  \begin{equation}\label{eq:u2superharm}
    \mathcal{L}u_2(y) = \mathcal{L}h(y) + \mathcal{L}\psi(y) - \mathcal{L}h_p(y)\le C_{\CTLAlpha} |\ln \rho(y)| + c_1 - C_{\CTLAlphaL+} \rho(y)^{p-\alpha} \le -1.
  \end{equation}

\medskip

  \textbf{Step 2. Translating super/subharmonic functions into super/submartingale properties for $X_t$.}

  We claim that the inequalities (\ref{eq:u1subharm}) and (\ref{eq:u2superharm}) imply that
  \begin{equation}\label{eq:u2supermar}
    t\mapsto u_2(X_{t\wedge \tau_{D(\delta_0,r_0)}}) +t\wedge \tau_{D(\delta_0,r_0)} \text{ is a non-negative bounded supermartingale},
  \end{equation}
  \begin{equation}\label{eq:MEdelta1}
    \mathbb{E}_x\left[\tau_{D(\delta_0,r_0)}\right] \le \rho(x)^{\alpha/2},
  \end{equation}
  and
  \begin{equation}\label{eq:u1submar}
    t\mapsto u_1(X_{t\wedge \tau_{D(\delta_0,r_0)}}) \text{ is a non-negative bounded submartingale}.
  \end{equation}
  Recall that $\tau_{D(\delta_0,r_0)}$ is the first exit time of $X$ upon leaving the set $D(\delta_0,r_0)$.

  Observe that if $v$ is a bounded $C^2$-function on $\mathbb{R}^d$ with bounded second order derivatives, then, by Markov property,
  \begin{equation}\label{eq:Mvmartingale}
    M_t^v :=v(X_t) - v(X_0) - \int_0^t \mathcal{L}v(X_s) ds   \  \text{ is a martingale.}
  \end{equation}
  Hence, if $u_1$ and $u_2$ are $C^2$-functions with bounded second order derivatives, then the above claims would just follow from (\ref{eq:Mvmartingale}), (\ref{eq:u1subharm}) and (\ref{eq:u2superharm}).  However, $u_1$ and $u_2$ are not $C^2$-functions. We shall approximate them by smooth functions.  Indeed, let $g$ be a mollifier, and $g_n(z):=2^{nd}g(2^nz)$, $z\in \mathbb{R}^d$ for $n\ge 1$. Define
  \begin{equation*}
    u_i^{(n)}(z) := g_n\ast u_i(z) = \int_{\mathbb{R}^d} g_n(y)u_i(z-y) dy,\qquad i=1,2.
  \end{equation*}
  Since
  \begin{equation*}
    \mathcal{L} u_i^{(n)} =g_n\ast \mathcal{L} u_i,\qquad i=1,2,
  \end{equation*}
  we have for any $n> m\ge 1$,
  \begin{equation*}
    \mathcal{L} u_1^{(n)} \ge 0,\quad \text{ and } \quad \mathcal{L}u_2^{(n)} \le -1,
  \end{equation*}
  on
  \begin{equation*}
    D_m(\delta_0,r_0) := \big\{y:\delta_0 - 2^{-m} > \rho(y) > 2^{-m}\text{ and }|\tilde{y}| < r_0 -2^{-m}\big\}.
  \end{equation*}
  Since each $u_i^{(n)}$ is a bounded smooth functions with bounded second order derivatives, it follows from (\ref{eq:Mvmartingale}) that for any $n> m\ge 1$,
  \begin{equation*}
    t\mapsto u_2^{(n)}(X_{t\wedge \tau_{D_m(\delta_0,r_0)}})+t\wedge \tau_{D_m(\delta_0,r_0)} \text{ is a non-negative supermartingale},
  \end{equation*}
  and
  \begin{equation*}
    t\mapsto u_1^{(n)}(X_{t\wedge \tau_{D_m(\delta_0,r_0)}}) \text{ is a non-negative bounded submartingale}.
  \end{equation*}
  Since each $u_i$ is bounded and continuous, $u_i^{(n)}$ converges to $u_i$ uniformly on $D_m(\delta_0,r_0)$, and hence,
  \begin{equation}\label{eq:u2supermar-2}
    t\mapsto u_2(X_{t\wedge \tau_{D_m(\delta_0,r_0)}})+t\wedge \tau_{D_m(\delta_0,r_0)} \text{ is a non-negative supermartingale},
  \end{equation}
  and
  \begin{equation*}
    t\mapsto u_1(X_{t\wedge \tau_{D_m(\delta_0,r_0)}}) \text{ is a non-negative bounded submartingale}.
  \end{equation*}
  Since $D_m(\delta_0,r_0)$ increases to $D(\delta_0,r_0)$ as $m\rightarrow \infty$, we obtain from above \eqref{eq:u2supermar} and (\ref{eq:u1submar}). Moreover, it follows from (\ref{eq:u2supermar-2}) that for each $n\ge 1$ and $t>0$,
  \begin{equation*}
    \mathbb{E}_x\left[u_2(X_{t\wedge \tau_{D_m(\delta_0,r_0)}})+t\wedge \tau_{D_m(\delta_0,r_0)}\right] \le u_2(x).
  \end{equation*}
  Since $u_2\ge 0$ and $D_m(\delta_0,r_0)$ increases to $D(\delta_0,r_0)$, by passing the above formula to the limit as $m\rightarrow \infty$ and then $t\rightarrow \infty$, we obtain
  \begin{equation*}
    \mathbb{E}_x\left[\tau_{D(\delta_0,r_0)}\right] \le u_2(x).
  \end{equation*}
  Note that $\tilde{x}=0$, $\psi(x) = 0$ and then, $u_2(x)\le \rho(x)^{\alpha/2}$. This together with the above inequality implies (\ref{eq:MEdelta1}). Consequently, (\ref{eq:meanEupper}) holds true.

\medskip

  \textbf{Step 3. Deriving the estimates of the exit distributions from super/submartingale properties}.

  Since $\psi(y) \ge 2^{p+2}$ for $|y-Q|>r_0/\sqrt{1+\Lambda^2}$ and $\phi(x) =0$, we have by (\ref{eq:u2supermar}),
  \begin{align*}
    \rho(x)^{\alpha/2} \ge&~ u_2(x)\\
    \ge&~ \mathbb{E}_x\left[u_2(X_{\tau_{D(\delta_0,r_0)}}); X_{\tau_{D(\delta_0,r_0)}} \in D\setminus D(r_0/\sqrt{1+\Lambda^2},r_0)\right]\\
    \ge&~ (2^{p+2}-1)\mathbb{P}_x\left(X_{\tau_{D(\delta_0,r_0)}} \in D\setminus D(r_0/\sqrt{1+\Lambda^2},r_0)\right).
  \end{align*}
  On the other hand, by (\ref{eq:u1submar}), we have
  \begin{equation*}
    \rho(x)^{\alpha/2} \le \rho(x)^{\alpha/2}+\rho(x)^p = u_1(x) \le \mathbb{E}_x\left[u_1(X_{\tau_{D(\delta_0,r_0)}})\right]\le 2\mathbb{P}_x\left(X_{\tau_{D(\delta_0,r_0)}} \in D\right).
  \end{equation*}
  Combining the above two formulas, we obtain
  \begin{align*}
    \mathbb{P}_x\left(X_{\tau_{D(\delta_0,r_0)}} \in D(r_0/\sqrt{1+\Lambda^2},r_0)\right) =&~ \mathbb{P}_x\left(X_{\tau_{D(\delta_0,r_0)}} \in D\right) - \mathbb{P}_x\left(X_{\tau_{D(\delta_0,r_0)}} \in D\setminus D(r_0/\sqrt{1+\Lambda^2},r_0)\right)\\
    \ge&~ \frac{2^{p+2}-3}{2(2^{p+2}-1)} \rho(x)^{\alpha/2},
  \end{align*}
  which implies (\ref{eq:exitLower}).

  Recall that $0\le h_p\le 1$. If $|y-Q|\ge r_0/\sqrt{1+\Lambda^2}$, then, $\psi(y) \ge 2^{p+2}$, we have
  \begin{equation*}
    u_2(y) = h(y)+\psi(y) -h_p(y) \ge 0+ 2^{p+2} - 1\ge 1,\qquad y \in B(Q,r_0/\sqrt{1+\Lambda^2})^c.
  \end{equation*}
  On the other hand, we have for $y\in B(Q,r_0)$ with $\delta_0 \le \rho(y) \le r_0$,
  \begin{equation*}
    u_2(y) \ge h(y)+\psi(y) -h_p(y) \ge \rho(y)^{\alpha/2} -\rho(y)^p \ge c_2,
  \end{equation*}
  where $c_2 = c_2(\delta_0,r_0,p) \in (0,1)$. It follows from the above two estimates that $u_2 \ge c_2$ on $D\setminus D(\delta_0,r_0)$. Therefore, by (\ref{eq:u2supermar}), we have
  \begin{equation*}
    \rho(x)^{\alpha/2} \ge u_2(x) \ge \mathbb{E}_x\left[u_2(X_{\tau_{D(\delta_0,r_0)}})\right] \ge c_2 \mathbb{P}_x\left(X_{\tau_{D(\delta_0,r_0)}} \in D\right),
  \end{equation*}
  which implies (\ref{eq:exitUpper}).
\end{proof}

\medskip

\section{Dirichlet heat kernel estimates}\label{S:4}

Throughout this section, $D\subset \mathbb{R}^d$ is a $C^{1,1}$ open set with characteristics $(R,\Lambda)$.

Recall that we use the following convention: for $u\in \mathbb{R}^d$, $a\in\mathbb{R}$ and $1\le i\le d$,
\begin{equation*}
  [u]_a^i := (u^{(1)},\cdots,u^{(i-1)},a,u^{(i+1)},\cdots,u^{(d)}).
\end{equation*}
Define
\begin{equation*}
  j(a,b) = \frac{\sC_{1,\alpha}}{|a-b|^{1+\alpha}} \quad \hbox{for } a \not= b \in \R,
\end{equation*}
where $\sC_{1,\alpha}$ is the positive constant in \eqref{e:1.2}.

With these notation, we rewrite the L\'evy system formula (\ref{eq:XLS-0}) as follows.  For any non-negative measurable function $f$ on $\mathbb{R}_+\times \mathbb{R}^d\times \mathbb{R}^d$ with $f(s, x, x)=0$ for any $s\geq 0$ and $x\in \R^d$ and for any stopping time $S$ with respect to the minimum augmented filtration generated by $X$, we have
\begin{equation}\label{eq:XLS}
  \mathbb{E}_x \left[\sum_{s\le S} f(s,X_{s-},X_s)\right] = \mathbb{E}_x \left[\int_0^S \sum_{i=1}^d \int_{\R} f(s,X_s,[X_s]^i_\theta) j(X_s^{(i)},\theta) d\theta ds\right].
\end{equation}

\subsection{Upper bound estimates}\label{S:4.1}

\begin{lemma}
  There is a constant $c = c(d,\alpha,R,\Lambda)>0$ such that for any $x\in D$,
  \begin{equation}\label{eq:tauD>1/4}
    \mathbb{P}_x\left(\tau_D \ge \frac{1}{4}\right) \le c\left(1\wedge \delta_D(x)^{\alpha/2}\right).
  \end{equation}
\end{lemma}

\begin{proof}
  Let $\delta_0,r_0$ be the constants from Lemma \ref{lem:exitDistri}. It suffices to prove (\ref{eq:tauD>1/4}) when $\delta_D(x) < \delta_0\wedge r_0 = \delta_0$. Indeed, let $Q \in \partial D$ be such that $\delta_D(x) = |x-Q|$, and $D(\delta_0,r_0)$ be the set defined in (\ref{eq:Dr1r2-def}). In this case, $x \in D_Q(\delta_0,r_0)$ and $\tilde{x} = 0$. It follows from (\ref{eq:exitUpper}) and (\ref{eq:meanEupper}) that
  \begin{align*}
    \mathbb{P}_x\left(\tau_D \ge \frac{1}{4}\right) \le&~ \mathbb{P}_x\left(\tau_D \ge \frac{1}{4}, \tau_{D(\delta_0,r_0)} \ge \frac{1}{4}\right) + \mathbb{P}_x\left(\tau_D \ge \frac{1}{4}, \tau_{D(\delta_0,r_0)} <\frac{1}{4}\right)\\
    \le&~ 4\mathbb{E}_x\left[\tau_{D(\delta_0,r_0)}\right] + \mathbb{P}_x\left(X_{\tau_{D(\delta_0,r_0)}} \in D\right)\\
    \le&~5C_{\CTexitU}\delta_D(x)^{\alpha/2}.
  \end{align*}
\end{proof}

\begin{lemma} \label{L:4.2}
  Let $U,U_1,U_3 \subset \mathbb{R}^d$ be three open sets with $U_1,U_3\subset U$ and  $\mathrm{dist}(U_1,U_3)>0$. Define $U_2:=U\setminus (U_1\cup U_3)$. We have, for any $t >0$, $x \in U_1$ and $y \in U_3$,
  \begin{equation}\label{eq:compUp}
  \begin{split}
    p_U(t,x,y) \le&~ \frac{2}{t} \mathbb{E}_x\left[\tau_{U_1}\right] \sup_{z\in U_1}p_U(t/2,z,y) + \mathbb{P}_x\left(X_{\tau_{U_1}} \in U_2\right) \sup_{t/2<s<t,z\in U_2} p_U(s,z,y)\\
    &~+ \int_0^{t/2} \int_{U_1} p_{U_1}(s,x,u) \left(\sum_{i=1}^d \int_{\mathbb{R}} p_U(t-s,[u]_\theta^i,y)d\theta \cdot \sup_{u \in U_1,[u]_a^i \in U_3} j(u^{(i)},a)\right)duds.
  \end{split}
  \end{equation}
\end{lemma}

\begin{proof}
  Fix $t > 0$ and $x \in U_1$. Let $0\le f \in L^1 (U)\cap L^\infty(U)$. By the strong Markov property  of $X$ and Proposition \ref{P:2.1}, we have
  \setlength{\arraycolsep}{1.5pt}
  \begin{eqnarray}
    P_t^Uf(x) & =& \mathbb{E}_x\left[f(X_t); t< \tau_U \right] \nonumber \\
    &=& \mathbb{E}_x\left[f(X_t); t<\tau_{U_1} \right] + \mathbb{E}_x\left[f(X_t); \tau_{U_1} \leq t< \tau_U\right]
    \nonumber \\
    &=&  P_t^{U_1}f(x) + \mathbb{E}_x\left[ f(X_t^U); \tau_{U_1} \leq t \right] \nonumber \\
    &=& P_t^{U_1}f(x) + \mathbb{E}_x\left[\mathbb{E}_{X_{\tau_{U_1}}^U} \left[f(X_{t-\tau_{U_1}}^U)\right]; \tau_{U_1}<t \right] +\mathbb{E}_x\left[f(X^U_t); \tau_{U_1} = t \right] \nonumber \\
    &=& P_t^{U_1}f(x) + \mathbb{E}_x\left[P_{t-\tau_{U_1}}^U f(X_{\tau_{U_1}}); \tau_{U_1}<t,~X_{\tau_{U_1}} \in U_2 \right] \nonumber \\
    && + \mathbb{E}_x\left[P_{t-\tau_{U_1}}^U f(X_{\tau_{U_1}}); \tau_{U_1}<t,~X_{\tau_{U_1}} \in U_3 \right]
         \nonumber \\
    &=:&  P_t^{U_1}f(x) + I+II .  \label{e:4.3}
  \end{eqnarray}
  Note that
  \begin{align*}
    I=&~ \mathbb{E}_x\left[\tau_{U_1}<t,~X_{\tau_{U_1}} \in U_2;P^{U}_{t-\tau_{U_1}}f(X_{\tau_{U_1}})\right]\\
    \le&~ \mathbb{P}_x\left(\tau_{U_1}<t,~X_{\tau_{U_1}} \in U_2\right)\sup_{0<s<t,~z\in U_2}P^{U}_sf(z)\\
    \le&~ \mathbb{P}_x\left(X_{\tau_{U_1}} \in U_2\right)\sup_{0<s<t,~z\in U_2}P^{U}_sf(z).
  \end{align*}
  Since $\mathrm{dist}(U_1,U_3)>0$, by (\ref{eq:XLS}), we obtain,
  \begin{align*}
    II =&~\mathbb{E}_x\left[\tau_{U_1}<t,~X_{\tau_{U_1}} \in U_3;P^U_{t-\tau_{U_1}}f(X_{\tau_{U_1}})\right]\\
    =&~\mathbb{E}_x\left[\int_0^t \1_{\{s<\tau_{U_1}\}} \cdot\left(\sum_{i=1}^d\int_{\mathbb{R}} \1_{\{[X_s]_\theta^i \in U_3\}}\cdot P^U_{t-s}f([X_s]_\theta^i)j(X_s^{(i)},\theta)d\theta\right)ds\right]\\
    =&~ \int_0^t \int_{U_1} p_{U_1}(s,x,u) \left(\sum_{i=1}^d\int_{\mathbb{R}} \1_{\{[u]_\theta^i \in U_3\}}\cdot P^U_{t-s}f([u]_\theta^i)j(u^{(i)},\theta)d\theta\right)du ds\\
    \le&~\int_0^t\int_{U_1} p_{U_1}(s,x,u) \left(\sum_{i=1}^d\int_{\mathbb{R}} P^U_{t-s}f([u]_\theta^i)d\theta \cdot \sup_{u\in U_1, [u]_a^i \in U_3}j(u^{(i)},a)\right)du ds.
  \end{align*}
  Thus we have by \eqref{e:4.3}
     \begin{align*}
    P_t^Uf(x) \le&~ P_t^{U_1}f(x) + \mathbb{P}_x\left(X_{\tau_{U_1}} \in U_2\right)\sup_{0<s<t,~z\in U_2}P^{U}_sf(z)\\
    &~+\int_0^t\int_{U_1} p_{U_1}(s,x,u) \left(\sum_{i=1}^d\int_{\mathbb{R}} P^U_{t-s}f([u]_\theta^i)d\theta \cdot \sup_{u\in U_1, [u]_a^i \in U_3}j(u^{(i)},a)\right)du ds.
  \end{align*}
  For any $y \in U_3$, by setting $f = p_U(t,\cdot,y)$ in the above inequality and using semigroup property, we obtain
  \begin{align*}
    p_U(2t,x,y) \le&~ \int_{U_1} p_{U_1}(t,x,z)p_U(t,z,y) dz + \mathbb{P}_x\left(X_{\tau_{U_1}} \in U_2\right)\sup_{t<s<2t,~z\in U_2}p_U(s,z,y)\\
    &~+\int_0^t\int_{U_1} p_{U_1}(s,x,u) \left(\sum_{i=1}^d\int_{\mathbb{R}} p_U(2t-s,[u]_\theta^i,y)d\theta \cdot \sup_{u\in U_1, [u]_a^i \in U_3}j(u^{(i)},a)\right)du ds\\
    \le&~ \mathbb{P}_x\left(\tau_{U_1} > t\right)\sup_{z\in U_1}p_U(t,z,y) + \mathbb{P}_x\left(X_{\tau_{U_1}} \in U_2\right)\sup_{t<s<2t,~z\in U_2}p_U(s,z,y)\\
    &~+\int_0^t\int_{U_1}  p_{U_1}(s,x,u) \left(\sum_{i=1}^d\int_{\mathbb{R}} p_U(2t-s,[u]_\theta^i,y)d\theta \cdot \sup_{u\in U_1, [u]_a^i \in U_3}j(u^{(i)},a)\right)du ds.
  \end{align*}
  By renaming $2t$ by $t$ in the above inequality, we finish the proof.
\end{proof}

\begin{lemma}
  There is a constant $c= c(d,\alpha,R,\Lambda)>0$ such that for all $x,y \in D$,
  \begin{equation}\label{eq:pD1/2Up}
    p_D(1/2,x,y) \le c\left(1\wedge \delta_D(x)^{\alpha/2}\right)p(1/2,x,y).
  \end{equation}
\end{lemma}
\begin{proof}
  It suffices to prove (\ref{eq:pD1/2Up}) when $\delta_D(x) < \delta_0$, where $\delta_0$ is the constant from Lemma \ref{lem:exitDistri}. Fix $x,y\in D$ with $\delta_D(x) < \delta_0$. Recall that $r_0 = r_0(R,\Lambda) = \frac{R}{4(1+\Lambda^2)}$. Take $r=4r_0$ for simplicity.

  \textbf{Case 1.}  For all $i=1,\cdots,d$, $|x^{(i)} - y^{(i)}| < r$. By semigroup property of $p_D$ and (\ref{eq:tauD>1/4}), we have
  \begin{align*}
    p_D(1/2,x,y) =&~ \int_D p_D(1/4,x,z)p_D(1/4,z,y)dz\\
    \le&~ \sup_{z\in D} p_D(1/4,z,y) \int_D p_D(1/4,x,z)dz\\
    \le&~ C_{\CTstable} 4^{d/\alpha} \mathbb{P}_x\left(\tau_D > 1/4\right)\\
    \le&~ c_1 \left(1\wedge \delta_D(x)^{\alpha/2}\right).
  \end{align*}
  On the other hand, since $|x^{(i)} - y^{(i)}| < r$, for all $i=1,\cdots,d$, we have by (\ref{eq:pEst}),
  \begin{equation*}
    p(1/2,x,y) \ge C_{\CTstable}^{-1} \prod_{i=1}^{d} \left(2^{1/\alpha} \wedge \frac{1/2}{r^{1+\alpha}}\right) > 0.
  \end{equation*}
  Combining the above two inequalities, we verify (\ref{eq:pD1/2Up}) in this case.

  \textbf{Case 2.}  There is some $1\le i\le d$  such that $|x^{(i)} - y^{(i)}| \ge r$. Let
  \begin{equation*}
    \mathcal{I} := \big\{i: |x^{(i)} - y^{(i)}| \ge r,~ 1\le i\le d\big\} ,
  \end{equation*}
  and $Q\in\partial D$ be such that $|x-Q| = \delta_D(x)$. Define
  \begin{align*}
    U_1 :=&~ D_Q(\delta_0,r_0)\quad \text{ (see \eqref{eq:Dr1r2-def} for the definition of the set $D_Q(\delta_0,r_0)$)},\\
    U_3 :=&~ \left\{z \in D: \exists~ i \in \mathcal{I}, \text{ such that } |z^{(i)} - x^{(i)}| > |x^{(i)} - y^{(i)}|/2\right\},\\
    U_2 :=&~ D \setminus U_1 \setminus U_3.
  \end{align*}
  Note that
  $U_1 \cap U_3 =\emptyset$, $x \in U_1$ by $\delta_D(x) < \delta_0$ and $y\in U_3$ by the definition of $\mathcal{I}$. By Lemma \ref{L:4.2} with $t = 1/2$ and $U=D$, we have
  \begin{equation}\label{eq:compUp-2}
  \begin{split}
    p_D(1/2,x,y) \le&~ 4\mathbb{E}_x\left[\tau_{U_1}\right]\sup_{z\in U_1}p_D(1/4,z,y) + \mathbb{P}_x\left(X_{\tau_{U_1}} \in U_2\right)\sup_{\frac{1}{4} < s < \frac{1}{2},z\in U_2} p_D(s,z,y)\\
    &~+ \int_0^{1/4} \int_{U_1} p_{U_1}(s,x,u)\left(\sum_{i=1}^d \int_{\mathbb{R}} p_D(1/2-s,[u]_\theta^i,y)d\theta \cdot\sup_{u \in U_1,[u]_a^i \in U_3} j(u^{(i)},a)\right)duds\\
    =:\!\!&~ \, I_1 + I_2 + I_3.
  \end{split}
  \end{equation}
  We estimate $I_1$, $I_2$, $I_3$ separately. Indeed, for any $i \in \mathcal{I}$ and $z \in U_1$, since $|z^{(i)} - x^{(i)}| < 2r_0$, we have
  \begin{equation*}
    |z^{(i)} - y^{(i)}| \ge |x^{(i)} - y^{(i)}| - |z^{(i)} - x^{(i)}| \ge |x^{(i)} - y^{(i)}| - 2r_0 \ge |x^{(i)} - y^{(i)}|/2.
  \end{equation*}
  This together with the upper bound of $p(1/4,z,y)$, the lower bound of $p(1/4,x,y)$ in (\ref{eq:pEst}) and (\ref{eq:meanEupper}) yields that
  \begin{align*}
    I_1 =&~ 4\mathbb{E}_x\left[\tau_{U_1}\right]\sup_{z\in U_1}p_D(1/4,z,y)\\
    \le&~ c_2\delta_D(x)^{\alpha/2} \sup_{z\in U_1} p(1/4,z,y)\\
    \le&~ c_2C_{\CTstable}\delta_D(x)^{\alpha/2} \sup_{z\in U_1}\prod_{i \in \mathcal{I}} \left(4^{1/\alpha}\wedge \frac{1/4}{|z^{(i)}-y^{(i)}|^{1+\alpha}}\right)\prod_{i \notin \mathcal{I}} \left(4^{1/\alpha}\wedge \frac{1/4}{|z^{(i)}-y^{(i)}|^{1+\alpha}}\right)\\
    \le&~ c_2C_{\CTstable}\delta_D(x)^{\alpha/2}  \prod_{i \in \mathcal{I}} \left(4^{1/\alpha}\wedge \frac{2^{1+\alpha}(1/4)}{|x^{(i)}-y^{(i)}|^{1+\alpha}}\right)\prod_{i \notin \mathcal{I}} \left(4^{1/\alpha}\right)\\
    \le&~ c_3\delta_D(x)^{\alpha/2} \prod_{i \in \mathcal{I}} \left(4^{1/\alpha}\wedge \frac{1/4}{|x^{(i)}-y^{(i)}|^{1+\alpha}}\right)\prod_{i \notin \mathcal{I}} \left(4^{1/\alpha}\wedge \frac{1/4}{|x^{(i)}-y^{(i)}|^{1+\alpha}}\right)\\
    \le&~ c_3C_{\CTstable}\delta_D(x)^{\alpha/2} p(1/4,x,y),
  \end{align*}
  where in the second to the last inequality, we have used the fact that for all $i \notin \mathcal{I}$,
  \begin{equation*}
    |x^{(i)}-y^{(i)}|^{1+\alpha} < r^{1+\alpha} < \infty.
  \end{equation*}
  By the definition of $U_2$, for any $z \in U_2$ and any $i \in \mathcal{I}$, we have
  \begin{equation*}
    |z^{(i)} - x^{(i)}| \le |x^{(i)} - y^{(i)}|/2,
  \end{equation*}
  and then
  \begin{equation*}
    |z^{(i)} - y^{(i)}| \ge |x^{(i)} - y^{(i)}| - |z^{(i)} - x^{(i)}| \ge |x^{(i)} - y^{(i)}| - |x^{(i)} - y^{(i)}|/2 = |x^{(i)} - y^{(i)}|/2.
  \end{equation*}
  This together with the upper bound of $p(s,z,y)$, the lower bound of $p(1/4,x,y)$ in (\ref{eq:pEst}) and (\ref{eq:exitUpper}) yields
  \begin{align*}
    I_2 =&~ \mathbb{P}_x\left(X_{\tau_{U_1}} \in U_2\right)\sup_{\frac{1}{4} < s < \frac{1}{2},z\in U_2} p_D(s,z,y)\\
    \le&~ \mathbb{P}_x\left(X_{\tau_{U_1}} \in D\right)\sup_{\frac{1}{4} < s < \frac{1}{2},z\in U_2} p(s,z,y)\\
    \le&~ c_4\delta_D(x)^{\alpha/2} \sup_{\frac{1}{4} < s < \frac{1}{2},z\in U_2}\prod_{i\in \mathcal{I}} \left(s^{-1/\alpha}\wedge \frac{s}{|z^{(i)}-y^{(i)}|^{1+\alpha}}\right)\prod_{i\notin \mathcal{I}} s^{-1/\alpha}\\
    \le&~ c_5\delta_D(x)^{\alpha/2} \prod_{i\in \mathcal{I}} \left(4^{1/\alpha}\wedge \frac{2^{1+\alpha}(1/2)}{|x^{(i)}-y^{(i)}|^{1+\alpha}}\right)\prod_{i\notin \mathcal{I}} \left(4^{1/\alpha}\wedge \frac{1/2}{|x^{(i)}-y^{(i)}|^{1+\alpha}}\right)\\
    \le&~ c_6\delta_D(x)^{\alpha/2} p(1/2,x,y).
  \end{align*}
  It remains to estimate $I_3$. Note that for all $u \in U_1$, we have
  \begin{equation*}
    |u^{(k)} - x^{(k)}| < 2r_0,\qquad k = 1,\cdots,d.
  \end{equation*}
  Hence, for all $i \notin \mathcal{I}$ and $u \in U_1$, by definition of $U_3$, it is not possible that there exists $a \in \mathbb{R}$ such that $[u]_a^i \in U_3$. In this case,
  \begin{equation}\label{eq:pD1/2I3-1}
    \sup_{u \in U_1,[u]_a^i \in U_3} j(u^{(i)},a) = 0.
  \end{equation}
  On the other hand, for any $i \in \mathcal{I}$ and $u\in U_1$, if $[u]_a^i \in U_3$ for some $a \in \mathbb{R}$, then the number $a$ must satisfy
  \begin{equation*}
    |a - x^{(i)}| > |x^{(i)} - y^{(i)}|/2 > 2r_0,
  \end{equation*}
   and so
  \begin{equation*}
    |a - u^{(i)}| \ge |a - x^{(i)}| - |x^{(i)} - u^{(i)}| \ge |x^{(i)} - y^{(i)}|/2 - r_0 \ge |x^{(i)} - y^{(i)}|/4.
  \end{equation*}
  In this case,
  \begin{equation}\label{eq:pD1/2I3-2}
    \sup_{u \in U_1,[u]_a^i \in U_3} j(u^{(i)},a) \le \sup_{u \in U_1,[u]_a^i \in U_3} \frac{\sC_{1,\alpha}}{|a - u^{(i)}|^{1+\alpha}} \le \frac{4^{1+\alpha}\sC_{1,\alpha}}{|x^{(i)} - y^{(i)}|^{1+\alpha}} \le c_7 \left(2^{1/\alpha}\wedge \frac{1/2}{|x^{(i)} - y^{(i)}|^{1+\alpha}}\right),
  \end{equation}
  and, for $k \in \mathcal{I}$ with $k\neq i$,
  \begin{equation}\label{eq:pD1/2I3-3}
    |u^{(k)} - y^{(k)}| \ge |x^{(k)} - y^{(k)}| - |u^{(k)} - x^{(k)}| \ge |x^{(k)} - y^{(k)}| - 2r_0 \ge |x^{(k)} - y^{(k)}|/2.
  \end{equation}
  Hence,
  \begin{align*}
    I_3 =&~ \int_0^{1/4} \int_{U_1} p_{U_1}(s,x,u)\left(\sum_{i=1}^d \int_{\mathbb{R}} p_D(1/2-s,[u]_\theta^i,y)d\theta \cdot\sup_{u \in U_1,[u]_a^i \in U_3} j(u^{(i)},a)\right)duds\\
    \le&~ \int_{1/4}^{1/2} \int_{U_1} p_{U_1}(1/2-s,x,u)\left(\sum_{i \in \mathcal{I}} \int_{\mathbb{R}} p(s,[u]_\theta^i,y)d\theta \cdot\sup_{u \in U_1,[u]_a^i \in U_3} j(u^{(i)},a)\right)duds   \quad (\text{by}~ (\ref{eq:pD1/2I3-1}))\\
    \le&~ c_7C_{\CTstable}\int_{1/4}^{1/2} \int_{U_1} p_{U_1}(1/2-s,x,u) \Bigg( \sum_{i \in \mathcal{I}} \sup_{\frac{1}{4} < s < \frac{1}{2}} \prod_{\genfrac{}{}{0pt}{4}{k \in \mathcal{I}}{k\neq i}} \left(s^{-1/\alpha}\wedge\frac{s}{|u^{(k)}-y^{(k)}|^{1+\alpha}}\right) \prod_{k \notin \mathcal{I}} \left(s^{-1/\alpha}\right)\\
    &~ \cdot\int_{\mathbb{R}} \left(s^{-1/\alpha}\wedge\frac{s}{|\theta-y^{(i)}|^{1+\alpha}} \right) d\theta \cdot\left(2^{1/\alpha}\wedge \frac{1/2}{|x^{(i)} - y^{(i)}|^{1+\alpha}}\right)\Bigg) duds   \quad (\text{by} ~ (\ref{eq:pD1/2I3-2}))\\
    \le&~ c_8\int_{1/4}^{1/2} \int_{U_1} p_{U_1}(1/2-s,x,u) \Bigg(\sum_{i \in \mathcal{I}} \prod_{\genfrac{}{}{0pt}{4}{k\in \mathcal{I}}{k\neq i}} \left(4^{1/\alpha}\wedge\frac{2^{1+\alpha}(1/2)}{|x^{(k)}-y^{(k)}|^{1+\alpha}}\right) \prod_{k \notin \mathcal{I}} \left(4^{1/\alpha}\right)
    \quad (\text{by} ~ (\ref{eq:pD1/2I3-3}))\\
    &~\cdot \left(2^{1/\alpha}\wedge \frac{1/2}{|x^{(i)} - y^{(i)}|^{1+\alpha}}\right)\Bigg) duds\\
    \le&~ c_8\int_0^{1/4} \int_{U_1} p_{U_1}(s,x,u) du ~ds \bigg( \sum_{i \in \mathcal{I}} \prod_{\genfrac{}{}{0pt}{4}{k\in \mathcal{I}}{k\neq i}} \left(4^{1/\alpha}\wedge\frac{2^{1+\alpha}(1/2)}{|x^{(k)}-y^{(k)}|^{1+\alpha}}\right) \prod_{k \notin \mathcal{I}} \left(4^{1/\alpha}\wedge\frac{1/2}{|x^{(k)}-y^{(k)}|^{1+\alpha}}\right)\\
    &~\cdot \left(2^{1/\alpha}\wedge \frac{1/2}{|x^{(i)} - y^{(i)}|^{1+\alpha}}\right)\bigg) \\
    \le&~ c_9\int_0^{1/4} \mathbb{P}_x\left(\tau_{U_1} > s\right) ds \cdot p(1/2,x,y)  \quad (\text{by the lower bound in (\ref{eq:pEst})})\\
    \le&~ c_9 \left(1 \wedge \mathbb{E}_x\left[\tau_{U_1}\right]\right) \cdot p(1/2,x,y)\\
    \le&~ c_9 \left(1 \wedge \delta_D(x)^{\alpha/2}\right) \cdot p(1/2,x,y) \quad \text{(by (\ref{eq:meanEupper}))}.
  \end{align*}
  Combining (\ref{eq:compUp-2}) and the estimates of $I_1,I_2,I_3$, we finish the case 2, and then this lemma.
\end{proof}

\begin{lemma}
  There is a constant $c= c(d,\alpha,R,\Lambda)>0$ such that for all $x,y \in D$,
  \begin{equation}\label{eq:pD1Up}
    p_D(1,x,y) \le c\left(1\wedge \delta_D(x)^{\alpha/2}\right)\left(1\wedge \delta_D(y)^{\alpha/2}\right)p(1,x,y).
  \end{equation}
\end{lemma}
\begin{proof}
  By the semigroup property of $p_D$ and (\ref{eq:pD1/2Up}), since $p_D(t,x,y)$ is symmetric in $x,y$, we have
  \begin{align*}
    p_D(1,x,y) =&~\int_{D}  p_D(1/2,x,z)p_D(1/2,z,y) dz\\
    \le&~ c \left(1\wedge \delta_D(x)^{\alpha/2}\right)\int_{\mathbb{R}^d} p(1/2,x,z)p(1/2,z,y) dz\left(1\wedge \delta_D(y)^{\alpha/2}\right)\\
    =&~ c\left(1\wedge \delta_D(x)^{\alpha/2}\right)\left(1\wedge \delta_D(y)^{\alpha/2}\right)p(1,x,y).
  \end{align*}
\end{proof}

\begin{lemma}\label{lem:pDtue}
  There is a constant $c= c(d,\alpha,R,\Lambda)>0$ such that for all $t \in (0,1]$, $x,y \in D$,
  \begin{equation*}
    p_D(t,x,y) \le c\left(1\wedge \frac{\delta_D(x)^{\alpha/2}}{\sqrt{t}}\right)\left(1\wedge \frac{\delta_D(y)^{\alpha/2}}{\sqrt{t}}\right)p(t,x,y).
  \end{equation*}
\end{lemma}

\begin{proof}
  Note that for $t\in (0,1]$, $D_t := t^{-1/\alpha}D=\{t^{-1/\alpha}z: z \in D\}$ is also a $C^{1,1}$ open set with the same characteristics of $D$. Hence, by scaling property (\ref{eq:DHKscale}), \eqref{eq:HKscale} and applying (\ref{eq:pD1Up}) for $D_t$, we obtain for all $x,y \in D$
  \begin{align*}
    p_D(t,x,y) =&~ t^{-d/\alpha}p_{D_t}(1,t^{-1/\alpha}x,t^{-1\alpha}y)\\
    \le&~ ct^{-d/\alpha}\left(1\wedge \delta_{D_t}(t^{-1/\alpha}x)^{\alpha/2}\right)\left(1\wedge \delta_{D_t}(t^{-1/\alpha}y)^{\alpha/2}\right)p(1,t^{-1/\alpha}x,t^{-1/\alpha}y)\\
    =&~c\left(1\wedge \frac{\delta_D(x)^{\alpha/2}}{\sqrt{t}}\right)\left(1\wedge \frac{\delta_D(y)^{\alpha/2}}{\sqrt{t}}\right)p(t,x,y).
  \end{align*}
\end{proof}

\medskip

\subsection{Lower bound estimates} \label{S:4.2}

The following near diagonal lower estimate of $p_B$ for balls $B$ is also called \emph{localized lower estimate} in some literatures (cf. \cite{GrigoryanHuHu.2017.JFA3311,GrigoryanHuHu.2018.AM433}).

\begin{lemma}\label{lem:LLE}
  Let $B :=B(x_0,r)$ be a ball with  radius $r > 0$. For any $a_1 > 0$, there exists $a_2 :=a_2(d,\alpha,a_1) \in (0, 1/a_1)$ such that for any $t^{1/\alpha} \le a_1 r$,
  \begin{equation*}
    p_B(t,x,y) \ge c t^{- {d}/{\alpha}},\qquad x,y \in B(x_0,a_2t^{1/\alpha}),
  \end{equation*}
  where $c=c(d,\alpha,a_1)>0$.
\end{lemma}

\begin{proof}
  \textbf{Step 1.} We first show that for any $x, y\in B$ and $t > 0$,
  \begin{equation}\label{eq:hkBcomp}
    p(t,x,y) \le p_B(t,x,y) + 2\sup_{s\in (t/2,t]} \sup_{\genfrac{}{}{0pt}{4}{z\in \{x, y\}}{w\in B^c}} p(s,z,w).
  \end{equation}
  By the Markov property of $X$, \eqref{eq:pU-def}, the symmetry of $p(t,x,y)$ and $p_B(t,x,y)$ in $x$ and $y$, and Theorem \ref{thm:qHolder}, we have for any $x, y\in B$,
  \begin{align*}
    p(2t,x,y) =&~ \int_{\R^d} p(t,x,z) p(t,z,y) dz \\
    =&~ \int_{\R^d}   p_B(t, x, z)  p(t, z, y) dz  + \int_{\R^d} \E_x \left[ p(t-\tau_B, X_{\tau_B}, z); \tau_B<t \right]  p(t, z, y)  dz  \\
    =&~ \int_{B}   p(t, y, z) p_B(t, z, x) dz    +   \E_x \left[ p(2t-\tau_B, X_{\tau_B}, y); \tau_B<t \right]      \\
    =&~ \int_{B}   p_B(t, y, z) p_B(t, z, x)  dz   + \int_{B} \E_y\left[ p(t-\tau_B, X_{\tau_B}, z); \tau_B<t \right]  p_B(t, z, x)  dz\\
    &~ +\, \E_x \left[ p(2t-\tau_B, X_{\tau_B}, y); \tau_B<t \right] \\
    \leq&~ \int_{B} p_B(t, y, z) p_B(t, z, x)  dz   + \int_{\R^d} \E_y \left[ p(t-\tau_B, X_{\tau_B}, z); \tau_B<t \right]  p (t, z, x)  dz\\
    &~ +  \,  \E_x \left[ p(2t-\tau_B, X_{\tau_B}, y); \tau_B<t \right]      \\
    =&~ p_B(2t, y, x)   +\E_y \left[ p(2t-\tau_B, X_{\tau_B}, x); \tau_B<t \right]   \\
    &~ + \,  \E_x \left[ p(2t-\tau_B, X_{\tau_B}, y); \tau_B<t \right]      \\
    \leq &~  p_B(2t, x, y)+ \sup_{s\in (t, 2t] } \sup_{z\in B^c}   p(s,z,x) + \sup_{s\in (t, 2t] } \sup_{z\in B^c}   p(s,z,y)  .
  \end{align*}
  This establishes the claim  (\ref{eq:hkBcomp}) after replacing $2t$ by $t$.

\medskip

  \textbf{Step 2.} We next show that there exists $a :=a(d,\alpha) \in (0,1)$ and $c_1 =c_1(d,\alpha)>0$ such that for any $t^{1/\alpha} \le a r$,
  \begin{equation}\label{eq:preLLE}
    p_B(t,x,y) \ge c_1 t^{- {d}/{\alpha}},\qquad x,y \in B(x_0,t^{1/\alpha}/2) \subset B.
  \end{equation}
  Indeed, we have by (\ref{eq:hkBcomp}) that for all $x,y \in B(x_0,t^{1/\alpha}/2)$,
  \begin{equation*}
    |x^{(i)} - y^{(i)}| \le |x-y| \le t^{1/\alpha},\quad i=1,2,\cdots,d,
  \end{equation*}
  and then, by (\ref{eq:pEst}),
  \begin{equation}\label{eq:hkB-lower}
    p(t,x,y) \ge C_{\CTstable}^{-1}\prod_{i=1}^d \left(t^{-1/\alpha}\wedge \frac{t}{(t^{1/\alpha})^{1+\alpha}}\right) = C_{\CTstable}^{-1}
    t^{-d/\alpha}.
  \end{equation}
  On the other hand, for all $t^{1/\alpha} \le a r$ (where $a$ is to be determined later), $z \in \{x, y\}$ and $w \in B^c$, we have
  \begin{equation*}
    |z-w| \ge  |x_0-w|-|x_0-z| \ge r - \frac{t^{1/\alpha}}{2}\ge  \frac{r}{2} \ge \frac{t^{1/\alpha}}{2a},
  \end{equation*}
  and then, there exists
   $k\in \{1,2,\cdots, d\}$  so that
  \begin{equation*}
    |z^{(k)} - w^{(k)}| \ge \frac{r}{2\sqrt{d}} \ge \frac{t^{1/\alpha}}{2\sqrt{d}a}.
  \end{equation*}
  Consequently,
  \begin{equation}\label{eq:hkB-upper}
    \begin{split}
      \sup_{s\in[t/2,t)}\sup_{z\in \{x,y\}, w\in B^c} p(s,z,w) \le&~ C_{\CTstable}  \bigg(\prod^d_{\genfrac{}{}{0pt}{4}{i=1}{i\neq k}} s^{-1/\alpha}\bigg)\left(\frac{s}{|z^{(k)} - w^{(k)}|^{1+\alpha}}\right)\\
      \le&~ C_{\CTstable}  2^{(d-1)/\alpha} t^{-(d-1)/\alpha}\left(\frac{t}{(t^{1/\alpha}/(2\sqrt{d} a))^{1+\alpha}}\right)\\
      =&~ c_2 a^{1+\alpha} t^{-d/\alpha},
    \end{split}
  \end{equation}
  where $c_2 := 2^{(d-1)/\alpha +1+\alpha} C_{\CTstable} d^{(1+\alpha)/2} > 0$. Combining (\ref{eq:hkBcomp}), (\ref{eq:hkB-lower}) and (\ref{eq:hkB-upper}), we obtain
  \begin{equation*}
    p_B(t,x,y) \ge (C_{\CTstable}^{-1} - 2c_2 a^{1+\alpha}) t^{-d/\alpha}.
  \end{equation*}
  Setting $a := (4C_{\CTstable}c_2)^{-1/(1+\alpha)}$, we obtain (\ref{eq:preLLE}) with $c_1 = (2C_{\CTstable})^{-1}$.

\medskip

  \textbf{Step 3.} When $a_1 \le a$,  this lemma follows directly from (\ref{eq:preLLE}) with $c = c_1$ and $a_2 = \frac{1}{2}$. So it suffices to consider the case that $a_1 > a$. Let
  \begin{equation*}
    n := \left[\left(\frac{a_1}{a}\right)^\alpha\right] + 1,\quad \text{ and } \quad a_2 :=  \frac{1}{2n^{1/\alpha}}.
  \end{equation*}
  For all $t^{1/\alpha} \le a_1 r$, we have
  \begin{equation*}
    \left(\frac{t}{n}\right)^{1/\alpha} \le \frac{a_1 r}{n^{1/\alpha}} \le \frac{a_1 r}{((a_1/a)^{\alpha})^{1/\alpha}} = a r < r,
  \end{equation*}
  and then,
  \begin{equation*}
    B(x_0, a_2t^{1/\alpha}) =  B(x_0,(t/n)^{1/\alpha}/2) \subset B.
  \end{equation*}
  Hence, by Step 2 and semigroup property of $p_B(t,x,y)$, we have for all $t^{1\alpha} \le a_1 r$ and $x,y \in B(x_0,a_2 t^{1/\alpha})$,
  \begin{align*}
    p_B(t,x,y) =&~ \int_{B^{n-1}} p_B(t/n ,x,z_1)p_B(t/n ,z_1,z_2)\cdots p_B(t/n ,z_{d-1},y) dz_1dz_2\cdots dz_{d-1}\\
    \ge&~ \int_{B(x_0,a_2t^{1/\alpha})^{n-1}} p_B(t/n ,x,z_1)p_B(t/n ,z_1,z_2)\cdots p_B(t/n ,z_{d-1},y) dz_1dz_2\cdots dz_{d-1}\\
    \ge&~ \left(\prod_{i=1}^n c_1 \left(t/n \right)^{-d/\alpha}\right) |B(x_0,a_2t^{1/\alpha})|^{n-1}\\
    =&~ c_1^n n^{dn/\alpha} t^{-dn/\alpha}\cdot |B(0,1)|^{n-1} (a_2 t^{1/\alpha})^{d(n-1)}\\
    =&\!\!: c t^{-d/\alpha},
  \end{align*}
  where $c= c_1^n n^{dn/\alpha} |B(0,1)|^{n-1} a_2^{d(n-1)} > 0$. This completes the proof of the lemma.
\end{proof}

\begin{lemma}\label{lem:pUnle-int}
  Let $U\subset \mathbb{R}^d$ be a non-empty open set. For any $a_1 > 0$, there exists $a_2 :=a_2(d,\alpha,,a_1) > 0$, $c=c(d,\alpha,a_1)>0$ such that for all $t > 0$ and $x,y \in U$ with $\delta_U(x)\wedge \delta_U(y) \ge a_1 t^{1/\alpha}$ and $|x-y| < a_2 t^{1/\alpha}$, we have
  \begin{equation*}
    p_U(t,x,y) \ge c t^{- {d}/{\alpha}}.
  \end{equation*}
\end{lemma}
\begin{proof}
  Fix $t > 0$ and let $r := a_1 t^{1/\alpha}$. Then, we have for $x \in U$ with $\delta_U(x) \ge a_1 t^{1/\alpha} = r$,
  \begin{equation*}
    t^{1/\alpha} \le a_1^{-1} r \quad \text{ and }\quad B:=B(x,r) \subset U.
  \end{equation*}
  By Lemma \ref{lem:LLE} (with $a_1$ being replaced by $a_1^{-1}$), there exists $a_2 :=a_2(d,\alpha,a_1) \in (0,a_1)$ and $c=c(d,\alpha,a_1)>0$ such that
  \begin{equation}\label{eq:LLEB}
    p_B(t,z,w) \ge c t^{-d/\alpha},\quad z,w\in B(x,a_2 t^{1/\alpha}).
  \end{equation}
  On the other hand, by shrinking $a_2$ if necessary, we have for all $x,y \in U$ with $\delta_U(x)\wedge \delta_U(y) \ge a_1 t^{1/\alpha}$ and $|x-y| < a_2 t^{1/\alpha}$,
  \begin{equation*}
    y\in B(x,a_2t^{1/\alpha}) \quad \text{ and }\quad B(x,a_2 t^{1/\alpha}) \subset B(x,r) \subset U.
  \end{equation*}
  Consequently, by (\ref{eq:LLEB}), we have
  \begin{equation*}
    p_U(t,x,y) \ge p_B(t,x,y) \ge c t^{-d/\alpha}.
  \end{equation*}
\end{proof}

\medskip

The following lemma gives the so-called \emph{survival estimate} in some literatures (cf. \cite{GrigoryanHuHu.2017.JFA3311}). Recall that $\tau_{B(x,r)}$ is the first exit time of $X$ for the ball $B(x,r)$, and  $\tau_{Q(x,r)}$ is the first exit time of $X$ for the cube $Q(x,r)$.

\begin{lemma}\label{lem:S}
  For any $t >0$, $r > 0$ and $x \in \mathbb{R}^d$, there exists $c = c(d,\alpha) >0$ such that
  \begin{equation}\label{eq:S-ball}
    \mathbb{P}_x\left(\tau_{B(x,r)} \le t\right) \le\frac{ct}{r^\alpha}
  \end{equation}
  and
  \begin{equation}\label{eq:S-cube}
    \mathbb{P}_x\left(\tau_{Q(x,r)} \le t\right) \le\frac{ct}{r^\alpha}.
  \end{equation}
\end{lemma}

\begin{proof}
  Fix $t >0$, $r > 0$ and $x \in \mathbb{R}^d$. By the strong Markov property (cf. \cite[p. 43-44]{BlumenthalGetoor.1968.313}), one can prove that
  \begin{equation*}
    \mathbb{P}_x\left(\tau_{B(x,r)} \le t\right) \le 2\sup_{s\le t} \sup_{y\in \mathbb{R}^d} \mathbb{P}_y\left(|X_s - X_0|\ge r/2\right)
  \end{equation*}
  (see also \cite[(3.1), p. 2494]{ChenKimSong.2012.AP2483}). By the
  scaling property (\ref{eq:HKscale}), we have
  \begin{align*}
    \mathbb{P}_y\left(|X_s - X_0|\ge r/2\right) =&~ \int_{|y-z|\ge r/2} p(s,y,z)dz= \int_{|z|\ge r/2} p(s,0,z)dz \\
    =&~ s^{-d/\alpha}\int_{|z|\ge r/2} p(1,0,s^{-1/\alpha}z)dz = \int_{|z|\ge r/(2s^{1/\alpha})} p(1,0,z)dz.
  \end{align*}
  Noting that for any $z =(z^{(1)}, z^{(2)}, \cdots, z^{(d)})$ with $|z|\ge r/(2s^{1/\alpha})$, there exists some $1\le k\le d$ such that $|z^{(k)}| \ge r/(2\sqrt{d}s^{1/\alpha})$, and then,
  \begin{equation*}
    \{|z| \ge r/(2s^{1/\alpha})\} \subset \bigcup_{k=1}^d \{|z^{(k)}| \ge r/(2\sqrt{d}s^{1/\alpha})\}.
  \end{equation*}
  Consequently, by (\ref{eq:ppalpha}) and \eqref{e:1.4}
  \begin{align*}
    \mathbb{P}_y\left(|X_s - X_0|\ge r/2\right) \le&~ \sum_{k=1}^d \int_{|z^{(k)}| \ge r/(2\sqrt{d}s^{1/\alpha})} p^{(1,\alpha)} (1,0,z^{(k)}) dz^{(k)}\\
    \le&~ c_1(d,\alpha) \int_{|\theta| \ge r/(2\sqrt{d}s^{1/\alpha})} \frac{d\theta}{|\theta|^{1+\alpha}}\\
    =&~ \frac{(2\sqrt{d})^\alpha c_1}{\alpha} \frac{s}{r^\alpha}.
  \end{align*}
  Combining the above three formulas, we obtain the first inequality (\ref{eq:S-ball}).

  Finally, since $B(x,r) \subset Q(x,r)$, we have $\{\tau_{Q(x,r)} \le t\} \subset \{\tau_{B(x,r)} \le t \}$ for all $t > 0$. Then, the second inequality (\ref{eq:S-cube}) follows from this and (\ref{eq:S-ball}).
\end{proof}

\begin{lemma}\label{lem:pUle-int}
  Assume that $U$ is an open set satisfying condition \hyperref[eq:conH]{$\bf (H_{\gamma})$} for some $\gamma\in (0,1]$. For any $a_1 > 0$, there exists $c=c(d,\alpha, \gamma,a_1)>0$ such that for all $t > 0$ and $x,y \in U$ with $\delta_U(x)\wedge \delta_U(y) \ge a_1 t^{1/\alpha}$, we have
  \begin{equation}\label{eq:pUle-int}
    p_U(t,x,y) \ge c p(t,x,y).
  \end{equation}
\end{lemma}

\begin{proof}
  Fix $t > 0$ and $x,y \in U$ with $\delta_U(x)\wedge \delta_U(y) \ge a_1 t^{1/\alpha}$. By Lemma \ref{lem:pUnle-int}, there exists $a_2 > 0$ and $c_1 = c_1(d,\alpha) > 0$ such that if $|x-y| < a_2 t^{1/\alpha}$, then
  \begin{equation*}
    p_U(t,x,y) \ge c_1 t^{-d/\alpha},
  \end{equation*}
  which together with (\ref{eq:pEst}) yields (\ref{eq:pUle-int}).

  It remains to consider the case when $|x-y| \ge a_2 t^{1/\alpha}$. Without loss of generality, we may assume that $a_2 < a_1$.

\medskip

  \textbf{Step 1.} Let $a \in (0,a_2/2]$ and $\delta \in (0,1)$ be determined later. By semigroup property of $p_U(t,x,y)$, we have
  \begin{align}\label{eq:pUle-int-1}
    p_U(t,x,y) =&~ \int_U p_U(\delta t,x,z) p_U((1-\delta)t,z,y) dz\notag\\
    \ge&~ \int_{B(y,a t^{1/\alpha})}p_U(\delta t,x,z) p_U((1-\delta)t,z,y) dz\\
    \ge&~ \inf_{z\in B(y,a t^{1/\alpha})} p_U((1-\delta)t,z,y) \mathbb{P}_x\left(X_{\delta t}^U \in B(y,a t^{1/\alpha})\right).\notag
  \end{align}
  Note that for all $z \in B(y,a t^{1/\alpha})$,
  \begin{equation*}
    \delta_U(z) \wedge \delta_U(y) \ge (a_1-a)t^{1/\alpha} \ge \frac{a_1}{2}t^{1/\alpha}.
  \end{equation*}
  Hence, by Lemma \ref{lem:pUnle-int}, there exists $a_3 \in (0,a_2)$ such that if
  \begin{equation}\label{eq:adelta-1}
    |z-y| < a t^{1/\alpha} < a_3((1-\delta)t)^{1/\alpha},
  \end{equation}
  then,
  \begin{equation*}
    \inf_{z\in B(y,a t^{1/\alpha})} p_U((1-\delta)t,z,y) \ge c_2 ((1-\delta)t)^{-d/\alpha} =  c_2 (1-\delta)^{-d/\alpha} t^{-d/\alpha},
  \end{equation*}
  for some $c_2 = c_2(d,\alpha,a_1) > 0$. Consequently, under condition (\ref{eq:adelta-1}), the above inequality together with (\ref{eq:pUle-int-1}) implies
  \begin{equation}\label{eq:pUle-int-2}
    p_U(t,x,y) \ge c_2 (1-\delta)^{-d/\alpha} t^{-d/\alpha} \mathbb{P}_x\left(X_{\delta t}^U \in B(y,a t^{1/\alpha})\right).
  \end{equation}
  Note that, inequality (\ref{eq:adelta-1}) can be achieved by choosing $a,\delta$ small enough such that
  \begin{equation}\label{eq:adelta-2}
    a < a_3 (1-\delta)^{1/\alpha}.
  \end{equation}

\medskip

  \textbf{Step 2} We  next derive a lower bound of $\mathbb{P}_x\left(X_{\delta t}^U \in B(y,a t^{1/\alpha})\right)$. Since $U$ satisfies condition \hyperref[eq:conH]{$\bf (H_{\gamma})$} and $\delta_U(x)\wedge \delta_U(y) \ge a_1 t^{1/\alpha}$, there exists a permutation $\{i_1,i_2,\cdots,i_d\}$ of $\{1,2,\cdots,d\}$, such that
  \begin{equation}\label{eq:pUnle-dBalls}
    B(\overline{xy}_k, \gamma a_1 t^{1/\alpha}) \subset U,\qquad k=1,2,\cdots,d,
  \end{equation}
  where $\overline{xy}_1 := [x]_{y^{(i_1)}}^{i_1}, \, \overline{xy}_2 := [\overline{xy}_1]_{y^{(i_2)}}^{i_2}, \cdots, \overline{xy}_d := [\overline{xy}_{d-1}]_{y^{(i_d)}}^{i_d}=y$. Set
  \begin{equation*}
    r := \frac{a t^{1/\alpha}}{\sqrt{d}},
  \end{equation*}
  where $a$ is chosen to be small enough such that
  \begin{equation}\label{eq:adelta-3}
    a < \frac{\gamma a_1}{2}.
  \end{equation}
  Then $Q_0:=Q(x,r) \subset U$ and
  \begin{equation}\label{eq:pUnle-Qks}
    Q_k := Q(\overline{xy}_k, r) \subset B_k := B(\overline{xy}_k, \sqrt{d}r) \subset U,\qquad k= 1,2,\cdots, d.
  \end{equation}

  In the rest of the proof, for a number $\lambda>0$ and a cube $Q:= Q(z,r)$, we use the notation $\lambda Q$ to denote the set $Q(z,\lambda r)$, that is $\lambda Q(z,r) = Q(z,\lambda r)$.

  By semigroup property, we have
  \begin{align}\label{eq:pUle-int-3}
    &~\mathbb{P}_x\left(X_{\delta t}^U \in B(y,a t^{1/\alpha})\right)\notag\\
    = &~ \int_{B(y,a t^{1/\alpha})} p_U(\delta t,x,z_d) dz_d \notag\\
    = &~ \int_{B(y,a t^{1/\alpha})} \left(\int_{U^{d-1}} p_U(\delta t/d,x,z_1) p_U(\delta t/d,z_1,z_2) \cdots p_U(\delta t/d,z_{d-1},z_d)dz_1 \cdots dz_{d-1}\right) dz_d  \notag\\
    \ge &~ \int_{Q_d} \int_{\left(2^{1-d}Q_1\right)\times\left(2^{2-d}Q_2\right) \times \cdots \times \left(2^{-1}Q_{d-1}\right)}  \notag\\
    &~\qquad p_U(\delta t/d,x,z_1) p_U(\delta t/d,z_1,z_2) \cdots p_U(\delta t/d,z_{d-1},z_d)dz_1\cdots dz_{d-1}dz_d  \notag\\
    \ge&~ \mathbb{P}_x\left(X_{\delta t/d}^U \in 2^{1-d}Q_1\right) \inf_{z_1\in 2^{1-d}Q_1} \mathbb{P}_{z_1}\left(X_{\delta t/d}^U \in 2^{2-d}Q_2\right)  \notag\\
    &~\times\cdots \times\inf_{z_{d-2}\in 2^{-2}Q_{d-2}} \mathbb{P}_{z_{d-2}} \left(X_{\delta t/d}^U \in 2^{-1}Q_{d-1}\right) \inf_{z_{d-1}\in 2^{-1}Q_{d-1}} \mathbb{P}_{z_{d-1}} \left(X_{\delta t/d}^U \in Q_d\right) \notag\\
    \ge&~ \inf_{z_0\in 2^{-d} Q_0}\mathbb{P}_{z_0} \left(X_{\delta t/d}^U \in 2^{1-d}Q_1\right) \inf_{z_1\in2^{1-d}Q_1} \mathbb{P}_{z_1} \left(X_{\delta t/d}^U \in 2^{2-d}Q_2\right) \notag\\
    &~ \times\cdots \times\inf_{z_{d-2}\in 2^{-2}Q_{d-2}} \mathbb{P}_{z_{d-2}} \left(X_{\delta t/d}^U \in 2^{-1}Q_{d-1}\right) \inf_{z_{d-1}\in 2^{-1}Q_{d-1}} \mathbb{P}_{z_{d-1}} \left(X_{\delta t/d}^U \in Q_d\right).
  \end{align}

\medskip

  \textbf{Step 3.}  We estimate the lower bound of each term on the right hand side of (\ref{eq:pUle-int-3}). In fact, they can be estimated similarly. We claim that for each $1\le k\le d$,
  \begin{equation}\label{eq:pUle-int-4}
    \inf_{z\in  2^{ k-1-d}Q_{k-1}} \mathbb{P}_z\left(X_{\delta t/d}^U \in 2^{k-d}Q_k\right) \ge c_3\left(1\wedge \frac{t^{1+1/\alpha}}{|y^{(i_k)} - x^{(i_k)}|^{1+\alpha}}\right),
  \end{equation}
  where the constant $c_3 > 0$ is independent of $t,x,y$.

  Fix $z \in 2^{(k-1)-d} Q_{k-1}$. Let $\widehat{Q} := Q(z,r/2^{d-(k-1)})$ and
  \begin{align*}
    \widetilde{Q} :=&~ \widehat{Q} + (y^{(i_k)} - x^{(i_k)}) e_{i_k} \\
    =&~ \left\{w \in \mathbb{R}^d:|w^{(i)} - z^{(i)}| < \frac{r}{2^{d-(k-1)}} \text{ for } i \neq i_k; |w^{(i_k)} - (z^{(i_k)}+y^{(i_k)} - x^{(i_k)})| < \frac{r}{2^{d-(k-1)}}\right\}.
  \end{align*}
  Then, we have $\widehat{Q} \subset 2^{k-d}Q_{k-1}$ and $\widetilde{Q} \subset 2^{k-d}Q_k$. Indeed, since $z \in 2^{(k-1)-d} Q_{k-1}$ and
  \begin{equation*}
    2^{k-d}Q_k = 2^{k-d}Q_{k-1} + (y^{(i_k)} - x^{(i_k)}) e_{i_k},
  \end{equation*}
  we have, for any $w \in \widetilde{Q}$ and $i \neq i_k$,
  \begin{align*}
    |w^{(i)} - \overline{xy}_k^{(i)}| \le&~ |w^{(i)} - z^{(i)}| + |z^{(i)} - \overline{xy}_k^{(i)}|\\
    =&~ |w^{(i)} - z^{(i)}| + |z^{(i)} - \overline{xy}_{k-1}^{(i)}|\\
    <&~ \frac{r}{2^{d-(k-1)}} + \frac{r}{2^{d-(k-1)}} = \frac{r}{2^{d-k}},
  \end{align*}
  and for $i_k$,
  \begin{align*}
    |w^{(i_k)} - \overline{xy}_k^{(i_k)}| =&~ |w^{(i_k)} - y^{(i_k)}|\\
    \le&~ |w^{(i_k)} - (z^{(i_k)}+y^{(i_k)} - x^{(i_k)})| + |z^{(i_k)} - x^{(i_k)}|\\
    =&~ |w^{(i_k)} - (z^{(i_k)}+y^{(i_k)} - x^{(i_k)})| + |z^{(i_k)} - \overline{xy}_{k-1}^{(i_k)}|\\
    <&~ \frac{r}{2^{d-(k-1)}} + \frac{r}{2^{d-(k-1)}} = \frac{r}{2^{d-k}}.
  \end{align*}
  This shows that $\widetilde{Q} \subset 2^{k-d}Q_k$.

  Since $\delta_U(z) \wedge \delta_U(w) \ge \frac{\gamma a_1}{2} t^{1/\alpha}$ for all $w \in \widetilde{Q}$ by  (\ref{eq:pUnle-dBalls}), (\ref{eq:pUnle-Qks}) and  (\ref{eq:adelta-3}),  it follows from Lemma \ref{lem:pUnle-int} that there exists $a_4 \in (0,\frac{1}{2})$ such that for all $|z-w| < a_4(\delta t/d)^{1/\alpha}$,
  \begin{equation}\label{eq:pUle_int-temp}
    p_U(\delta t/ d,z,w) \ge c_4 \left(\delta t/ d\right)^{-d/\alpha}.
  \end{equation}
  Let $c_5$ be the constant in (\ref{eq:S-cube}), and choose $a,\delta$ small enough such that
  \begin{equation}\label{eq:adelta-4}
    (c_5\vee1)\frac{\delta}{d} < \left(\frac{a}{2^{d+2}\sqrt{d}}\right)^\alpha.
  \end{equation}
  Let
  \begin{equation*}
    a_5 := a_4 \left( {\delta}/{d}\right)^{1/\alpha}.
  \end{equation*}

  \textbf{Case 1.} $|y^{(i_k)} - x^{(i_k)}| \le a_5 t^{1/\alpha}$. In this case, by (\ref{eq:adelta-4}),
  \begin{equation*}
    |y^{(i_k)} - x^{(i_k)}| \le a_5 t^{1/\alpha} = a_4 \left(\frac{\delta}{d}\right)^{1/\alpha}t^{1/\alpha} \le \frac{1}{2}\cdot\frac{a t^{1/\alpha}}{2^{d+2}\sqrt{d}} < \frac{r}{2^{d-k+2}}.
  \end{equation*}
  Hence, by the definitions of $\widehat{Q}$ and $\widetilde{Q}$, we have
  \begin{equation*}
    B(z,a_5t^{1/\alpha}) \subset 2^{-1}\widehat{Q} \subset \widehat{Q}\cap \widetilde{Q}
  \end{equation*}
  (See Figure \ref{fig:shift-small}).
  \begin{figure}[htp]
    \begin{tikzpicture}[
      declare function = {
          a = 0;
          b = 0;
          sr = 1.5;
          z1 = a+sr/1.2;
          z2 = b+sr/1.6;
          shiftx = sr/3;
          br = 2.3*sr/5;
          },
        scale = 1.5
      ]
      \draw [style = dashed,blue] (a-sr,b-sr) rectangle (a+sr,b+sr);
      \fill (a,b) circle (1pt);
      \node [below] at (a,b) {\large $\overline{xy}_{k-1}\quad$};
      \node [above left] at (a-sr,b-sr) {\large $2^{k-1-d}Q_{k-1}$};

      \draw [style = dashed,red] (a+shiftx-sr,b-sr) rectangle (a+shiftx+sr,b+sr);
      \fill (a+shiftx,b) circle (1pt);
      \node [below] at (a+shiftx,b) {\large $\overline{xy}_k$};
      \node [above right] at (a+shiftx+sr,b-sr) {\large $2^{k-1-d}Q_k$};

      \draw [blue] (z1-sr,z2-sr) rectangle (z1+sr,z2+sr);
      \fill (z1,z2) circle (1pt);
      \node [below] at (z1,z2) {\large $z$};
      \node [left] at (z1-sr,z2+sr) {\large $\widehat{Q}$};

      \draw [red](z1+shiftx-sr,z2-sr) rectangle (z1+shiftx+sr,z2+sr);
      \fill (z1+shiftx,z2) circle (1pt);
      \node [right] at (z1+shiftx+sr,z2+sr) {\large $\widetilde{Q}$};

      \draw (z1,z2) circle (br);
      \draw [->] (z1 + 0.1, z2) -- (z1+shiftx-0.1,z2);
      \draw [->] (a + 0.1, b) -- (a+shiftx-0.1,b);

      \node [above] at (z1,z2+br) {\large $B(z,a_5t^{1/\alpha})$};
    \end{tikzpicture}
    \caption{The case when $|y^{(i_k)} - x^{(i_k)}| \le a_5 t^{1/\alpha}$}
    \label{fig:shift-small}
  \end{figure}
  ~\\
  Consequently, by (\ref{eq:pUle_int-temp}), we have
  \begin{align*}
    \mathbb{P}_z\left(X_{\delta t/d}^U \in 2^{k-d}Q_k\right) \ge&~ \mathbb{P}_z\left(X_{\delta t/d}^U \in \widetilde{Q}\right)\\
    =&~ \int_{\widetilde{Q}} p_U(\delta t/ d,z,w) dw\\
    \ge&~ \int_{B(z,a_5t^{1/\alpha})} p_U(\delta t/ d,z,w) dw\\
    \ge&~ c_4 \left(\delta t/ d\right)^{-d/\alpha} |B(z,a_5t^{1/\alpha})|\\
    \ge&~ c_4 \left(\delta t/ d\right)^{-d/\alpha} \left(a_5t^{1/\alpha}\right)^d|B(0,1)|\\
    =&~ c_4 \left({\delta}/{d}\right)^{-d/\alpha} a_5^d|B(0,1)|\\
    \ge&~ c_4 \left(\frac{\delta}{d}\right)^{-d/\alpha} a_5^d|B(0,1)| \left(1\wedge \frac{t^{1+1/\alpha}}{|y^{(i_k)} - x^{(i_k)}|^{1+\alpha}}\right),
  \end{align*}
  which is (\ref{eq:pUle-int-4}).

  \textbf{Case 2.} $|y^{(i_k)} - x^{(i_k)}| > a_5 t^{1/\alpha}$. Without loss of generality, we may and do assume that $y^{(i_k)} > x^{(i_k)}$. In this case, let \begin{equation*}
    A := \left\{\theta \in \mathbb{R}: \left(\frac{r}{2^{d-k+2}}-|y^{(i_k)} - x^{(i_k)}|\right)\vee \left(-\frac{r}{2^{d-k+2}}\right) < \theta - (z^{(i_k)}+y^{(i_k)} - x^{(i_k)}) < \frac{r}{2^{d-k+2}}\right\},
  \end{equation*}
  and, for any $w \in 2^{-1}\widetilde{Q}\setminus\overline{\left(2^{-1}\widehat{Q}\right)}$, we have
  \begin{equation}\label{eq:wrange}
    2^{-1}\widetilde{Q}\setminus\overline{\left(2^{-1}\widehat{Q}\right)} = \left\{w \in \mathbb{R}^d: |w^{(i)} - z^{(i)}| < \frac{r}{2^{d-k+2}} \text{ for } i \neq i_k; w^{(i_k)} \in A\right\}
  \end{equation}
  (See Figure \ref{fig:shift-big}).
  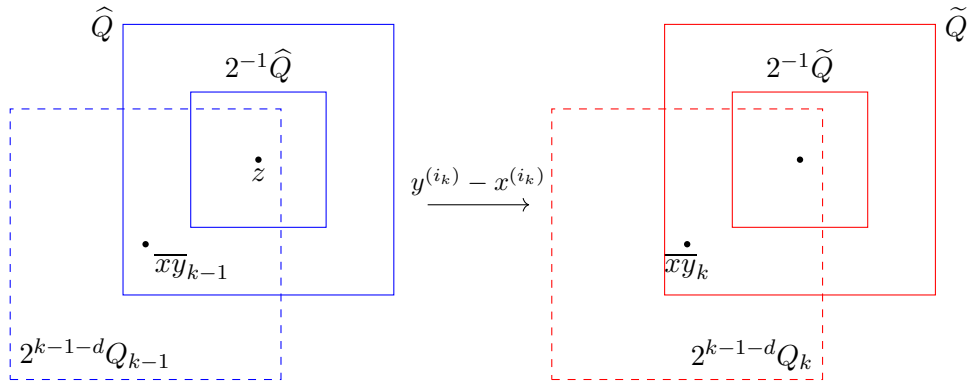
\begin{figure}[htp]
    \begin{tikzpicture}[
      declare function = {
          a = 0;
          b = 0;
          sr = 1.5;
          z1 = a+sr/1.2;
          z2 = b+sr/1.6;
          shiftx = 4*sr;
          br = 2.5*sr/5;
          },
        scale = 1.2
      ]
      \draw [style = dashed,blue] (a-sr,b-sr) rectangle (a+sr,b+sr);
      \fill (a,b) circle (1pt);
      \node [below right] at (a,b) {\large $\overline{xy}_{k-1}$};
      \node [above right] at (a-sr,b-sr) {\large $2^{k-1-d}Q_{k-1}$};

      \draw [style = dashed,red] (a+shiftx-sr,b-sr) rectangle (a+shiftx+sr,b+sr);
      \fill (a+shiftx,b) circle (1pt);
      \node [below] at (a+shiftx,b) {\large $\overline{xy}_k$};
      \node [above left] at (a+shiftx+sr,b-sr) {\large $2^{k-1-d}Q_k$};

      \draw [blue] (z1-sr,z2-sr) rectangle (z1+sr,z2+sr);
      \fill (z1,z2) circle (1pt);
      \node [below] at (z1,z2) {\large $z$};
      \node [left] at (z1-sr,z2+sr) {\large $\widehat{Q}$};

      \draw [red](z1+shiftx-sr,z2-sr) rectangle (z1+shiftx+sr,z2+sr);
      \fill (z1+shiftx,z2) circle (1pt);
      \node [right] at (z1+shiftx+sr,z2+sr) {\large $\widetilde{Q}$};

      \draw [blue] (z1-sr/2,z2-sr/2) rectangle (z1+sr/2,z2+sr/2);
      \node [above] at (z1,z2+sr/2) {\large $2^{-1}\widehat{Q}$};

      \draw [red] (z1+shiftx-sr/2,z2-sr/2) rectangle (z1+shiftx+sr/2,z2+sr/2);
      \node [above] at (z1+shiftx,z2+sr/2) {\large $2^{-1}\widetilde{Q}$};

      \draw [->] (z1 + shiftx/3.2, z2-sr/3) -- (z1+shiftx/2,z2-sr/3) node[pos =0.5,sloped,above] {$y^{(i_k)} - x^{(i_k)}$};
    \end{tikzpicture}
    \caption{The case when $|y^{(i_k)} - x^{(i_k)}| > a_5 t^{1/\alpha}$}
    \label{fig:shift-big}
  \end{figure}

  We are to apply (\ref{eq:XLS}) to estimate $\mathbb{P}_z(X_{\delta t/d}^U \in 2^{k-d}Q_k)$. Indeed, by (\ref{eq:S-cube}) and (\ref{eq:adelta-4}), we have
  \begin{equation}\label{eq:S-cube-temp}
    \mathbb{P}_z \left(\tau_{2^{-1}\widehat{Q}} \le \delta t/ d\right) \le \frac{c_5\delta t}{d}\frac{1}{(r/2^{d-k+2})^\alpha} < \left(\frac{a}{2^{d+2}\sqrt{d}}\right)^\alpha \cdot t\cdot \frac{2^{(d-k+2)\alpha}\sqrt{d}^\alpha}{a^\alpha t} < \frac{1}{2^\alpha} <1.
  \end{equation}
  This implies that
  \begin{equation}\label{eq:pUnle-exitlow}
    \mathbb{E}_z\left[\frac{\delta t}{2d}\wedge \tau_{2^{-1}\widehat{Q}}\right]\ge \frac{\delta t}{2d}\mathbb{P}_z\left(\tau_{2^{-1}\widehat{Q}}> \frac{\delta t}{2d}\right) \ge\frac{(1-2^{-\alpha})\delta t}{2d}.
  \end{equation}
  Denote by $\sigma^U_{2^{-1}\widetilde{Q}}$ the first hitting time of $X^U$ for the set
    $\tfrac{1}{2}\widetilde{Q}$:
  \begin{equation*}
    \sigma^U_{2^{-1}\widetilde{Q}} := \inf\{s>0:X^U_s \in 2^{-1}\widetilde{Q}\}.
  \end{equation*}
  Since $(2^{-1}\widehat{Q}) \cup (2^{-1}\widetilde{Q}) \subset U$, the above inequality together with (\ref{eq:XLS}) yields that
  \begin{align*}
    &~\mathbb{P}_z\left(\sigma^U_{2^{-1}\widetilde{Q}} < \delta t/ d\right) \\
    \ge &~ \mathbb{P}_z\left(X^U_{\left(\frac{\delta t}{2d}\right)\wedge\tau^U_{2^{-1}\widehat{Q}}} \in 2^{-1}\widetilde{Q}\setminus\overline{\left(2^{-1}\widehat{Q}\right)}\right)\\
    =&~ \mathbb{P}_z\left(X_{\left(\frac{\delta t}{2d}\right)\wedge\tau_{2^{-1}\widehat{Q}}} \in 2^{-1}\widetilde{Q}\setminus\overline{\left(2^{-1}\widehat{Q}\right)}\right)\\
    =&~\mathbb{E}_z\!\!\left[\int_0^{\left(\frac{\delta t}{2d}\right)\wedge\tau_{2^{-1}\widehat{Q}}} \sum_{i=1}^d \int_{\mathbb{R}^d} \!\!\! \1_{2^{-1}\widetilde{Q}\setminus\overline{\left(2^{-1}\widehat{Q}\right)}}(X^{(1)}_s,X^{(2)}_s, \cdots,X^{(i-1)}_s,w^{(i)},X^{(i+1)}_s,\cdots, X^{(d)}_s) j(X^{(i)}_s,w^{(i)})dw^{(i)}ds\right]\\
    \ge&~\mathbb{E}_z\!\!\left[\int_0^{\left(\frac{\delta t}{2d}\right)\wedge\tau_{2^{-1}\widehat{Q}}} \!\!\!\int_{\mathbb{R}} \!\!\! \1_{2^{-1}\widetilde{Q}\setminus\overline{\left(2^{-1}\widehat{Q}\right)}}(X^{(1)}_s,X^{(2)}_s, \cdots,X^{(i_k-1)}_s,w^{(i_k)},X^{(i_k+1)}_s,\cdots, X^{(d)}_s) j(X^{(i_k)}_s,w^{(i_k)})dw^{(i_k)}ds\right].
  \end{align*}
  For any $i =1,2,\cdots,d$ and $s < \left(\frac{\delta t}{2d}\right)\wedge\tau_{2^{-1}\widehat{Q}}$, we have
  \begin{equation*}
    |X^{(i)}_s - z^{(i)}| < \frac{r}{2^{d-k+2}}.
  \end{equation*}
  And for any $|w^{(i_k)} -(z^{(i_k)} + y^{(i_k)} - x^{(i_k)})| < \frac{r}{2^{d-k+2}}$, we have
  \begin{align*}
    |X^{(i_k)}_s - w^{(i_k)}| \le&~ |X^{(i_k)}_s - z^{(i_k)}| + |w^{(i_k)} -(z^{(i_k)} + y^{(i_k)} - x^{(i_k)})| + |y^{(i_k)} - x^{(i_k)}|\\
    \le&~ \frac{r}{2^{d-k+2}}+\frac{r}{2^{d-k+2}} +  |y^{(i_k)} - x^{(i_k)}|\\
    \le&~ a t^{1/\alpha} + |y^{(i_k)} - x^{(i_k)}|\\
    \le&~ (a/a_5 + 1)|y^{(i_k)} - x^{(i_k)}|.
  \end{align*}
  Combining the above three formulas, (\ref{eq:pUnle-exitlow}) and (\ref{eq:wrange}), we obtain
  \begin{align*}
    &~ \mathbb{P}_z\left(\sigma^U_{2^{-1}\widetilde{Q}} < \delta t/ d\right) \\
    \ge&~ \mathbb{E}_z \left[\int_0^{\left(\frac{\delta t}{2d}\right) \wedge\tau_{2^{-1}\widehat{Q}}} \!\!\int_{\mathbb{R}} \!\! \1_{2^{-1}\widetilde{Q}\setminus\overline{\left(2^{-1}\widehat{Q}\right)}}(X^{(1)}_s, \cdots, X^{(i_k-1)}_s,w^{(i_k)},X^{(i_k+1)}_s,\cdots, X^{(d)}_s) \cdot j(X^{(i_k)}_s,w^{(i_k)})dw^{(i_k)}ds\right]\\
    \ge &~\mathbb{E}_z\left[\int_0^{\left(\frac{\delta t}{2d}\right)\wedge\tau_{2^{-1}\widehat{Q}}} \int_{\{w^{(i_k)}\in A\}} j(X^{(i_k)}_s,w^{(i_k)})dw^{(i_k)}ds\right]\\
    \ge &~\mathbb{E}_z\left[\int_0^{\left(\frac{\delta t}{2d}\right)\wedge\tau_{2^{-1}\widehat{Q}}} \int_{\{w^{(i_k)}\in A\}} dw^{(i_k)}\frac{\sC_{1,\alpha}}{(a/a_5 + 1)^{1+\alpha}|y^{(i_k)} - x^{(i_k)}|^{1+\alpha}}ds\right]\\
    =&~\mathbb{E}_z\left[\left(\frac{\delta t}{2d}\right)\wedge\tau_{2^{-1}\widehat{Q}}\right] \frac{\sC_{1,\alpha}}{(a/a_5 + 1)^{1+\alpha}|y^{(i_k)} - x^{(i_k)}|^{1+\alpha}}\cdot |\{w^{(i_k)}\in A\}| \\
    =&~\mathbb{E}_z\left[\left(\frac{\delta t}{2d}\right)\wedge\tau_{2^{-1}\widehat{Q}}\right]\frac{\sC_{1,\alpha}}{(a/a_5 + 1)^{1+\alpha}|y^{(i_k)} - x^{(i_k)}|^{1+\alpha}}\\
    &~\times\left(\frac{r}{2^{d-k+2}}-\left(\frac{r}{2^{d-k+2}}-|y^{(i_k)} - x^{(i_k)}|\right)\vee \left(-\frac{r}{2^{d-k+2}}\right)\right)\quad \text{(by (\ref{eq:wrange}))}\\
    \ge&~\frac{(1-2^{-\alpha})\delta t}{2d} \frac{\sC_{1,\alpha}}{(a/a_5 + 1)^{1+\alpha}|y^{(i_k)} - x^{(i_k)}|^{1+\alpha}} \left(|y^{(i_k)} - x^{(i_k)}| \wedge \frac{r}{2^{d-k+1}}\right)\quad  \text{(by (\ref{eq:pUnle-exitlow}))}\\
    \ge&~\frac{(1-2^{-\alpha})\delta t}{2d} \frac{\sC_{1,\alpha}}{(a/a_5 + 1)^{1+\alpha}|y^{(i_k)} - x^{(i_k)}|^{1+\alpha}} \left(a_5 t^{1/\alpha} \wedge  \frac{a t^{1/\alpha}}{2^d\sqrt{d}}\right)\\
    \ge&~ \frac{c_6 t^{1+1/\alpha}}{|y^{(i_k)} - x^{(i_k)}|^{1+\alpha}}.
  \end{align*}
  Furthermore, by the strong Markov property, the above inequality and  an inequality similar to (\ref{eq:S-cube-temp}) yields that
  \begin{align*}
    \mathbb{P}_z\left(X_{\delta t/d}^U \in 2^{k-d}Q_k\right) \ge&~ \mathbb{P}_z\left(X_{\delta t/d}^U \in \widetilde{Q}\right)\\
    \ge&~ \mathbb{P}_z\left(X^U \text{ hits } 2^{-1}\widetilde{Q} \text{ before time } \delta t/d \text{ and stays in } \widetilde{Q} \text{ for at least } \delta t/d \text{ units of time}\right)\\
    \ge&~ \mathbb{P}_z\left(\sigma^U_{2^{-1}\widetilde{Q}} < \delta t/d;~\tau_{\widetilde{Q}}^U\circ \theta_{\sigma^U_{2^{-1}\widetilde{Q}}} > \delta t/d\right)\\
    =&~ \mathbb{P}_z\left(\sigma^U_{2^{-1}\widetilde{Q}} < \delta t/d;~ \mathbb{E}_{X_{\sigma^U_{2^{-1}\widetilde{Q}}}} \left[\tau_{\widetilde{Q}}^U> \delta t/d\right]\right)\\
    =&~ \mathbb{P}_z\left(\sigma^U_{2^{-1}\widetilde{Q}} < \delta t/d\right) \inf_{w \in 2^{-1}\widetilde{Q}} \mathbb{P}_w \left(\tau_{\widetilde{Q}}^U> \delta t/d\right)\\
    =&~ \mathbb{P}_z\left(\sigma^U_{2^{-1}\widetilde{Q}} < \delta t/d\right) \inf_{w \in 2^{-1}\widetilde{Q}} \mathbb{P}_w \left(\tau_{\widetilde{Q}} > \delta t/d\right)\\
    =&~ \mathbb{P}_z\left(\sigma^U_{2^{-1}\widetilde{Q}} < \delta t/d\right) \inf_{w \in 2^{-1}\widetilde{Q}} \mathbb{P}_w \left(\tau_{Q(w,r/2^{d-k+2})} > \delta t/d\right)\\
    \ge&~ \frac{c_6 t^{1+1/\alpha}}{|y^{(i_k)} - x^{(i_k)}|^{1+\alpha}} \cdot \left(1-\frac{1}{2^\alpha}\right)  \quad
    \text{(by an inequality similar to (\ref{eq:S-cube-temp}))},
  \end{align*}
  which gives (\ref{eq:pUle-int-4}).

\medskip

  \textbf{Step 4.} Note that $\{i_1,i_2,\cdots,i_d\}$ is a permutation of $\{1,2,\cdots,d\}$. By choosing $a,\delta$ small enough such that all the conditions (\ref{eq:adelta-2}), (\ref{eq:adelta-3}), (\ref{eq:adelta-4}) are satisfied, it follows from (\ref{eq:pUle-int-2}), (\ref{eq:pUle-int-3}) and (\ref{eq:pUle-int-4}) that
  \begin{align*}
    p_U(t,x,y) \ge&~ c_2 (1-\delta)^{-d/\alpha} t^{-d/\alpha} \prod_{k=1}^d c_3 \left(1\wedge \frac{t^{1+1/\alpha}}{|y^{(i_k)} - x^{(i_k)}|^{1+\alpha}}\right)\\
    =&~ c_2c_3^d (1-\delta)^{-d/\alpha} \prod_{k=1}^d \left(t^{-1/\alpha}\wedge \frac{t}{|y^{(k)} - x^{(k)}|^{1+\alpha}}\right),
  \end{align*}
  which together with (\ref{eq:pEst}) finishes the proof.
\end{proof}

\begin{rem}\label{rem:para}  \rm
  Suppose that $t > 0$, $a \in \mathbb{R}$ and $U\subset \mathbb{R}^d$ is  an open set. Let $Q_1:=Q(x,c_1t^{1/\alpha})$ and $Q_2:=Q(x+ae_i,c_1t^{1/\alpha})$ be two cubes with $Q_1\cup Q_2 \subset U$, where $1\le i\le d$, $c_1>0$ is a small constant and $e_i$ is the unit vector in the positive $x^{(i)}$-direction. Then, by the same arguments that lead to \eqref{eq:pUle-int-4}, we can in fact prove the following more general inequality: there exists $c_2 > 0$ independent of $t,a,x$ such that
  \begin{equation*}
    \inf_{z\in2^{-1}Q_1} \mathbb{P}_z\left(X_{t}^U \in Q_2\right) \ge c_2\left(1\wedge \frac{t^{1+1/\alpha}}{|a|^{1+\alpha}}\right).
  \end{equation*}
\end{rem}

\begin{lemma}\label{lem:pDTle}
  Let $D\subset \mathbb{R}^d$ be a $C^{1,1}$ open set with characteristics $(R,\Lambda)$ satisfying condition \hyperref[eq:conH]{$\bf (H_{\gamma})$} for some $\gamma \in (0,1]$. There exist constants $c=c(d,\alpha,R,\Lambda,\gamma) > 0$ and $t_* = t_*(d,\alpha,R,\Lambda,\gamma) > 0$ such that for all $x,y \in D$, we have
  \begin{equation*}
    p_D(3t_*,x,y) \ge c\left(1\wedge \frac{\delta_D(x)^{\alpha/2}}{\sqrt{t_*}}\right)\left(1\wedge \frac{\delta_D(y)^{\alpha/2}}{\sqrt{t_*}}\right) p(t_*,x,y).
  \end{equation*}
\end{lemma}

\begin{proof}
  Let $\delta_0,r_0$ be the constants from Lemma \ref{lem:exitDistri}. For any $z \in D$, we can choose $Q_z \in \partial D$ such that $|z-Q_z| = \delta_D(z)$. Let $D_{Q_z}(\delta_0,r_0)$ be the set from Lemma \ref{lem:exitDistri}. Define
  \begin{equation*}
    \delta_1 := \frac{\delta_0}{2\sqrt{1+\Lambda^2}},\quad \text{ and } \quad r_1 := \frac{r_0}{\sqrt{1+\Lambda^2}}.
  \end{equation*}
  Note that
  \begin{equation*}
    \delta_0 < \frac{r_0}{2\sqrt{1+\Lambda^2}} = \frac{r_1}{2} < \frac{r_0}{2}.
  \end{equation*}
  For the above $z$, we use the following notation:
  \begin{align*}
    &~ U^0_z := D_{Q_z}(r_1,r_0) \setminus D_{Q_z}(\delta_0,r_0),\\
    &~ U_z:= \bigcup_{u\in U^0_z} B(u,\delta_1),\\
    &~ E_z:= \begin{cases}
      B(z,\delta_1), &z \notin D_{Q_z}(\delta_0,r_0)~ \ (\text{ i.e. }\delta_D(z) \ge \delta_0),\\
      U_z,& z \in D_{Q_z}(\delta_0,r_0)~\  (\text{ i.e. }\delta_D(z) < \delta_0).
    \end{cases}
  \end{align*}
  Let $C_{\CTexitL}$, $C_{\CTexitU}$ be the constants from Lemma \ref{lem:exitDistri}, and define
  \begin{equation}\label{eq:def-tast}
    t_* = \frac{2C_{\CTexitU}}{C_{\CTexitL}} \vee 1.
  \end{equation}

  By semigroup property, we have for any $x,y\in D$,
  \begin{equation}\label{eq:pDle-1}
  \begin{split}
    p_D(3t_*,x,y) =&~ \int_{D\times D} p_D(t_*,x,u) p_D(t_*,u,v) p_D(t_*,v,y) dudv\\
    \ge&~ \int_{E_x \times E_y} p_D(t_*,x,u) p_D(t_*,u,v) p_D(t_*,v,y) dudv\\
    \ge&~ \inf_{u\in E_x,v\in E_y} p_D(t_*,u,v) \int_{E_x} p_D(t_*,x,u) du \int_{E_y} p_D(t_*,v,y) dv\\
    =&~ I_1 \cdot I_2 \cdot I_3.
  \end{split}
  \end{equation}
  We estimate $I_1$, $I_2$, $I_3$ separately. For $1\le i\le d$ and $(u,v) \in E_x\times E_y$, if $|y^{(i)} - x^{(i)}| < 3r_0$, then,
  \begin{equation*}
    |u^{(i)} - v^{(i)}| \le |u^{(i)} - x^{(i)}| + |x^{(i)} - y^{(i)}| + |y^{(i)} - v^{(i)}| < 3r_0 + 3r_0 + 3r_0 = 9r_0,
  \end{equation*}
  and
  \begin{equation*}
    t_*^{-1/\alpha}\wedge \frac{t_*}{|u^{(i)} - v^{(i)}|^{1+\alpha}} \ge t_*^{-1/\alpha}\wedge \frac{t_*}{(9r_0)^{1+\alpha}} \ge \left(1\wedge \frac{t_*^{1+1/\alpha}}{(9r_0)^{1+\alpha}}\right)\left(t_*^{-1/\alpha}\wedge \frac{t_*}{|x^{(i)} - y^{(i)}|^{1+\alpha}}\right).
  \end{equation*}
  If $|y^{(i)} - x^{(i)}| \ge 3r_0$, then,
  \begin{equation*}
    |u^{(i)} - v^{(i)}| \le |u^{(i)} - x^{(i)}| + |x^{(i)} - y^{(i)}| + |y^{(i)} - v^{(i)}| < 3r_0 + |x^{(i)} - y^{(i)}| + 3r_0 \le 3|x^{(i)} - y^{(i)}|,
  \end{equation*}
  and
  \begin{equation*}
    t_*^{-1/\alpha}\wedge \frac{t_*}{|u^{(i)} - v^{(i)}|^{1+\alpha}} \ge \left(t_*^{-1/\alpha}\wedge \frac{t_*}{(3|x^{(i)} - y^{(i)}|)^{1+\alpha}}\right) \ge \frac{1}{3^{1+\alpha}}\left(t_*^{-1/\alpha}\wedge \frac{t_*}{|x^{(i)} - y^{(i)}|^{1+\alpha}}\right).
  \end{equation*}
  Combining the above two cases, we obtain
  \begin{equation}\label{eq:uivi>xiyi}
    t_*^{-1/\alpha}\wedge \frac{t_*}{|u^{(i)} - v^{(i)}|^{1+\alpha}} \ge c_1\left(t_*^{-1/\alpha}\wedge \frac{t_*}{|x^{(i)} - y^{(i)}|^{1+\alpha}}\right)
  \end{equation}
  for some $c_1 = c_1(d,\alpha,R,\Lambda) > 0$. Since for all $(u,v) \in E_x\times E_y$,
  \begin{equation}\label{eq:detlauvLower}
    \delta_D(u)\wedge \delta_D(v) \ge \frac{\rho_{Q_x}(u)\wedge\rho_{Q_y}(v)}{\sqrt{1+\Lambda^2}} \ge \frac{\delta_0 - \delta_1}{\sqrt{1+\Lambda^2}}  \ge \delta_1 = (\delta_1t_*^{1/\alpha})\cdot t_*^{-1/\alpha},
  \end{equation}
  the above inequality (\ref{eq:uivi>xiyi}) together with Lemma \ref{lem:pUle-int} (with $a_1 = \delta_1t_*^{1/\alpha}$) yields that, there exists $c_2 =c_2(d,\alpha,R,\Lambda,\gamma) > 0$ such that
  \begin{align*}
    I_1 =&~ \inf_{u\in E_x,v\in E_y} p_D(t_*,u,v)\\
    \ge&~ c_2 \inf_{u\in E_x,v\in E_y} p(t_*,u,v)\\
    \ge&~ c_2 C_{\CTstable}^{-1} \prod_{i=1}^d \left(t_*^{-1/\alpha}\wedge \frac{t_*}{|u^{(i)} - v^{(i)}|^{1+\alpha}}\right)\\
    \ge&~ c_1^dc_2 C_{\CTstable}^{-1} \prod_{i=1}^d \left(t_*^{-1/\alpha}\wedge \frac{t_*}{|x^{(i)} - y^{(i)}|^{1+\alpha}}\right)\\
    \ge&~ c_1^dc_2 C_{\CTstable}^{-2}p(t_*,x,y).
  \end{align*}

  We next estimate the lower bound of $I_2$. The estimate for $I_3$ is similar so it will be omitted. If $\delta_D(x) \ge \delta_0$, then $E_x = B(x,\delta_1)$ and, $\delta_D(x)\wedge\delta_D(u) \ge \delta_1 = (\delta_1t_*^{1/\alpha})t_*^{-1/\alpha}$ for $u \in E_x$ by (\ref{eq:detlauvLower}). Hence, by Lemma \ref{lem:pUle-int} with $a_1 = \delta_1t_*^{1/\alpha}$, we obtain
  \begin{align*}
    I_2 =&~ \int_{E_x} p_D(t_*,x,u) du =\int_{B(x,\delta_1)} p_D(t_*,x,u) du\\
    \ge&~ c_3 \int_{B(x,\delta_1)} p(t_*,x,u) du = c_3 \int_{B(0,\delta_1)} p(t_*,0,u) du\\
    \ge&~ c_3 \int_{B(0,\delta_1)} p(t_*,0,u) du \left(1\wedge \frac{\delta_D(x)^{\alpha/2}}{\sqrt{t_*}}\right),
  \end{align*}
  where $\int_{B(0,\delta_1)} p(t_*,0,u) du$ is a positive constant independent of $x$ and $D$.

  Recall that $\sigma^D_{U_x^0}:=\inf\{t>0: X^D_t \in U_x^0\}$ is the first entry time  of $X^D$ for the set $U_x^0$. If $\delta_D(x) < \delta_0$, then $E_x = U_x = \cup_{u \in U_x^0} B(u,\delta_1)$, and by Markov property, we have
  \begin{align*}
    I_2 = &~\int_{E_x}p_D(t_*,x,u) du = \mathbb{P}_x \left(X^D_{t_*}  \in U_x\right) \\
    \ge&~ \mathbb{P}_x\Big(X^D \text{ hits $U_x^0$ before time } t_*
     \hbox{ and stays in $U_x$ for at least $t_*$ units of time}\Big)\\
    =&~ \mathbb{P}_x\Big(\sigma^D_{U_x^0} < t_*, ~\tau^D_{U_x} \circ \theta_{\sigma^D_{U_x^0}} > t_*\Big)\\
    =&~ \mathbb{E}_y\bigg[\sigma^D_{U_x^0} < t_*,\mathbb{P}_{X^D_{\sigma^D_{U_x^0}}} \Big(\tau^D_{U_x} > t_*\Big)\bigg]\\
    \ge&~ \inf_{z \in U_x^0}\mathbb{P}_{z} \Big(\tau^D_{U_x} > t_*\Big)\cdot\mathbb{P}_x\Big(\sigma^D_{U_x^0} < t_*\Big).
  \end{align*}
  Similar to (\ref{eq:detlauvLower}), we have for $z \in U_x^0$,
  \begin{equation*}
    \delta_D(z) \ge \frac{\rho_{Q_x}(z)}{\sqrt{1+\Lambda^2}} \ge \frac{\delta_0}{\sqrt{1+\Lambda^2}} = 2\delta_1,
  \end{equation*}
  and then,
  \begin{equation*}
    B:=B(z,\delta_1)\subset U_x \subset D.
  \end{equation*}
  By Lemma \ref{lem:LLE} with $t=t_*$, $r = \delta_1$ and $a_1 = t_*^{1/\alpha}/\delta_1$, there exists $a_2 = a_2(d,\alpha,R,\Lambda) > 0$ and $c_4 = c_4(d,\alpha,R,\Lambda) > 0$ such that $B(z,a_2t_*^{1/\alpha}) \subset B$ and
  \begin{equation*}
    p_B(t_*,z,v) \ge c_4 t_*^{-d/\alpha},\qquad v \in B(z,a_2t_*^{1/\alpha}).
  \end{equation*}
  Hence,
  \begin{align*}
    \mathbb{P}_z\left(\tau^D_{U_x} > t_*\right) =&~ \mathbb{P}_z\left(\tau_{U_x} > t_*\right) = \int_{U_x} p_{U_x}(t_*,z,v) dv\\
    \ge&~ \int_B p_B(t_*,z,v) dv \ge \int_{B(z,a_2t_*^{1/\alpha})} p_B(t_*,z,v) dv\\
    \ge&~ c_4 t_*^{-d/\alpha} |B(z,a_2t_*^{1/\alpha})| = c_4 t_*^{-d/\alpha} \cdot a_2^dt_*^{d/\alpha} |B(0,1)| \\
    =&~ c_4 a_2^d |B(0,1)| >0.
  \end{align*}
  On the other hand, by the definition  (\ref{eq:def-tast}) of $t_*$, (\ref{eq:exitLower}) and (\ref{eq:exitUpper}), we obtain
  \begin{align*}
    \mathbb{P}_x\Big(\sigma^D_{U_x^0} < t_*\Big) \ge&~ \mathbb{P}_x\left(\tau^D_{D_{Q_x(\delta_0,r_0)}} < t_*, X^D_{\tau^D_{D_{Q_x(\delta_0,r_0)}}} \in U_x^0\right)\\
    =&~ \mathbb{P}_x\Big(\tau^D_{D_{Q_x(\delta_0,r_0)}} < t_*, X^D_{\tau^D_{D_{Q_x(\delta_0,r_0)}}} \in U_x^0, \tau_{D_{Q_x(\delta_0,r_0)}} < \tau_D\Big)\\
    =&~ \mathbb{P}_x\Big(\tau_{D_{Q_x(\delta_0,r_0)}} < t_*, X_{\tau_{D_{Q_x(\delta_0,r_0)}}} \in U_x^0\Big)\\
    =&~ \mathbb{P}_x\Big(X_{\tau_{D_{Q_x(\delta_0,r_0)}}} \in U_x^0\Big) - \mathbb{P}_x\Big(X_{\tau_{D_{Q_x(\delta_0,r_0)}}} \in U_x^0, \tau_{D_{Q_x(\delta_0,r_0)}} \ge t_*\Big)\\
    \ge&~ \mathbb{P}_x\Big(X_{\tau_{D_{Q_x(\delta_0,r_0)}}} \in U_x^0\Big) - \mathbb{P}_x\Big(\tau_{D_{Q_x(\delta_0,r_0)}} \ge t_*\Big)\\
    =&~ \mathbb{P}_x\Big(X_{\tau_{D_{Q_x(\delta_0,r_0)}}} \in D_{Q_x(r_1,r_0)}\Big) - \mathbb{P}_x\Big(\tau_{D_{Q_x(\delta_0,r_0)}} \ge t_*\Big)\\
    \ge&~ \mathbb{P}_x\Big(X_{\tau_{D_{Q_x(\delta_0,r_0)}}} \in D_{Q_x(r_1,r_0)}\Big) - \frac{\mathbb{E}_x[\tau_{D_{Q_x(\delta_0,r_0)}}]}{t_*}\Big)\\
    \ge&~ C_{\CTexitL} \delta_D(x)^{\alpha/2} - \frac{C_{\CTexitU} \delta_D(x)^{\alpha/2}}{t_*}\\
    =&~ \left(C_{\CTexitL}\sqrt{t_*} - \frac{C_{\CTexitU}}{\sqrt{t_*}}\right) \frac{\delta_D(x)^{\alpha/2}}{\sqrt{t_*}}\\
    \ge&~ \sqrt{\frac{C_{\CTexitL} C_{\CTexitU}}{2}} \frac{\delta_D(x)^{\alpha/2}}{\sqrt{t_*}}\\
    \ge&~ \sqrt{\frac{C_{\CTexitL} C_{\CTexitU}}{2}} \left(1\wedge\frac{\delta_D(x)^{\alpha/2}}{\sqrt{t_*}}\right).
  \end{align*}
  Combining the above four inequalities, we obtain
  \begin{equation*}
    I_2 \ge c_5 \left(1\wedge \frac{\delta_D(x)^{\alpha/2}}{\sqrt{t_*}}\right).
  \end{equation*}
  Finally, combining the estimates of $I_1,I_2,I_3$ and (\ref{eq:pDle-1}), we finish the proof.
\end{proof}

Note that the constant $t_*$ from Lemma \ref{lem:pDTle} is greater than or equal to $1$.

\begin{lemma}\label{lem:pDtle}
  Let $D\subset \mathbb{R}^d$ be a $C^{1,1}$ open set with characteristics $(R,\Lambda)$ satisfying condition \hyperref[eq:conH]{$\bf (H_{\gamma})$} for some $\gamma \in (0,1]$. There is a constant $c= c(d,\alpha,R,\Lambda,\gamma)>0$ such that for all $t \in (0,3t_*]$, $x,y \in D$,
  \begin{equation*}
    p_D(t,x,y) \ge c\left(1\wedge \frac{\delta_D(x)^{\alpha/2}}{\sqrt{t}}\right)\left(1\wedge \frac{\delta_D(y)^{\alpha/2}}{\sqrt{t}}\right)p(t,x,y).
  \end{equation*}
\end{lemma}
\begin{proof}
  For $t\in (0,3t_*]$, let
  \begin{equation*}
    \lambda := \left(\frac{3t_*}{t}\right)^{1/\alpha} \ge 1.
  \end{equation*}
  Note that $\lambda D=\{\lambda z: z \in D\}$ is also a $C^{1,1}$ open set with the same characteristics of $D$ and satisfies condition \hyperref[eq:conH]{$\bf (H_{\gamma})$}. Then, by the scaling property (\ref{eq:DHKscale}), (\ref{eq:HKscale}), (\ref{eq:pEst}) and Lemma \ref{lem:pDTle} for $\lambda D$, we obtain for all $x,y \in D$
  \begin{align*}
    p_D(t,x,y) =&~ \lambda^d  p_{\lambda D}(\lambda^\alpha t,\lambda x,\lambda y)\\
    =&~ \lambda^d  p_{\lambda D}(3t_*,\lambda x,\lambda y)\\
    \ge&~ c\lambda^d \left(1\wedge \frac{\delta_{\lambda D}(\lambda x)^{\alpha/2}}{\sqrt{t_*}}\right) \left(1\wedge \frac{\delta_{\lambda D}(\lambda y)^{\alpha/2}}{\sqrt{t_*}}\right) p(t_*,\lambda x,\lambda y)\\
    =&~ c \left(1\wedge \frac{\lambda^{\alpha/2}\delta_D(x)^{\alpha/2}}{\sqrt{t_*}}\right) \left(1\wedge \frac{\lambda^{\alpha/2} \delta_D(y)^{\alpha/2}}{\sqrt{t_*}}\right)\lambda^d p(3^{-1}\lambda^\alpha t,\lambda x,\lambda y)\\
    =&~ c \left(1\wedge \frac{\sqrt{3}\delta_D(x)^{\alpha/2}}{\sqrt{t}}\right) \left(1\wedge \frac{\sqrt{3} \delta_D(y)^{\alpha/2}}{\sqrt{t}}\right) p(3^{-1}t,x,y)\\
    \ge &~ \frac{c}{3^d C_{\CTstable}^2} \left(1\wedge \frac{\delta_D(x)^{\alpha/2}}{\sqrt{t}}\right) \left(1\wedge \frac{\delta_D(y)^{\alpha/2}}{\sqrt{t}}\right) p(t,x,y).
  \end{align*}
\end{proof}

We need the following lemma to estimate Dirichlet heat kernel $p_D(t,x,y)$ for large time.

\begin{lemma}\label{lem:pUrawUp}
  Suppose that $U\subset \mathbb{R}^d$ is a bounded open set and let $\kappa \ge \mathrm{diam}(U)$. Then, there are two positive constants $c_i= c_i(d,\alpha,\kappa), i=1,2$ such that
  \begin{align*}
    p_U(t,x,y) \le c_1 e^{-c_2t},\quad t > 1, ~x,y \in U.
  \end{align*}
\end{lemma}

\begin{proof}
  Since $\{|z| \ge \kappa\} \supset \{z = (z^{(1)},\cdots,z^{(d)}): |z^{(1)}| \ge \kappa\vee 1\}$, by (\ref{eq:ppalpha}), for every $x\in U$, we have
  \begin{align*}
    \mathbb{P}_x(\tau_U \le 1) \ge&~ \mathbb{P}_x(X_1\in \mathbb{R}^d \setminus U) = \int_{\mathbb{R}^d\setminus U} p(1,x,z) dz
    \ge  \int_{|z| \ge \kappa} p(1,0,z)dz\\
    \ge&~ \int_{|z^{(1)}| \ge \kappa\vee 1} p(1,0,z)dz
    =  \int_{|z^{(1)}| \ge \kappa\vee 1} p^{(1,\alpha)} (1,0,z^{(1)})dz^{(1)} := c_3(\alpha,\kappa).
  \end{align*}
  Then,
  \begin{align*}
    \sup_{x\in U} \int_U p_U(1,x,y)dy &= \sup_{x\in U} \mathbb{P}_x(\tau_U > 1) \le (1-c_3) =: c_4 < 1.
  \end{align*}
  For $t \in (0,1]$ and $x\in U$, we have $\int_U p_U(t,x,z) dz \le \int_U p(t,x,z) dz \le 1\le ee^{-t}$. When $t\in (n,n+1]$ for some integer $n \ge 1$, we set $x_0 = x$ and then obtain by semigroup property,
  \begin{align*}
    \int_U p_U(t,x,z) dz &= \int_{U^{n+1}} \prod_{k=1}^n p_U(1,x_{k-1},x_k) p_U(t-n,x_n,z) dx_1\cdots dx_ndz\\
    &= \int_{U^{n+1}} p_U(t-n,x_n,z)dz ~p_U(1,x_{n-1},x_n)dx_n \cdots p_U(1,x_{0},x_1)dx_1\\
    &\le c_4^n \le c_4^{-1}c_4^t.
  \end{align*}
  Combining the above two inequalities, we have for all $(t,x)\in (0,\infty)\times U$,
  \begin{equation*}
    \int_U p_U(t,x,y) dy \le (e\vee c_4^{-1})e^{-c_2 t},
  \end{equation*}
  where $c_2 := \ln(c_4^{-1}\wedge e)$. Note that by (\ref{eq:pEst}), $p_U(1,z,y) \le p(1,z,y)\le C_{\CTstable}$ for all $z,y\in U$. Thus, for all $t > 1$ and $z,y\in U$ we have
  \begin{align*}
    p_U(t,x,y) =&~ \int_U p_U(t-1,x,z) p_U(1,z,y)dz\\
    \le&~ C_{\CTstable}\int_U p_U(t-1,x,z)dz \le C_{\CTstable}(e\vee c_4^{-1})e^{-c_2 (t-1)}.
  \end{align*}
  Setting $c_1 = C_{\CTstable}(e\vee c_4^{-1})e^{c_2}$, we finish the proof.
\end{proof}

\medskip

\subsection{Proof of Theorem \ref{T:1.1}}\label{Subs:4.3}

\begin{proof}[Proof of Theorem \ref{T:1.1}]
  (i) By Lemmas \ref{lem:pDtue} and \ref{lem:pDtle}, we obtain (\ref{eq:pDEst-small-1}) and (\ref{eq:pDEst-small-2}) for all $t \in (0,1]$. By the semigroup property and (\ref{eq:pDEst-small-1}) for $t \in (1/2,1]$, we have for $x,y \in D$,
  \begin{align*}
    p_D(2t,x,y) =&~ \int_D p_D(t,x,z) p_D(t,z,y) dz\\
    \le&~ c_1 \int_D \left(1\wedge\frac{\delta_D(x)^{\alpha/2}}{\sqrt{t}}\right) p(t,x,z) \left(1\wedge\frac{\delta_D(z)^{\alpha/2}}{\sqrt{t}}\right) \left(1\wedge\frac{\delta_D(z)^{\alpha/2}}{\sqrt{t}}\right) p(t,z,y) \left(1\wedge\frac{\delta_D(y)^{\alpha/2}}{\sqrt{t}}\right) dz\\
    \le&~ c_1 \left(1\wedge\frac{\delta_D(x)^{\alpha/2}}{\sqrt{t}}\right) \left(1\wedge\frac{\delta_D(y)^{\alpha/2}}{\sqrt{t}}\right) \int_{\mathbb{R}^d} p(t,x,z)p(t,z,y) dz\\
    =&~ c_1 \left(1\wedge\frac{\delta_D(x)^{\alpha/2}}{\sqrt{t}}\right) \left(1\wedge \frac{\delta_D(y)^{\alpha/2}}{\sqrt{t}}\right) p(2t,x,y)\\
    \le&~ 2c_1 \left(1\wedge\frac{\delta_D(x)^{\alpha/2}}{\sqrt{2t}}\right) \left(1\wedge \frac{\delta_D(y)^{\alpha/2}}{\sqrt{2t}}\right) p(2t,x,y).
  \end{align*}
  This shows that (\ref{eq:pDEst-small-1}) holds for $t \in (1,2]$. For general $T> 0$, one can repeat the above arguments for $\lceil \log_2T\rceil$ times to prove (\ref{eq:pDEst-small-1}) for $t \in (0,T]$. Here for $a\in \R$, the notation $\lceil a\rceil$ stands for the smallest integer greater than or equal to $a$.

  (ii) We next estimate the lower bound of $p_D(2t,x,y)$. Note that $D$ is a $C^{1,1}$ open set with characteristics $(R,\Lambda)$. Fix $x\in D$ and $t \in (1/2,1]$, and let $Q\in \partial D$ be such that $\delta_D(x) = |x-Q|$. Define
  \begin{equation*}
    x_0 =
    \begin{cases}
      Q + \frac{R}{4|x-Q|}(x-Q), &\quad \text{ if } \delta_D(x) < \frac{R}{4},\\
      x, &\quad \text{ if } \delta_D(x) \ge \frac{R}{4},
    \end{cases}
  \end{equation*}
  and, let $r := R/8$. We have
  \begin{equation*}
    B := B(x_0,r) \subset B(x_0,R/4) \subset D,
  \end{equation*}
  and $\delta_D(z) \ge r$ for all $z \in B$. By semigroup property of $p_D(2t,x,y)$ and (\ref{eq:pDEst-small-2}) for $t \in (1/2,1]$, we have
  \begin{equation}\label{eq:pD2le-1}
  \begin{split}
    p_D(2t,x,y) =&~ \int_D p_D(t,x,z) p_D(t,z,y) dz\\
    \ge&~ c_2 \int_D \left(1\wedge \frac{\delta_D(x)^{\alpha/2}}{\sqrt{t}}\right) p(t,x,z) \left(1\wedge \frac{\delta_D(z)^{\alpha/2}}{\sqrt{t}}\right) \left(1\wedge \frac{\delta_D(z)^{\alpha/2}}{\sqrt{t}}\right)p(t,z,y)\left(1\wedge \frac{\delta_D(y)^{\alpha/2}}{\sqrt{t}}\right) dz\\
    \ge&~ c_2 \left(1\wedge \frac{\delta_D(x)^{\alpha/2}}{\sqrt{t}}\right) \left(1\wedge \frac{\delta_D(y)^{\alpha/2}}{\sqrt{t}}\right) \int_{B} \left(1\wedge \delta_D(z)^{\alpha/2}\right)^2 p(t,x,z)p(t,z,y) dz\\
    \ge&~ c_2 r^\alpha \left(1\wedge \frac{\delta_D(x)^{\alpha/2}}{\sqrt{t}}\right) \left(1\wedge \frac{\delta_D(y)^{\alpha/2}}{\sqrt{t}}\right) \int_B p(t,x,z)p(t,z,y) dz.
  \end{split}
  \end{equation}
  Note that $x \in B(x_0,R/4)$ and for all $z \in B=B(x_0,r)$,
  \begin{equation*}
    |x-z| < R/4+R/8 = 3R/8.
  \end{equation*}
  We have, by (\ref{eq:pEst}), for $t \in (1/2,1]$,
  \begin{equation}\label{eq:pD2le-2}
    p(t,x,z) \ge C_{\CTstable}^{-1} \prod_{i=1}^d \left(t^{-\frac{1}{\alpha}}\wedge \frac{t}{|x^{(i)}-z^{(i)}|^{1+\alpha}}\right) \ge C_{\CTstable}^{-1} \prod_{i=1}^d \left(1\wedge \frac{1/2}{(3R/8)^{1+\alpha}}\right) = C_{\CTstable}^{-1} > 0.
  \end{equation}
  On the other hand, for $1\le i\le d$, if $|x^{(i)} - y^{(i)}| < R/2$, then
  \begin{equation*}
    |z^{(i)}-y^{(i)}| \le |z^{(i)}-x^{(i)}| + |x^{(i)}-y^{(i)}| \le 3R/8+R/2 < R,
  \end{equation*}
  and
  \begin{align*}
    t^{-\frac{1}{\alpha}}\wedge \frac{t}{|z^{(i)}-y^{(i)}|^{1+\alpha}} \ge&~ (2t)^{-\frac{1}{\alpha}}\wedge \frac{2t}{2R^{1+\alpha}}\\
    =&~ \left(1\wedge \frac{(2t)^{1+1/\alpha}}{2R^{1+\alpha}}\right) (2t)^{-\frac{1}{\alpha}}\\
    \ge&~ \frac{1}{2} \left((2t)^{-\frac{1}{\alpha}}\wedge \frac{2t}{|x^{(i)}-y^{(i)}|^{1+\alpha}}\right).
  \end{align*}
  If $|x^{(i)} - y^{(i)}| \ge R/2$, then
  \begin{equation*}
    |z^{(i)}-y^{(i)}| \le |z^{(i)}-x^{(i)}| + |x^{(i)}-y^{(i)}| \le 3R/8 + |x^{(i)}-y^{(i)}| \le 2|x^{(i)}-y^{(i)}|,
  \end{equation*}
  and
  \begin{align*}
    t^{-\frac{1}{\alpha}}\wedge \frac{t}{|z^{(i)}-y^{(i)}|^{1+\alpha}} \ge&~ (2t)^{-\frac{1}{\alpha}}\wedge \frac{2t}{2|z^{(i)}-y^{(i)}|^{1+\alpha}}\\
    \ge&~ (2t)^{-\frac{1}{\alpha}}\wedge \frac{2t}{2(2|x^{(i)}-y^{(i)}|)^{1+\alpha}}\\
    \ge&~ \frac{1}{2^{2+\alpha}}\left((2t)^{-\frac{1}{\alpha}}\wedge \frac{2t}{|x^{(i)}-y^{(i)}|^{1+\alpha}}\right).
  \end{align*}
  Combining the above two cases, we always have
  \begin{equation*}
    t^{-\frac{1}{\alpha}}\wedge \frac{t}{|z^{(i)}-y^{(i)}|^{1+\alpha}} \ge c_3 \left((2t)^{-\frac{1}{\alpha}}\wedge \frac{2t}{|x^{(i)}-y^{(i)}|^{1+\alpha}}\right),
  \end{equation*}
  for some constant $c_3 >0$ independent of $t,x,y$. Hence, by (\ref{eq:pEst}), we have
  \begin{align*}
    p(t,z,y) \ge&~ C_{\CTstable}^{-1} \prod_{i=1}^d \left(t^{-\frac{1}{\alpha}}\wedge \frac{t}{|z^{(i)}-y^{(i)}|^{1+\alpha}}\right)\\
    \ge&~ c_3^dC_{\CTstable}^{-1} \prod_{i=1}^d \left((2t)^{-\frac{1}{\alpha}}\wedge \frac{2t}{|x^{(i)}-y^{(i)}|^{1+\alpha}}\right)\\
    \ge&~ c_3^dC_{\CTstable}^{-2} p(2t,x,y).
  \end{align*}
  Combining this, (\ref{eq:pD2le-1}) and (\ref{eq:pD2le-2}), we have
  \begin{align*}
    p_D(2t,x,y) \ge&~ c_2 r^\alpha \left(1\wedge \frac{\delta_D(x)^{\alpha/2}}{\sqrt{t}}\right) \left(1\wedge \frac{\delta_D(y)^{\alpha/2}}{\sqrt{t}}\right) \int_B C_{\CTstable}^{-1} \cdot c_3^dC_{\CTstable}^{-2} p(2t,x,y) dz\\
    \ge&~ c_2 c_3^d  r^\alpha C_{\CTstable}^{-3} \left(1\wedge \frac{\delta_D(x)^{\alpha/2}}{\sqrt{t}}\right) \left(1\wedge \frac{\delta_D(y)^{\alpha/2}}{\sqrt{t}}\right) |B|~ p(2t,x,y) \\
    \ge&~ c_2 c_3^d  r^{d+\alpha} C_{\CTstable}^{-3} \left(1\wedge \frac{\delta_D(x)^{\alpha/2}}{\sqrt{2t}}\right) \left(1\wedge \frac{\delta_D(y)^{\alpha/2}}{\sqrt{2t}}\right) |B(0,1)| ~p(2t,x,y).
  \end{align*}
  Therefore, we have proved (\ref{eq:pDEst-small-2}) for $t \in (1,2]$. Similarly, one can repeat the above arguments for $\lceil\log_2 T\rceil$ times to prove (\ref{eq:pDEst-small-2}) for $t \in (0,T]$.

  (iii) Now, assume in addition that $D$ is a \emph{bounded} $C^{1,1}$ open set and satisfies \hyperref[eq:conH]{$\bf (H_{\gamma})$} for some $\gamma \in (0,1]$. Recall that $\mathcal{L}^D$ is the infinitesimal generator of the semigroup $\{P^D_t,t\ge0\}$ on $L^2(D,dx)$. Since for each $t > 0$, the heat kernel $p_D(t,x,y)$ is bounded on $D \times D$, it follows from Jentzsch's Theorem (\cite[Theorem V.6.6, p. 337]{Schaefer.1974.376}) that the value $-\lambda_1(D) = \sup (\sigma(\mathcal{L}^D))$ is an eigenvalue of multiplicity 1 for $\mathcal{L}^D$ and that the eigenfunction $\phi_D$ associated with $\lambda_1(D)$ can be chosen to be strictly positive with $\|\phi_D\|_{L^2(D)}=1$. In the rest of the proof, we write $\lambda_1(D)$ as $\lambda_1$ for simplicity.

\medskip

  \textbf{Step 1.} We first prove the second inequality in (\ref{eq:eigenBnd}). Since $\phi_D$ is the eigenfunction of $\mathcal{L}^D$ associated with $\lambda_1$, we have for all $t > 0$ and $x\in D$,
  \begin{equation}\label{eq:phit}
    \phi_D(x) = e^{\lambda_1 t} P^D_t \phi_D(x) = e^{\lambda_1 t}\int_D p_D(t,x,y) \phi_D(y) dy.
  \end{equation}
  Setting $t = 1/4$ in (\ref{eq:phit}), by (\ref{eq:pDEst-small-1}) with $T = 1$,  \eqref{eq:pEst} and H\"{o}lder inequality, we have for all $x\in D$,
  \begin{align}\label{eq:phiUpdelta}
    \phi_D(x) \le&~ c_4 e^{\frac{1}{4}\lambda_1} (1\wedge2\delta_D(x)^{\alpha/2}) \int_D p(1/4,x,y) \phi_D(y) dy\notag\\
    \le&~ 2c_4 e^{\frac{1}{4}\lambda_1}(1\wedge \delta_D(x)^{\alpha/2}) \sqrt{\int_D p(1/4,x,y)^2 dy} \cdot \|\phi_D\|_{L^2(D)} \notag\\
    \le&~ 2c_4 e^{\frac{1}{4}\lambda_1}(1\wedge \delta_D(x)^{\alpha/2})\sqrt{p(1/2,x,x)}  \notag\\
    \le&~ 2c_4 \sqrt{2^{d/\alpha}C_{\CTstable}} e^{\frac{1}{4}\lambda_1}(1\wedge \delta_D(x)^{\alpha/2}) \notag\\
    =:\!\!&~ \,c_5 e^{\frac{1}{4}\lambda_1}(1\wedge \delta_D(x)^{\alpha/2}),
  \end{align}
  where $c_i = c_i(d,\alpha,R,\Lambda,\gamma)>0$, $i = 4,5$.

\smallskip

  On the other hand, set $\kappa := {\rm diam}(D) $ to be the diameter of $D$ for simplicity. Setting $t = 1$ in (\ref{eq:phit}), by (\ref{eq:pDEst-small-2}) with $T = 1$ and \eqref{eq:pEst}, we have for all $x \in D$,
  \begin{align*}
    \phi_D(x) \ge&~ c_4^{-1} (1\wedge \delta_D(x)^{\alpha/2}) e^{\lambda_1}\int_D (1\wedge \delta_D(y)^{\alpha/2}) p(1,x,y) \phi_D(y) dy\\
    \ge&~ c_4^{-1} C_{\CTstable}^{-1} \left(1\wedge \kappa^{-d(1+\alpha)}\right) (1\wedge \delta_D(x)^{\alpha/2}) e^{\lambda_1}\int_D (1\wedge \delta_D(y)^{\alpha/2}) \phi_D(y) dy\\
    =:\!\! &~ \, c_6 (1\wedge \delta_D(x)^{\alpha/2}) e^{\lambda_1}\int_D (1\wedge \delta_D(y)^{\alpha/2}) \phi_D(y) dy,
  \end{align*}
  where $c_6 = c_6 (d,\alpha,R,\Lambda,\gamma,\kappa)>0$. Combining this and (\ref{eq:phiUpdelta}), we have for $x\in D$,
  \begin{align}\label{eq:phiLowdelta}
    \phi_D(x) \ge&~ c_6 (1\wedge \delta_D(x)^{\alpha/2}) e^{\lambda_1}\int_D c_5^{-1} e^{-\frac{1}{4}\lambda_1} (\phi_D(y))^2 dy\notag\\
    \ge&~ c_6c_5^{-1} e^{\frac{3}{4}\lambda_1}(1\wedge \delta_D(x)^{\alpha/2}) \int_D (\phi_D(y))^2  dy\notag\\
    \ge&~ c_6c_5^{-1} e^{\frac{3}{4}\lambda_1}(1\wedge \delta_D(x)^{\alpha/2}).
  \end{align}
  Combining (\ref{eq:phiUpdelta}) and (\ref{eq:phiLowdelta}), we have for $x\in D$,
  \begin{equation*}
    c_6c_5^{-1} e^{\frac{3}{4}\lambda_1}(1\wedge \delta_D(x)^{\alpha/2}) \le \phi_D(x) \le  c_5 e^{\frac{1}{4}\lambda_1}(1\wedge \delta_D(x)^{\alpha/2}),
  \end{equation*}
  which implies that
  \begin{equation}\label{eq:lambda0Bnd}
    \lambda_1 \le 2\ln(c_5^2/c_6) < \infty.
  \end{equation}
  This is exactly the second inequality in (\ref{eq:eigenBnd}).

\medskip

  \textbf{Step 2.} We next show  (\ref{eq:pDEst-large}). In view of (\ref{eq:phiUpdelta}), (\ref{eq:phiLowdelta}) and (\ref{eq:lambda0Bnd}), we have for $x \in D$,
  \begin{equation}\label{eq:phiEqvdelta}
    c_7^{-1} (1\wedge\delta_D(x)^{\alpha/2}) \le \phi_D(x) \le c_7 (1\wedge\delta_D(x)^{\alpha/2}),
  \end{equation}
  where $c_7= c_7(d,\alpha,R,\Lambda,\gamma,\kappa) \ge 1$.

  Now, multiplying (\ref{eq:phit}) by $\phi_D(x)$ and integrating it over $D$ with respect to $dx$, we have for $t > 0$,
  \begin{equation*}
    1 = \int_D (\phi_D(x))^2 dx = e^{\lambda_1t}\int_{D\times D} \phi_D(x) p_D(t,x,y) \phi_D(y) dx dy,
  \end{equation*}
  which implies that
  \begin{equation*}
    \int_{D\times D} \phi_D(x) p_D(t,x,y) \phi_D(y) dx dy = e^{-\lambda_1t}.
  \end{equation*}
  Combining this and (\ref{eq:phiEqvdelta}), we have for $t > 0$,
  \begin{equation}\label{eq:deltaLambdaEqv}
    c_7^{-2} e^{-\lambda_1t} \le \int_{D\times D}(1\wedge\delta_D(x)^{\alpha/2}) p_D(t,x,y) (1\wedge\delta_D(y)^{\alpha/2}) dx dy \le c_7^2 e^{-\lambda_1t} .
  \end{equation}
  For any $T > 0$, set $t_0 = (T\wedge 1)/4$. Note that, by (\ref{eq:pDEst-small-1}) and (\ref{eq:pDEst-small-2})(with $T=1$), there is a constant $c_8 = c_8(d,\alpha,R,\Lambda,\gamma) \ge 1$ such that for $(u,v) \in D\times D$,
  \begin{align}
    p_D(t_0,u,v) \le&~ c_8 t_0^{-d/\alpha-1} (1\wedge\delta_D(u)^{\alpha/2})(1\wedge\delta_D(v)^{\alpha/2}),\label{eq:smallUp}\\
    p_D(t_0,u,v) \ge&~ c_8^{-1}\left(t_0^{-1/\alpha}\wedge (t_0\kappa^{-(1+\alpha)})\right)^d (1\wedge\delta_D(u)^{\alpha/2}) (1\wedge\delta_D(v)^{\alpha/2}).\label{eq:smallLow}
  \end{align}
  Combining (\ref{eq:smallUp}), the semigroup property of $p_D$ and (\ref{eq:deltaLambdaEqv}), we have for all $(t,x,y) \in (T,\infty)\times D\times D$,
  \begin{align}\label{eq:LargeUp}
    p_D(t,x,y) =&~ \int_{D\times D} p_D(t_0,x,u) p_D(t-2t_0,u,v) p_D(t_0,v,y) dudv\notag\\
    \le&~ c_8^2 t_0^{-2d/\alpha-2} \delta_D(x)^{\alpha/2}\delta_D(y)^{\alpha/2} \int_{D\times D}(1\wedge\delta_D(u)^{\alpha/2}) p_D(t-2t_0,u,v) (1\wedge\delta_D(v)^{\alpha/2}) du dv\notag\\
    \le&~ c_7^2c_8^2 t_0^{-2d/\alpha-2} e^{-(t-2t_0)\lambda_1} \delta_D(x)^{\alpha/2}\delta_D(y)^{\alpha/2}.
  \end{align}
  Similarly, by (\ref{eq:smallLow}), the semigroup property of $p_D$ and (\ref{eq:deltaLambdaEqv}), we have for all $(t,x,y) \in (T,\infty)\times D\times D$,
  \begin{align}\label{eq:LargeLow}
    p_D(t,x,y) =&~ \int_{D\times D} p_D(t_0,x,u)p_D(t-2t_0,u,v)p_D(t_0,v,y) dudv\notag\\
    \ge&~ c_8^{-2}\left(t_0^{-1/\alpha}\wedge (t_0\kappa^{-(1+\alpha)})\right)^{2d} (\kappa\vee 1)^{-\alpha} \delta_D(x)^{\alpha/2} \delta_D(y)^{\alpha/2}\notag\\
    &\qquad \qquad \cdot \int_{D\times D}(1\wedge\delta_D(u)^{\alpha/2}) p_D(t-2t_0,u,v) (1\wedge\delta_D(v)^{\alpha/2}) du dv\notag\\
    \ge&~ c_8^{-2}\left(t_0^{-1/\alpha}\wedge (t_0\kappa^{-(1+\alpha)})\right)^{2d} (\kappa\vee 1)^{-\alpha} e^{-(t-2t_0)\lambda_1}\delta_D(x)^{\alpha/2} \delta_D(y)^{\alpha/2}.
  \end{align}
  Combining (\ref{eq:lambda0Bnd}), (\ref{eq:LargeUp}) and (\ref{eq:LargeLow}), we obtain (\ref{eq:pDEst-large}).

\medskip

  \textbf{Step 3.} For the first inequality in (\ref{eq:eigenBnd}), by Lemma \ref{lem:pUrawUp}, there exist $c_i = c_i(d,\alpha,\kappa) > 0$, $i = 9,10$ such that
  \begin{equation*}
    p_D(t,x,y) \le c_9 e^{-c_{10} t},\qquad t> 0, ~x,y \in D.
  \end{equation*}
  This together with the first inequality in (\ref{eq:deltaLambdaEqv}) yields that for all $t > 1$,
  \begin{align*}
    c_7^{-2} e^{-\lambda_1t} \le&~ \int_{D\times D}(1\wedge\delta_D(x)^{\alpha/2}) p_D(t,x,y) (1\wedge\delta_D(y)^{\alpha/2}) dx dy\\
    \le&~ c_9 e^{-c_{10} t} \int_{D\times D}(1\wedge\delta_D(x)^{\alpha/2}) (1\wedge\delta_D(y)^{\alpha/2}) dx dy\\
    \le&~ c_9 e^{-c_{10} t} |D|^2.
  \end{align*}
  Rewriting the above inequality, we have for all $t > 1$,
  \begin{equation*}
    e^{(c_{10} - \lambda_1)t} \le c_7^2 c_9 |D|^2 < \infty.
  \end{equation*}
  Since the above inequality holds for all $t > 1$, we obtain that
  \begin{equation*}
    \lambda_1 \ge c_{10},
  \end{equation*}
  which is exactly the first inequality in (\ref{eq:eigenBnd}).
\end{proof}

\section{Irreducibility}\label{S:irred}

Let $D\subset \mathbb{R}^d$ be a non-empty open set. In this section, we present a proof for Theorem \ref{thm:irreducible} on the irreducibility of $X^D$ and a proof for Corollary \ref{C:1.4}. For $z\in D$, define
\begin{align} \label{e:5.1}
  U_z :=&~ \{w \in D: \text{there exist finitely many points $\{x_i\}_{i=0}^N\subset D$ with $x_0=z$ and $x_N = w$ so that} \nonumber \\
  &~\quad \qquad \quad \text{ each   pair $(x_{i-1}, x_i)$, $1\leq i\leq N$,  has only one different coordinate}\}.
\end{align}
Note that $U_z\supset B(z, r)$ for any $r>0$ such that $B(z, r)\subset D$.

\begin{lemma}\label{L:5.1}
  For each $z\in D$, $U_z$ is both open and closed in $D$.
\end{lemma}

\proof
  For any $w\in  U_z$, there is $r > 0$ so that $B(w,r)\subset D$ and so $B(w,r) \subset U_z$ by the definition of $U_z$. This shows that $U_z$ is an open subset of $D$. If $U_z=D$, then $U_z$ is clearly closed in $D$. Suppose now $U_z\not= D$. Then for any $w\in D\setminus U_z$, there is some $r>0$ so that $B(w, r)\subset D$. Note that $B(w,r) \cap U_z = \emptyset$, as otherwise there would be some $w_1\in B(w,r)\cap U_z$ which would imply that $w\in U_z$. Hence  $U_z$ is a closed subset of $D$.
\qed

\begin{thm}\label{T:5.2}
  Suppose that $x_0 $ and $y_0$ are two distinct points in $D$ so that $y_0\in U_{x_0}$. Then $p_D(t, x_0, y_0)>0$ for every $t>0$.
\end{thm}

\proof
  Fix an arbitrary $t>0$. Set $c_1= t^{-1/\alpha}(\delta_D(x_0)\wedge \delta_D(y_0)) > 0$. According to Lemma \ref{lem:pUnle-int}, there is some constant $c_2>0$ so that $p_D(t,x_0,y_0) > 0$ whenever $|x_0-y_0|\le c_2 t^{1/\alpha}$.
  Suppose that $|x_0-y_0| >  c_2 t^{1/\alpha}$. We have by \eqref{eq:pUle-int-1},
  \begin{equation}\label{eq:irre-1}
    p_D(t,x_0,y_0) \ge \inf_{z\in B(y_0,c_3 t^{1/\alpha})} p_D( t/2, z, y_0) \mathbb{P}_{x_0} \left(X_{  t/2}^D \in B(y_0,c_3 t^{1/\alpha})\right),
  \end{equation}
  where $c_3 \in (0,(c_1\wedge c_2)/2)$  is to be chosen sufficiently small later. Note that for $z \in B(y_0,c_3 t^{1/\alpha})$,
  \begin{equation*}
    \delta_D(z) \wedge \delta_D(y_0) \ge (c_1-c_3)t^{1/\alpha} \ge \frac{c_1}{2}t^{1/\alpha}.
  \end{equation*}
  By Lemma \ref{lem:pUnle-int}, there exists $c_4 >0$ and $c_5>0$ depending on $d,\alpha$ and $c_1$ such that
  \begin{equation*}
    p_D(t/2, z,y_0) \ge c_5 t^{-d/\alpha}
  \end{equation*}
  for every $z$ with $|z-y_0| < c_4t ^{1/\alpha}$.  Taking $c_3\in (0,(c_1\wedge c_2)/2)$ small enough so that $c_3<c_4$, we get by \eqref{eq:irre-1}
  that
  \begin{equation}\label{eq:irre-3}
    p_D(t,x_0,y_0) \ge c_5  t^{-d/\alpha} \mathbb{P}_{x_0} \left(X_{t/2}^D \in B(y_0,c_3 t^{1/\alpha})\right).
  \end{equation}

  We next show that $\mathbb{P}_{x_0}\left(X_{t/2}^D \in B(y_0, c_3t^{1/\alpha})\right) > 0$. Let $\{x_i\}_{i=1}^N$ be a finite sequence of points in $D$ in the definition for $y_0\in U_{x_0}$. Define
  \begin{equation*}
    r := \frac{1}{4\sqrt{d}} \min\left\{\min_{1\le i\le N} \{|x_{i-1} - x_i|\}, \min_{0\le i\le N} \{\delta_D(x_i)\}, c_3 t^{1/\alpha}\right\},
  \end{equation*}
  and
  \begin{equation*}
    Q_i := Q(x_i,r)  \quad  \hbox{for } 0\le i\le N.
  \end{equation*}
  For each $0\le i\le N$, $Q_i \subset D$ and $\delta_D(z) \ge c_3 t^{1/\alpha}$ for all $z \in Q_i$. Note also that $Q_N \subset B(y_0,c_3 t^{1/\alpha})$.

  For $\lambda>0$ and a cube $Q:= Q(z,\bar{r})$, we use the notation $\lambda Q$ to denote the cube $Q(z,\lambda \bar{r})$, that is $\lambda Q(z,\bar{r}) = Q(z,\lambda \bar{r})$. We have as that for \eqref{eq:pUle-int-3},
  \begin{align}\label{eq:irre-4}
    &~\mathbb{P}_{x_0}\left(X_{t/2}^D \in B(y_0,c_3 t^{1/\alpha})\right) \notag\\
    = &~ \int_{B(y_0,c_3 t^{1/\alpha})} p_D(t/2,x_0,z_N) dz_N \notag\\
    = &~ \int_{B(y_0,c_3 t^{1/\alpha})} \left(\int_{D^{N-1}}  p_D({t/(2N)},x_0,z_1) p_D({t/(2N)},z_1,z_2) \cdots p_D({t/(2N)},z_{N-1},z_N)dz_1\cdots dz_{N-1}\right) dz_N  \notag\\
    \ge &~ \int_{Q_N} \int_{\left(2^{1-N}Q_1\right)\times\left(2^{2-N}Q_2\right) \times \cdots \times \left(2^{-1}Q_{N-1}\right)}  \notag\\
    &~\qquad p_D({t/(2N)},x_0,z_1) p_D({t/(2N)},z_1,z_2) \cdots p_D({t/(2N)},z_{N-1},z_N)dz_1\cdots dz_{N-1}dz_N  \notag\\
    \ge&~ \mathbb{P}_{x_0}\left(X_{t/(2N)}^D \in 2^{1-N}Q_1\right)\inf_{z_1\in2^{1-N}Q_1} \mathbb{P}_{z_1}\left(X_{t/(2N)}^D \in 2^{2-N}Q_2\right)  \notag\\
    &~ \times\cdots \times\inf_{z_{N-2}\in 2^{-2}Q_{N-2}} \mathbb{P}_{z_{N-2}}\left(X_{t/(2N)}^D \in 2^{-1}Q_{N-1}\right)\inf_{z_{N-1}\in2^{-1}Q_{N-1}} \mathbb{P}_{z_{N-1}}\left(X_{t/(2N)}^D \in Q_N\right) \notag\\
    \ge&~ \inf_{z_0\in 2^{-N} Q_0}\mathbb{P}_{z_0}\left(X_{t/(2N)}^D \in 2^{1-N}Q_1\right)\inf_{z_1\in2^{1-N}Q_1} \mathbb{P}_{z_1}\left(X_{t/(2N)}^D \in 2^{2-N}Q_2\right) \notag\\
    &~ \times\cdots \times\inf_{z_{N-2}\in2^{-2}Q_{N-2}} \mathbb{P}_{z_{N-2}}\left(X_{t/(2N)}^D \in 2^{-1}Q_{N-1}\right)\inf_{z_{N-1}\in2^{-1}Q_{N-1}} \mathbb{P}_{z_{N-1}}\left(X_{t/(2N)}^D \in Q_N\right).
  \end{align}
  We next estimate the lower bound of the right hand side of the above inequality. Let $1\leq k \leq N$. Note that the centers  $x_{k-1}$ and $x_k$ of $Q_{k-1}$ and $Q_k$ differ by only one coordinate. Thus there exists some $1\le i_k \le d$ and  $a_k\neq 0$ so that
  \begin{equation*}
    Q_k = Q_{k-1} + a_k e_{i_k}.
  \end{equation*}
  Hence, by Remark \ref{rem:para} (with $c_1 = r t^{-1/\alpha}$), we have
  \begin{equation*}
    \inf_{z\in2^{k-1-N} Q_{k-1}} \mathbb{P}_z\left(X_{t/(2N)}^D \in 2^{k-N}Q_k\right) \ge c_6\left(1\wedge \frac{t^{1+1/\alpha}}{|a_k|^{1+\alpha}}\right) > 0,
  \end{equation*}
  where the constant $c_6 > 0$ that may depend on $t$. This together with \eqref{eq:irre-3}-\eqref{eq:irre-4} yields that $p_D(t,x_0,y_0)>0$ when $|x_0-y_0| \ge c_2 t^{1/\alpha}$.  Combining the above two cases, we see that $p_D(t,x_0,y_0)>0$ for any $t>0$.
\qed

\medskip

We next establish a converse of  Theorem \ref{T:5.2}.

\begin{thm}\label{T:5.3}
  Suppose that $x_0\in D$ and $y_0\in D\setminus U_{x_0}$. Then $p_D(t,x,y)=0$ on
  $(0,\infty)\times U_{x_0} \times U_{y_0}$. In particular, $p_D(t,x_0,y_0)=0$ for all $t>0$.
\end{thm}

\proof
  Suppose $x_0\in D$ and $y_0\in  D\setminus U_{x_0}$. Then $U_{x_0}\cap U_{y_0} =\emptyset$ and
  \begin{equation}\label{eq:irre-6}
    z + ae_i \notin D \setminus U_{x_0}  \quad \hbox{for  every } a\in \R,  \  z\in U_{x_0} \hbox{ and } 1\leq i \leq d.
  \end{equation}
  Since $U_{x_0}$ is both open and closed by Lemma \ref{L:5.1},  we have  for every $x\in U_{x_0}$, with $\tau:= \tau_{U_{x_0}}$,
  $$
    X_{\tau-} \in  U_{x_0} \quad \hbox{and} \quad X_{\tau} \in D\setminus U_{x_0} \qquad  \P_x \hbox{-a.s. on }  \{\tau <\tau_D\}.
  $$
  By the strong Markov property of $X^D$,  \eqref{eq:XLS} and \eqref{eq:irre-6},  we have for every  $t > 0$ and $x \in U_{x_0}$ and $y \in U_{y_0}$,
  \begin{equation}\label{e:5.7}
    p_D(t,x, y) = \mathbb{E}_x\left[p_D(t-\tau,X^D_{\tau},y); \tau <t\right] =0.
  \end{equation}
\qed

\medskip

\begin{proof}[Proof of Theorem \ref{thm:irreducible}]
  (i) (Sufficient condition)  Suppose that the property \eqref{e:1.13} holds.  We have by Theorem \ref{T:5.2} that $p_D(t,x,y)>0$ for every $t>0$ and $x, y\in D$. This in particular implies that $X^D$ is irreducible as for any non-empty open subset $U\subset D$,  $\P_x(X_t^D \in U)= \int_U p_D(t,x,y) dy >0$ for every $t>0$ and $x\in D$.

\smallskip

  (ii) (Necessary condition) Suppose that $X^D$ is irreducible in $D$. Recall that for $z\in D$, $U_z$ is the  set defined by \eqref{e:5.1}. Were there two distinct $x, y\in D$ so that $U_x\cap U_y=\emptyset$, we would have by Theorem \ref{T:5.3}
  $$
    p_D(t, x, z)=0 \quad \hbox{for all } t>0 \hbox{ and }  z\in U_y.
  $$
  Consequently,
  $$
    \E_x \int_0^{\tau_D} \1_{U_y}(X_s) ds = \int_0^\infty \int_{U_y} p_D (t, x, z) dz dt =0.
  $$
  Since $X^D$ is irreducible, we have, by Lemma \ref{L:5.1} that $\P_x(\sigma_{U_y} <\tau_D)>0$. Since $U_y$ is a non-empty open subset of $D$, $\P_x$-a.s. on $\{\sigma_{U_y} <\tau_D\}$, $X^D$ spends positive Lebesgue amount of the time in $U_y$ in view of the right continuity of the sample paths of $X$. Thus $\E_x \int_0^{\tau_D} \1_{U_y}(X_s) ds >0$. This contradiction proves that $U_x\cap U_y \not= \emptyset$ for every $x\not=y $ in $D$. In other words, $U_x=D$ for every $x\in D$. Thus the property \eqref{e:1.13} holds.

\smallskip

  This completes the proof of the theorem.
\end{proof}

\smallskip

\begin{proof}[Proof of Corollary \ref{C:1.4}]
  (i) This follows directly from  Theorem \ref{thm:irreducible} and the connectedness of $D_1$ as $p_D(t,x,y)$ $\geq$ $p_{D_1}(t,x,y) >0$ for $(t,x,y) \in (0, \infty)\times D_1\times D_1$.

\smallskip

  (ii) If there is some $x_0\in D_1$ and $y_0\in D_2$ so that $y_0\in U_{x_0}$ (which is equivalent to $x_0\in U_{y_0}$), then, since $D_1$ and $D_2$ are connected, $U_z  \supset D_1\cup D_2$ for every $z\in D_1\cup D_2$. In this case, we have by Theorem \ref{T:5.2} that  $p(t,x,y) >0$ for any $(t,x,y) \times D_1 \times D_2$. Otherwise, $U_x\cap U_y =\emptyset$  for any $x\in D_1$ and $y\in D_2$. We then have by Theorem \ref{T:5.3} that $p_D(t,x,y)=0$ for every $(t,x,y)\in (0, \infty)\times D_1\times D_2$.
\end{proof}

\medskip

\section{Examples}\label{S:6}

In this section, we present three more examples and present a proof for Theorem \ref{T:1.6}. The first two show that the lower bound estimate (\ref{eq:pDEst-small-2}) in Theorem \ref{T:1.1}(ii) may fail for some smooth bounded connected open sets that do not satisfy condition \hyperref[eq:conH]{$\bf (H_{\gamma})$} for any $\gamma \in (0, 1]$. The third presents a bounded $C^{1,1}$ open set that does not satisfy the irreducibility condition \eqref{e:1.13} but for which we can derive two-sided sharp estimates for its Dirichlet heat kernel.

Recall that
\begin{equation*}
  j(a,b) = \frac{\sC_{1,\alpha}}{|a-b|^{1+\alpha}} \quad \hbox{for } a \not= b \in \R,
\end{equation*}
where $\sC_{1,\alpha}$ is the positive constant in \eqref{e:1.2}.

\bigskip

\begin{figure}[htp]
  \begin{tikzpicture}[
      scale = 0.8
    ]

    \draw [rounded corners = 10pt] (1,-1) -- (-5,-1) -- (-5,9) -- (5,9) -- (5,3) -- (3,3) -- (3,7) -- (-3,7) -- (-3,1) -- (1,1) -- cycle;

    \coordinate (a) at (1,-1);
    \coordinate (b) at (-5,-1);
    \coordinate (c) at (-5,9);
    \coordinate (d) at (5,9);
    \coordinate (e) at (5,3);
    \coordinate (f) at (3,3);
    \coordinate (g) at (3,7);
    \coordinate (h) at (-3,7);
    \coordinate (i) at (-3,1);
    \coordinate (j) at (1,1);

    \fill (0,0) circle (1.5pt);
    \draw [->] (2,0) -- (0.15,0);
    \node [right] at (2,0) {\Large $x = (0,0)$};

    \fill (4,4) circle (1.5pt);
    \draw [->] (6,4) -- (4.15,4);
    \node [right] at (6,4) {\Large $y = (4,4)$};

    \begin{scope}
      \clip [rounded corners = 10pt] (1,-1) -- (-5,-1) -- (-5,9) -- (5,9) -- (5,3) -- (3,3) -- (3,7) -- (-3,7) -- (-3,1) -- (1,1) -- cycle;
      \fill [pattern= north east lines] (-1,-1) -- (-5,-1) -- (-5,1) -- (-1,1) -- cycle;
      \fill [pattern= north east lines] (-1,7) -- (-1,9) -- (1,9) -- (1,7) -- cycle;
    \end{scope}

    \node [below] at (0,-1) {\Large $U_1$};
    \node [below] at ($(i)+2*(0,-1)$) {\Large $U_2$};
    \node at (-4,5) {\Large $U_3$};
    \node [below] at (0,7) {\Large $U_4$};
    \node at (4,6) {\Large $U_5$};

    \draw [gray,very thin] ($(b)+0.5*(-1,0)$) -- ($(b)+1.5*(-1,0)$);
    \draw [gray,very thin] ($(c)+0.5*(-1,0)$) -- ($(c)+1.5*(-1,0)$);
    \draw [<->,gray,very thin] ($(b)+1*(-1,0)$) -- ($(c)+1*(-1,0)$) node[pos =0.5,sloped,above] {10};

    \draw [gray,very thin] ($(d)+0.5*(0,1)$) -- ($(d)+1.5*(0,1)$);
    \draw [gray,very thin] ($(c)+0.5*(0,1)$) -- ($(c)+1.5*(0,1)$);
    \draw [<->,gray,very thin] ($(d)+1*(0,1)$) -- ($(c)+1*(0,1)$) node[pos =0.5,sloped,above] {10};

    \draw [gray,very thin]($(e)+0.3*(0,-0.8)$) -- ($(e)+1.5*(0,-0.8)$);
    \draw [gray,very thin] ($(f)+0.3*(0,-0.8)$) -- ($(f)+1.5*(0,-0.8)$);
    \draw [<->,gray,very thin] ($(e)+0.9*(0,-0.8)$) -- ($(f)+0.9*(0,-0.8)$) node[pos =0.5,sloped,below] {2};

    \coordinate (k) at (-1,1);
    \draw [gray,very thin] ($(k)+0.3*(0,0.8)$) -- ($(k)+1.5*(0,0.8)$);
    \draw [gray,very thin] ($(j)+0.3*(0,0.8)$) -- ($(j)+1.5*(0,0.8)$);
    \draw [<->,gray,very thin] ($(k)+0.9*(0,0.8)$) -- ($(j)+0.9*(0,0.8)$) node[pos =0.5,sloped,above] {2};

  \end{tikzpicture}
  \caption{The set $U \subset \mathbb{R}^2$}
  \label{fig:cExample}
\end{figure}
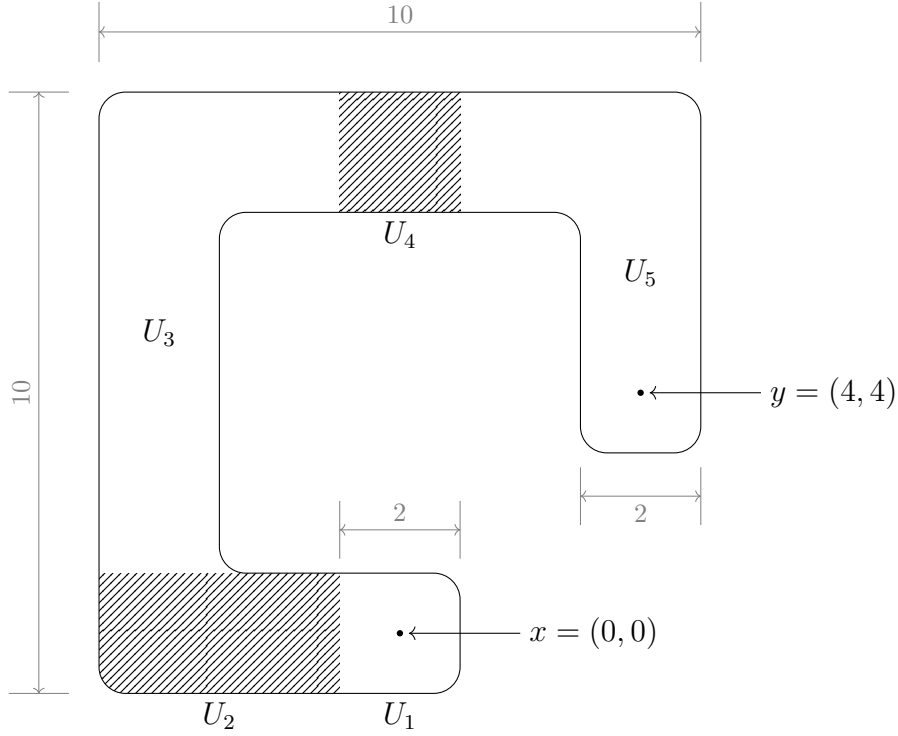

\begin{exm}\label{em:1}\rm
  Let $U_i \subset \mathbb{R}^d$, $i=1,2,\cdots,5$ be sets and $x,y\in \mathbb{R}^d$ be two points as shown in Figure \ref{fig:cExample}. Set
  \begin{equation*}
    U =\bigcup_{i=1}^5 U_i \subset \mathbb{R}^d,
  \end{equation*}
  which is a bounded connected smooth open set in $\R^2$ that does not satisfy condition \hyperref[eq:conH]{$\bf (H_{\gamma})$}
  for any $\gamma \in (0, 1]$ as swapping any coordinate of $x=(0,0)$ by that of $y=(4,4)$ results a point  falling outside $D$.
\end{exm}

\noindent \textbf{Claim:} The lower bound estimate (\ref{eq:pDEst-small-2}) in Theorem \ref{T:1.1}(ii) fails for this open set $U$.
\begin{proof}[Proof of the Claim]
  Fix $x = (0,0)$ and $y=(4,4)$. Since
  \begin{equation*}
    |x^{(i)} - y^{(i)}| = 4 > 0  \qquad \hbox{for }  i=1,2,
  \end{equation*}
  if the inequality (\ref{eq:pDEst-small-2}) in Theorem \ref{T:1.1} (ii) does hold for the set $U$ shown in Figure \ref{fig:cExample}, then there exists $c_1 > 0$ such that
  \begin{equation}\label{eq:cExaLower}
    p_U(t,x,y) \ge c_1 t^2,\qquad t \in (0,1].
  \end{equation}
  We will show that there exists $c_2 > 0$ such that
  \begin{equation}\label{eq:cExaUp}
    p_U(t,x,y) \le c_2 t^3 \quad \hbox{for } t \in (0,1],
  \end{equation}
  which will contradict (\ref{eq:cExaLower}) and finish our claim.

  Indeed, since $U_3 \cap U_1 = \emptyset$, we have by (\ref{eq:XLS}),
  \begin{align*}
     \mathbb{P}_x\left(X_{\tau_{U_1}} \in U_3\right) =&~ \mathbb{E}_x\left[\int_0^{\tau_{U_1}} \sum_{i=1}^2 \int_{\mathbb{R}^d} \1_{U_3}([X_s]^i_\theta) j(X_s^{(i)},\theta)d\theta ds \right].
  \end{align*}
  Note that by the definitions of $U_3$ and $U_1$, we have that for any $z \in U_3$ and $w \in U_1$,
  \begin{equation*}
    z^{(1)} < w^{(1)}\quad \text{ and } \quad z^{(2)} > w^{(2)}.
  \end{equation*}
  Hence, since $X_s \in U_1$ for $s <\tau_{U_1}$, it is not possible that $[X_s]^i_\theta \in U_3$ for any $\theta \in \mathbb{R}$. This together with the above identity implies that
  \begin{equation}\label{eq:cannotjump}
     \mathbb{P}_x\left(X_{\tau_{U_1}} \in U_3\right) =0.
  \end{equation}
  Similarly, one can prove that
  \begin{equation*}
     \mathbb{P}_x\left(X_{\tau_{U_1}} \in U_5\right) =0.
  \end{equation*}
  By the above two identities and the strong Markov property, we obtain, for  almost every  $w \in U$,
  \begin{equation}\label{eq:cExa-pUUp-1}
  \begin{split}
    p_U(t,x,w) =&~ \mathbb{E}_x\left[p_U(t-\tau_{U_1},X^U_{\tau_{U_1}},w); \tau_{U_1} < t\right]\\
    =&~ \mathbb{E}_x\left[p_U(t-\tau_{U_1},X^U_{\tau_{U_1}},w); \tau_{U_1} < t,X^U_{\tau_{U_1}} \in U_2\cup U_4\right]\\
    =&~ \mathbb{E}_x\left[p_U(t-\tau_{U_1},X_{\tau_{U_1}},w); \tau_{U_1} < t,X_{\tau_{U_1}} \in U_2 \cup U_4\right].
  \end{split}
  \end{equation}
  By the positions of $U_2$, $U_4$ and $y$, for any $z \in U_2 \cup U_4$, we have
  \begin{equation*}
    |z^{(i)} - y^{(i)}| \ge 3,\quad i=1,2,
  \end{equation*}
  and, then, by (\ref{eq:pEst}),
  \begin{align*}
    p_U(t-\tau_{U_1},X_{\tau_{U_1}},y)\cdot\1_{\{\tau_{U_1} < t, X_{\tau_{U_1}} \in U_2 \cup U_4\}} \le&~ C_\CTstable \prod_{i=1}^2\frac{t-\tau_{U_1}}{(X_{\tau_{U_1}})^{(i)} - y^{(i)}|^{1+\alpha}}\1_{\{\tau_{U_1} < t,X_{\tau_{U_1}} \in U_2 \cup U_4\}}\\
    \le&~ \frac{C_\CTstable t^2}{3^{2(1+\alpha)}} \1_{\{\tau_{U_1} < t\}}.
  \end{align*}
  Combining (\ref{eq:cExa-pUUp-1}), the continuity of $p_U$ and dominated convergence theorem, we obtain for any $t \in (0,1]$,
  \begin{align*}
    p_U(t,x,y) \le&~ C_\CTstable \mathbb{E}_x\left[\prod_{i=1}^2\frac{t-\tau_{U_1}}{(X_{\tau_{U_1}})^{(i)} - y^{(i)}|^{1+\alpha}}; \tau_{U_1} < t,X_{\tau_{U_1}} \in U_2 \cup U_4\right]\\
    \le&~ \frac{C_\CTstable t^2}{3^{2(1+\alpha)}} \mathbb{P}_x\left(\tau_{U_1} < t\right).
  \end{align*}
  Furthermore, by (\ref{eq:S-ball}) with $r=1$, we have
  \begin{equation*}
    \mathbb{P}_x\left(\tau_{U_1} < t\right) \le \mathbb{P}_x\left(\tau_{B(x,1)} < t\right) \le c_3 t.
  \end{equation*}
  Combining the above two inequalities, we get \eqref{eq:cExaUp} with $c_2 := \frac{c_3C_\CTstable}{3^{2(1+\alpha)}}$.
\end{proof}

\bigskip

In the above example, the domain $U$ is connected but not convex. The ideas used in the above example can be refined to show that the lower bound estimate \eqref{eq:pDEst-small-2} for the Dirichlet heat kernel may still fail for some smooth bounded convex domains.

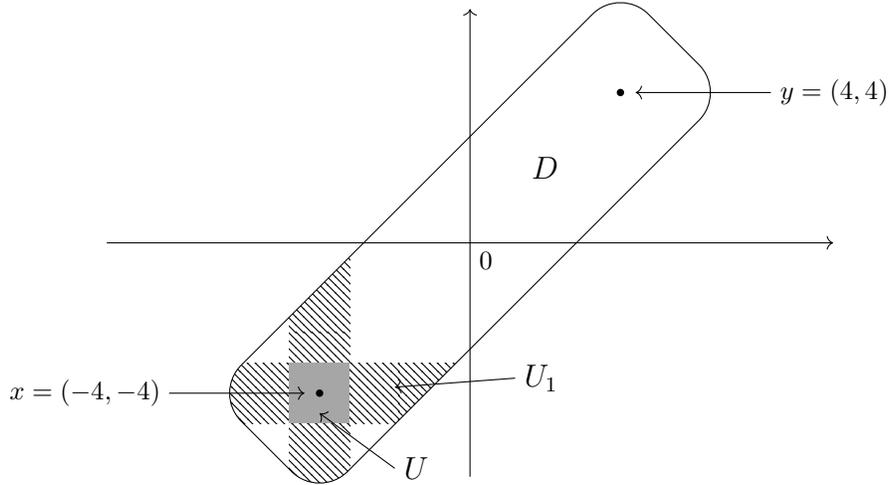
\begin{figure}[htp]
  \begin{tikzpicture}[
    declare function = {
        r = 1;
        s = 0;
        a = 0.2;
        },
      scale = 2,
    ]

    \draw [->] ({-1*r*sqrt(2)-1},0) -- ({r*sqrt(2)+1},0);
    \draw [->] (0,{-1.1*r*sqrt(2)}) -- (0,{1.1*r*sqrt(2)});
    \draw [rotate around={45:(0,0)}, rounded corners = 15pt] ({-1*(r*sqrt(2)+0.5)},{-1*(s+0.5)}) rectangle ({r*sqrt(2)+0.5},s+0.5);
    \node at ({r/2},{r/2}) {\Large$D$};
    \node [below right] {$0$};

     \coordinate (x) at ({-1*r},{-1*r});
    \begin{scope}
      \clip[rotate around={45:(0,0)}, rounded corners = 15pt] ({-1*(r*sqrt(2)+0.5)},{-1*(s+0.5)}) rectangle ({r*sqrt(2)+0.5},s+0.5);
      \fill[pattern= north west lines] ($(x)+(-r,{-1*a})$) rectangle ($(x)+(r,a)$);
    \end{scope}
    \begin{scope}
      \clip[rotate around={45:(0,0)}, rounded corners = 15pt] ({-1*(r*sqrt(2)+0.5)},{-1*(s+0.5)}) rectangle ({r*sqrt(2)+0.5},s+0.5);
      \fill[pattern= north west lines] ($(x)+({-1*a},-r)$) rectangle ($(x)+(a,r)$);
    \end{scope}

    \fill[color = gray!70] ($(x)+({-1*a},{-1*a})$) rectangle ($(x)+(a,a)$);

    \fill (x) circle (0.7pt);
    \draw [->] ($(x)+(-1,0)$) -- ($(x)+(-0.1,0)$);
    \node [left] at ($(x)+(-1,0)$) {$x=(-4,-4)$};

    \draw [->] ($(x)+({r/2},{-r/2})$) -- ($(x)+(0,{-a/1.5})$);
    \node [right] at ($(x)+({r/2},{-r/2})$) {\Large $U$};

    \draw [->] ($(x)+(1.3*r,0.1*r)$) -- ($(x)+(0.5*r,{a/5})$);
    \node [right] at ($(x)+(1.3*r,0.1*r)$) {\Large $U_1$};

    \coordinate (y) at (r,r);
    \fill (y) circle (0.7pt);
    \draw [->] ($(y)+(1,0)$) -- ($(y)+(0.1,0)$);
    \node [right] at ($(y)+(1,0)$) {$y=(4,4)$};
  \end{tikzpicture}
  \caption{The convex set $D\subset \mathbb{R}^2$ that does not satisfy condition $\bf (H_{\gamma})$ for any $\gamma \in (0,1)$}
  \label{fig:convexD}
\end{figure}

\begin{exm}\label{em:2}\rm
  Let $D$ be the tilted rectangle with rounded corners shown in Figure \ref{fig:convexD}. The points $x,y$ has the coordinates $(-4,-4)$ and $(4,4)$ respectively. The connected open set $D\subset \mathbb{R}^2$ that does not satisfy condition $\bf (H_{\gamma})$ for any $\gamma \in (0,1)$ as swapping any coordinate of $x$ by that of $y$ results a point falling outside $D$.
\end{exm}

\noindent \textbf{Claim:} The lower bound estimate  \eqref{eq:pDEst-small-2}  in Theorem \ref{T:1.1}(ii) fails for this convex open set $D$.

\begin{proof}[Proof of the Claim]
  The proof is almost the same to that in Example \ref{em:1}. Indeed, since
  \begin{equation*}
    |x^{(i)} - y^{(i)}| = 8 > 0,\qquad i=1,2,
  \end{equation*}
  if the inequality (\ref{eq:pDEst-small-2}) in Theorem \ref{T:1.1} (ii) does hold for the set $D$ shown in Figure \ref{fig:cExample}, then there exists $c_1 > 0$ such that
  \begin{equation}\label{eq:em2-Lower}
    p_U(t,x,y) \ge c_1 t^2,\qquad t \in (0,1].
  \end{equation}
  We will show that in fact  there exists $c_2 > 0$ such that
  \begin{equation}\label{eq:em2-Up}
    p_U(t,x,y) \le c_2 t^3,\qquad t \in (0,1],
  \end{equation}
  which will establish the claim by contradiction.

  Let $U := Q(x,1)$ be the open square with gray color and $U_1$ be the dash region in Figure \ref{fig:convexD}. Set $U_2 := D\setminus (U\cup U_1)$. Since $U \cap U_2 = \emptyset$, by (\ref{eq:XLS}), we obtain
  \begin{align*}
     \mathbb{P}_x\left(X_{\tau_{U}} \in U_2\right) =&~ \mathbb{E}_x\left[\int_0^{\tau_{U}} \sum_{i=1}^2 \int_{\mathbb{R}^d} \1_{U_2}([X_s]^i_\theta) j(X_s^{(i)},\theta)d\theta ds \right].
  \end{align*}
  Note that by the positions of $U$ and $U_2$, we have that for any $z \in U$ and $w \in U_2$,
  \begin{equation*}
    z^{(1)} \neq w^{(1)}\quad \text{ and } \quad z^{(2)} \neq w^{(2)}.
  \end{equation*}
  Hence, since $X_s \in U$ for $s <\tau_{U}$, it is not possible that $[X_s]^i_\theta \in U_2$ for any $\theta \in \mathbb{R}$. This together with the above identity implies that
  \begin{equation}\label{eq:em2-cannotjump}
     \mathbb{P}_x\left(X_{\tau_{U}} \in U_2\right) =0.
  \end{equation}
  Moreover, by the above identity and the strong Markov property, we obtain, for almost every $w \in U$, we have
  \begin{equation}\label{eq:em2-pUUp-1}
  \begin{split}
    p_D(t,x,w) =&~ \mathbb{E}_x\left[p_D(t-\tau_U,X^D_{\tau_U},w); \tau_U < t\right]\\
    =&~ \mathbb{E}_x\left[p_D(t-\tau_U,X^D_{\tau_U},w); \tau_U < t,X^D_{\tau_U} \in U_1\right]\\
    =&~ \mathbb{E}_x\left[p_D(t-\tau_U,X_{\tau_U},w); \tau_U < t,X_{\tau_U} \in U_1\right].
  \end{split}
  \end{equation}
  By the positions of $U_1$ and $y$, for any $z \in U_1$,
  \begin{equation*}
    |z^{(i)} - y^{(i)}| > 6,\quad i=1,2,
  \end{equation*}
  and, then, by (\ref{eq:pEst}),
  \begin{align*}
    p_D(t-\tau_U,X_{\tau_U},y)\cdot\1_{\{\tau_U < t,X_{\tau_U} \in U_1\}} \le&~ C_\CTstable \prod_{i=1}^2\frac{t-\tau_U}{(X_{\tau_U})^{(i)} - y^{(i)}|^{1+\alpha}}\1_{\{\tau_U < t,X_{\tau_U} \in U_1\}}\\
    \le&~ \frac{C_\CTstable t^2}{3^{2(1+\alpha)}} \1_{\{\tau_U < t\}}.
  \end{align*}
  This together with (\ref{eq:em2-pUUp-1}), the continuity of $p_D$ and dominated convergence theorem, yields for any $t \in (0,1]$,
  \begin{align*}
    p_D(t,x,y) \le&~ C_\CTstable \mathbb{E}_x\left[\prod_{i=1}^2\frac{t-\tau_U}{(X_{\tau_U})^{(i)} - y^{(i)}|^{1+\alpha}}; \tau_U < t,X_{\tau_U} \in U_1\right]\\
    \le&~ \frac{C_\CTstable t^2}{3^{2(1+\alpha)}} \mathbb{P}_x\left(\tau_U < t\right).
  \end{align*}
  Furthermore, by (\ref{eq:S-ball}) with $r=1$, we have
  \begin{equation*}
    \mathbb{P}_x\left(\tau_U < t\right) \le \mathbb{P}_x\left(\tau_{B(x,1)} < t\right) \le c_3 t.
  \end{equation*}
  Combining the above two inequalities, we obtain (\ref{eq:cExaUp}) with $c_2 := \frac{c_3C_\CTstable}{3^{2(1+\alpha)}}$. We finish the proof.
\end{proof}

\bigskip

The following is an example of a bounded smooth open set $D\subset \R^2$ that does not satisfy the irreducibility condition \eqref{e:1.13} but for which we can derive two-sided sharp estimates for the Dirichlet heat kernel in $D$.

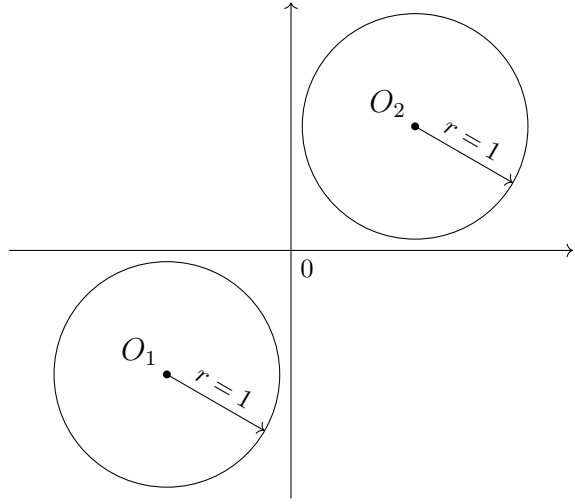
\begin{figure}[htp]
   \begin{tikzpicture}[
        scale = 1.5
    ]
    \draw [->] (-2.5,0) -- (2.5,0);
    \draw [->] (0,-2.2) -- (0,2.2);
    \node [below right] {$0$};

    \draw (-1.1,-1.1) circle (1);
    \fill (-1.1,-1.1) circle (1pt);
    \draw [->] (-1.1,-1.1) -- (-1.1+0.866,-1.1-0.5)  node[pos =0.5,sloped,above] {$r=1$};
    \node [above left] at (-1.1, -1.1) {\Large $O_1$};
    \draw (1.1,1.1) circle (1);
    \fill (1.1,1.1) circle (1pt);
    \draw [->] (1.1,1.1) -- (1.1+0.866,1.1-0.5) node[pos =0.5,sloped,above] {$r=1$};
    \node [above left] at (1.1, 1.1) {\Large $O_2$};
  \end{tikzpicture}
  \caption{The set $D := B(O_1,r) \cup B(O_2,r)$ with $r = 1$ on $\mathbb{R}^2$}
  \label{fig:diagonal}
\end{figure}

\begin{exm}\label{em:3}\rm
  Let $r>0$ and $D$ be the union of two disjoint balls sitting in diagonal quadrants:
  \begin{equation}\label{eq:D2-def}
    D := B(O_1,r) \cup B(O_2,r),
  \end{equation}
  where the two points $O_1, O_2 \in \mathbb{R}^d$ satisfy
  \begin{equation*}
    O_1^{(i)} < -r\quad \text{ and } \quad O_2^{(i)} > r  \quad  \hbox{for } 1\le i \le d,
  \end{equation*}
  (see Figure \ref{fig:diagonal}).
\end{exm}

The open set $D$ clearly does not satisfy the condition \eqref{e:1.13} as for any $x \in O_1$ and $y\in O_2$, swapping any coordinate of $x$ by that of $y$ results a point falling outside $D$. So $X^D$ is not irreducible. It follows from Corollary \ref{C:1.4} that the following holds with $B_1 := B(O_1,r)$  and $B_2:= B(O_2,r)$.
\begin{enumerate}
  \item[\rm (i)] $p_D(t,x,y) = 0$ for all $ t > 0$ and $x, y$ that are not in the same connected component of $D$.

  \item[\rm (ii)] For $i=1, 2$, $p_D(t,x,y) = p_{B_i}(t,x,y) $ for all $t > 0$ and $x,y \in B_i$, and $p_{B_i}(t,x,y)$ has the two-sided estimates given by Theorem \ref{T:1.1} with $B_i$ in place of $D$ there.
\end{enumerate}

\medskip

  The above example clearly can be extended to more general open sets $D$ that is the union of two disjoint $C^{1.1}$-smooth connected open subsets $O_1$ and $O_2$ that at most one of them has non-empty intersection with any line that is parallel to the coordinate axes.

\medskip

  We conclude this paper by presenting the proof of Theorem \ref{T:1.6}, using the techniques from Section  \ref{S:4}.

\medskip

\begin{proof}[Proof of Theorem \ref{T:1.6}.]
  Fix $T > 0$. Note that $p_D(t,x,y)$ is symmetric. For any $x,y \in D$, if $y$ does not belong to $A_1\cup A_4$, then $x,y$ belong to the case in (i).

  (i). For any $x \in A_i$, $y \in A_j$ with $|i-j|\le 2$, we have that $x,y$ belong to some $C^{1,1}$ open subset $U$ of $D$ and $U$ satisfies the condition \hyperref[eq:conH]{$\bf (H_{1})$}. In this case, we also have
  \begin{equation}\label{eq:dDx=dUx-square}
    \delta_D(x) = \delta_U(x),\quad \delta_D(y) = \delta_U(y).
  \end{equation}
  For example, for $x\in A_1$ and $y \in A_3$, we can set $U = A_1 \cup A_2 \cup A_3$. For $x \in A_4$ and $y \in A_2$, we can set $U = A_2 \cup A_3\cup A_4$.

  The upper bound of $p_D(t,x,y)$ in (\ref{eq:pDEst-exm-square-1}) follows directly from Theorem \ref{T:1.1} (i) since $D$ is a $C^{1,1}$ open set. For the lower bound of $p_D(t,x,y)$ in (\ref{eq:pDEst-exm-square-1}), note that $U\subset D$ is a $C^{1,1}$ open set satisfying the condition \hyperref[eq:conH]{$\bf (H_{1})$}. By Theorem \ref{T:1.1} (ii) and (\ref{eq:dDx=dUx-square}), we have
  \begin{align*}
    p_D(t,x,y) \ge&~ p_U(t,x,y)\\
    \ge&~ c_1 \left(1\wedge \frac{\delta_U(x)^{\alpha/2}}{\sqrt{t}}\right)\left(1\wedge \frac{\delta_U(y)^{\alpha/2}}{\sqrt{t}}\right)p(t,x,y)\\
    =&~c_1\left(1\wedge \frac{\delta_D(x)^{\alpha/2}}{\sqrt{t}}\right)\left(1\wedge \frac{\delta_D(y)^{\alpha/2}}{\sqrt{t}}\right)p(t,x,y).
  \end{align*}

  (ii). Fix $x =(x^{(1)}, x^{(2)}) \in A_1$ and $y =(y^{(1)}, y^{(2)}) \in A_4$. Similar to (\ref{eq:cannotjump}), one can apply (\ref{eq:XLS}) to prove that
  \begin{equation*}
    \mathbb{P}_x\left(X_{\tau_{A_1}} \in A_3 \cup A_4\right) = 0.
  \end{equation*}
  Hence, by (\ref{eq:pEst}), the function
  \begin{equation*}
    p_D(t-\tau_{A_1},X_{\tau_{A_1}},y)\1_{\{\tau_{A_1} < t\}}
  \end{equation*}
  is uniformly bounded in $y \in A_4$, and then, by the strong Markov property of $X$, the continuity of $p_D$ and dominated convergence theorem, we have
  \begin{align*}
    p_D(t,x,y) =&~ \mathbb{E}_x\left[p_D(t-\tau_{A_1},X_{\tau_{A_1}},y);\tau_{A_1} < t\right]\\
    =&~ \mathbb{E}_x\left[p_D(t-\tau_{A_1},X_{\tau_{A_1}},y);\tau_{A_1} < t,X_{\tau_{A_1}} \in A_2\right].
  \end{align*}
  Furthermore, by (\ref{eq:XLS}) again, we have for $t > 0$,
  \begin{align}\label{eq:4squares-SMP}
    p_D(t,x,y) =&~ \mathbb{E}_x\left[p_D(t-\tau_{A_1},X_{\tau_{A_1}},y);\tau_{A_1} < t,X_{\tau_{A_1}} \in A_2\right]\notag\\
    =&~\mathbb{E}_x\left[\int_0^t \1_{\{s<\tau_{A_1}\}} \cdot\left(\sum_{i=1}^2\int_{\mathbb{R}} \1_{\{[X_s]_\theta^i \in A_2\}}\cdot p_D(t-s,[X_s]_\theta^i,y)j(X_s^{(i)},\theta)d\theta\right)ds\right]\notag\\
    =&~\int_0^t \int_{A_1} p_{A_1}(s,x,u)\left(\int_{\mathbb{R}} \1_{\{[u]_\theta^1 \in A_2\}}\cdot p_D(t-s,[u]^1_\theta,y) j(u^{(1)},\theta)d\theta\right) du ds.
  \end{align}
  Note that
  \begin{enumerate}
    \item[$\bullet$] for any $u = (u^{(1)}, u^{(2)}) \in A_1$ and $\theta \in \mathbb{R}$ with $[u]^1_\theta \in A_2$, we have
        $$
          |u^{(1)} - \theta| \geq 1, \quad  |y^{(1)} - \theta| \geq 1   \quad \hbox{and} \quad  |u^{(2)} - y^{(2)}|\geq 1;
        $$

    \item[$\bullet$] by \eqref{eq:pDEst-exm-square-1}, for any $[u]^1_\theta \in A_2$ and $0<s<t$, we have
        \begin{align*}
          p_D(t-s,[u]^1_\theta,y) \le&~ C_{\CTstable} C_{\CTexmDa} \left(1\wedge \frac{\delta_D(y)^{\alpha/2}}{\sqrt{t-s}}\right) \left((t-s)^{-\frac{1}{\alpha}}\wedge\frac{t-s}{|\theta-y^{(1)}|^{1+\alpha}}\right) \left((t-s)^{-\frac{1}{\alpha}}\wedge\frac{t-s}{|u^{(2)}-y^{(2)}|^{1+\alpha}}\right)\\
          \le&~ C_{\CTstable} C_{\CTexmDa} \left(1\wedge \frac{\delta_D(y)^{\alpha/2}} {\sqrt{t-s}}\right) t^2.
        \end{align*}
  \end{enumerate}
  Combining the above two inequalities with (\ref{eq:4squares-SMP}) and Theorem \ref{T:1.1} (i), we obtain, for all $t \in (0,T]$,
  \begin{align*}
    p_D(t,x,y) \le&~ c_2\!\!\int_0^t \!\!\int_{A_1}\!\! p_{A_1}(s,x,u) \left(\int_{O_2^{(1)}-1}^{O_2^{(1)}+1} d\theta\right) \left(1\wedge \frac{\delta_D(y)^{\alpha/2}}{\sqrt{t-s}}\right) t^2 du ds\\
    \le&~ 2c_2 t^2 \int_0^t \left(\int_{A_1} p(s,x,u) du \right) \left(1\wedge\frac{\delta_{A_1}(x)^{\alpha/2}}{\sqrt{s}}\right) \left(1\wedge \frac{\delta_D(y)^{\alpha/2}}{\sqrt{t-s}}\right) ds        \\
    \le&~ 2c_2 t^2 \left( \int_0^{t/2} + \int_{t/2}^t\right) \left(1\wedge\frac{\delta_{D}(x)^{\alpha/2}}{\sqrt{s}}\right) \left(1\wedge \frac{\delta_D(y)^{\alpha/2}}{\sqrt{t-s}}\right) ds  \\
    =&~ c_3 t^3 \left(1\wedge\frac{\delta_D(x)^{\alpha/2}}{\sqrt{t}}\right) \left(1\wedge\frac{\delta_D(y)^{\alpha/2}}{\sqrt{t}}\right).
  \end{align*}
  which is exactly the upper bound in \eqref{eq:pDEst-exm-square-2}.

\medskip

  We apply (\ref{eq:4squares-SMP}) to establish the lower bound of $p_D(t,x,y)$. Fix $t \in (0,T]$. Let $\delta < (16T^{1/\alpha})^{-1}$ and $Q_x \in \partial D$ be such that $\delta_D(x) = |x-Q_x|$. Set
  \begin{equation*}
    x_0 =
    \begin{cases}
      Q_x + \frac{2\delta t^{1/\alpha}}{|x-Q_x|}(x-Q_x), &\text{ if } \delta_D(x) < 2\delta t^{1/\alpha},\\
      x, &\text{ if } \delta_D(x) \ge 2\delta t^{1/\alpha}.
    \end{cases}
  \end{equation*}
  Define
  \begin{equation*}
    E_x = B(x_0,\delta t^{1/\alpha}).
  \end{equation*}
  Note that $E_x \subset A_1 \subset D$.
  Observe that
  \begin{enumerate}
    \item[$\bullet$] for any $s \in (0,t)$ and $u \in E_x$, we have $\delta_D(u) \ge \delta t^{1/\alpha}$ and so  by Theorem \ref{T:1.1} (ii),
        \begin{align*}
          & p_{A_1}(s,x,u) \\
          \ge&~ c_4\left(1\wedge\frac{\delta_{A_1}(x)^{\alpha/2}}{\sqrt{s}}\right) \left(1\wedge\frac{\delta_{A_1}(u)^{\alpha/2}}{\sqrt{s}}\right) \left(s^{-\frac{1}{\alpha}}\wedge \frac{s}{|x^{(1)} - u^{(1)}|^{1+\alpha}}\right)\left(s^{-\frac{1}{\alpha}}\wedge \frac{s}{|x^{(2)} - u^{(2)}|^{1+\alpha}}\right)\\
          \ge&~ c_4\left(1\wedge\frac{\delta_D(x)^{\alpha/2}}{\sqrt{t}}\right)(1\wedge\delta^{\alpha/2})
           \left(s^{-\frac{1}{\alpha}}\wedge \frac{s}{(3\delta t^{1/\alpha})^{1+\alpha}}\right)^2\\
                    \ge&~ c_5\left(1\wedge\frac{\delta_D(x)^{\alpha/2}}{ \sqrt{t}}\right) s^2     t^{-2-2/\alpha};
        \end{align*}
    \item[$\bullet$] for any $u \in E_x$ and $\theta \in [O_2^{(1)}-\frac{1}{2},  O_2^{(1)} +\frac{1}{2}]$, we have $[u]^1_\theta \in A_2$, $|u^{(1)}-\theta| \le 3$ and
        \begin{equation*}
          \delta_D([u]^1_\theta) \ge c_6 t^{1/\alpha}.
        \end{equation*}
        Furthermore, by \eqref{eq:pDEst-exm-square-1},   we have for any $s \in (0,\frac{t}{2}]$,
        \begin{align*}
          &~p_D(t-s,[u]^1_\theta,y)  \\
          \ge&~ c_7 \left(1\wedge \frac{\delta_D([u]^1_\theta)^{\alpha/2}}{\sqrt{t-s}}\right) \left(1\wedge \frac{\delta_D(y)^{\alpha/2}}{\sqrt{t-s}}\right) \left((t-s)^{-\frac{1}{\alpha}}\wedge\frac{t-s}{|\theta-y^{(1)}|^{1+\alpha}}\right) \left((t-s)^{-\frac{1}{\alpha}}\wedge\frac{t-s}{|u^{(2)}-y^{(2)}|^{1+\alpha}}\right)\\
          \ge&~  c_8 \left(1\wedge\frac{\delta_D(y)^{\alpha/2}}{\sqrt{t}}\right) t^2.
        \end{align*}
  \end{enumerate}
  Combining the above two inequalities with (\ref{eq:4squares-SMP}), we have for any $t \in (0,T]$,
  \begin{align*}
    p_D(t,x,y) \ge&~\int_0^{t/2} \int_{E_x} p_{A_1}(s,x,u)\left(\int_{\mathbb{R}} \1_{\{\theta \in [O_2^{(1)} -\frac{1}{2}, O_2^{(1)} +\frac{1}{2}]\}} \cdot p_D(t-s,[u]^1_\theta,y) j(u^{(1)},\theta)d\theta\right) du ds\\
    \ge&~ c_9 \int_0^{t/2} \int_{E_x} \left(1\wedge\frac{\delta_D(x)^{\alpha/2}}{\sqrt{t}}\right) s^2 t^{-2-2/\alpha}.
    \left( \int_{O_2^{(1)} -\frac{1}{2}}^{O_2^{(1)} +\frac{1}{2}} \left(1\wedge \frac{\delta_D(y)^{\alpha/2}}{\sqrt{t}}\right) t^2 \frac{1}{|u^{(1)} - \theta|^{1+\alpha}} d\theta \right) du ds\\
    \ge&~ c_{10} t^{-2/\alpha} \left(1\wedge\frac{\delta_D(x)^{\alpha/2}}{\sqrt{t}}\right) \left(1\wedge \frac{\delta_D(y)^{\alpha/2}}{\sqrt{t}}\right) |E_x| \int_0^{t/2} s^2 ds  \\
    =&~ c_{11} t^{-2/\alpha} \left(1\wedge\frac{\delta_D(x)^{\alpha/2}}{\sqrt{t}}\right) \left(1\wedge \frac{\delta_D(y)^{\alpha/2}}{\sqrt{t}}\right) (\delta t^{1/\alpha})^2 \, t^3  \\
    =& ~c_{12} t^3 \left(1\wedge\frac{\delta_D(x)^{\alpha/2}}{\sqrt{t}}\right) \left(1\wedge \frac{\delta_D(y)^{\alpha/2}}{\sqrt{t}}\right),
  \end{align*}
  which is the lower bound in (\ref{eq:pDEst-exm-square-2}).

  On the other hand, for $x=(x^{(1)}, x^{(2)})\in A_1$ and $y=(y^{(1)}, y^{(2)})\in A_4$, we have
  $$
    1\leq |x^{(k)}-y^{(k)} |\leq 3+3+2=8 \quad \hbox{for } k=1, 2.
  $$
  Hence by \eqref{eq:pEst},
  $$
    p(t,x,y) \stackrel{c_{13}}{\asymp} \prod_{k=1}^2 \left(t^{-1/\alpha}\wedge\frac{t}{|x^{(k)}-y^{(k)}|^{1+\alpha}}\right) \stackrel{c_{14}}{\asymp} t^2 \quad \hbox{for any } t\in (0, T],
  $$
  where $c_{13}$ and $c_{14}$ are positive constants that depend only on $(\alpha, T)$. So we get \eqref{e:6.12} from (\ref{eq:pDEst-exm-square-2}). This completes the proof of the theorem.
\end{proof}

  \textbf{Acknowledgements} We thank Russell Lyons for suggesting the name `rectilinear stable process'
over our previous terminology `cylindrical stable process'.  We also thank the  referee for helpful  comments.

%\bibliographystyle{abbrv}
%\bibliography{ref}

\begin{thebibliography}{10}

\bibitem{BaeKim.2020.PA661}
J.~Bae and P.~Kim.
\newblock On estimates of transition density for subordinate {B}rownian motions
  with {G}aussian components in {$C^{1,1}$}-open sets.
\newblock {\em Potential Anal.}, 52(4):661--687, 2020.

\bibitem{BassChen.2010.MZ489}
R.~F. Bass and Z.-Q. Chen.
\newblock Regularity of harmonic functions for a class of singular stable-like
  processes.
\newblock {\em Math. Z.}, 266(3):489--503, 2010.

\bibitem{BendikovGrigoryanHuEtAl.2021.ASNSPCS5399}
A.~Bendikov, A.~Grigor'yan, E.~Hu, and J.~Hu.
\newblock Heat kernels and non-local {D}irichlet forms on ultrametric spaces.
\newblock {\em Ann. Sc. Norm. Super. Pisa Cl. Sci. (5)}, 22(1):399--461, 2021.

\bibitem{BlumenthalGetoor.1960.TAMS263}
R.~M. Blumenthal and R.~K. Getoor.
\newblock Some theorems on stable processes.
\newblock {\em Trans. Amer. Math. Soc.}, 95:263--273, 1960.

\bibitem{BlumenthalGetoor.1968.313}
R.~M. Blumenthal and R.~K. Getoor.
\newblock {\em Markov processes and potential theory}.
\newblock Pure and Applied Mathematics, Vol. 29. Academic Press, New
  York-London, 1968.

\bibitem{BogdanBurdzyChen.2003.PTRF89}
K.~Bogdan, K.~Burdzy, and Z.-Q. Chen.
\newblock Censored stable processes.
\newblock {\em Probab. Theory Related Fields}, 127(1):89--152, 2003.

\bibitem{BogdanGrzywnyRyznar.2010.AP1901}
K.~Bogdan, T.~Grzywny, and M.~Ryznar.
\newblock Heat kernel estimates for the fractional {L}aplacian with {D}irichlet
  conditions.
\newblock {\em Ann. Probab.}, 38(5):1901--1923, 2010.

\bibitem{BogdanGrzywnyRyznar.2014.SPA3612}
K.~Bogdan, T.~Grzywny, and M.~Ryznar.
\newblock Dirichlet heat kernel for unimodal {L}\'{e}vy processes.
\newblock {\em Stochastic Process. Appl.}, 124(11):3612--3650, 2014.

\bibitem{ChenKimWang.2022.MA373}
X.~Chen, P.~Kim, and J.~Wang.
\newblock Two-sided {D}irichlet heat kernel estimates of symmetric stable
  processes on horn-shaped regions.
\newblock {\em Math. Ann.}, 384(1-2):373--418, 2022.

\bibitem{ChenCroydonKumagai.2015.AP1594}
Z.-Q. Chen, D.~A. Croydon, and T.~Kumagai.
\newblock Quenched invariance principles for random walks and elliptic
  diffusions in random media with boundary.
\newblock {\em Ann. Probab.}, 43(4):1594--1642, 2015.

\bibitem{ChenKimSong.2010.JEMSJ1307}
Z.-Q. Chen, P.~Kim, and R.~Song.
\newblock Heat kernel estimates for the {D}irichlet fractional {L}aplacian.
\newblock {\em J. Eur. Math. Soc. (JEMS)}, 12(5):1307--1329, 2010.

\bibitem{ChenKimSong.2011.JLMS258}
Z.-Q. Chen, P.~Kim, and R.~Song.
\newblock Heat kernel estimates for {$\Delta+\Delta^{\alpha/2}$} in {$C^{1,1}$}
  open sets.
\newblock {\em J. Lond. Math. Soc. (2)}, 84(1):58--80, 2011.

\bibitem{ChenKimSong.2012.AP2483}
Z.-Q. Chen, P.~Kim, and R.~Song.
\newblock Dirichlet heat kernel estimates for fractional {L}aplacian with
  gradient perturbation.
\newblock {\em Ann. Probab.}, 40(6):2483--2538, 2012.

\bibitem{ChenKimSong.2014.PLMS390}
Z.-Q. Chen, P.~Kim, and R.~Song.
\newblock Dirichlet heat kernel estimates for rotationally symmetric {L}\'{e}vy
  processes.
\newblock {\em Proc. Lond. Math. Soc. (3)}, 109(1):90--120, 2014.

\bibitem{ChenKimSong.2016.JRAM111}
Z.-Q. Chen, P.~Kim, and R.~Song.
\newblock Dirichlet heat kernel estimates for subordinate {B}rownian motions
  with {G}aussian components.
\newblock {\em J. Reine Angew. Math.}, 711:111--138, 2016.

\bibitem{ChenKimSongEtAl.2012.TAMS4169}
Z.-Q. Chen, P.~Kim, R.~Song, and Z.~Vondra{\v{c}}ek.
\newblock Boundary {H}arnack principle for {$\Delta+\Delta^{\alpha/2}$}.
\newblock {\em Trans. Amer. Math. Soc.}, 364(8):4169--4205, 2012.

\bibitem{ChenKumagai.2003.SPA27}
Z.-Q. Chen and T.~Kumagai.
\newblock Heat kernel estimates for stable-like processes on {$d$}-sets.
\newblock {\em Stochastic Process. Appl.}, 108(1):27--62, 2003.

\bibitem{ChenKumagai.2008.PTRF277}
Z.-Q. Chen and T.~Kumagai.
\newblock Heat kernel estimates for jump processes of mixed types on metric
  measure spaces.
\newblock {\em Probab. Theory Related Fields}, 140(1-2):277--317, 2008.

\bibitem{ChenTokle.2011.PTRF373}
Z.-Q. Chen and J.~Tokle.
\newblock Global heat kernel estimates for fractional {L}aplacians in unbounded
  open sets.
\newblock {\em Probab. Theory Related Fields}, 149(3-4):373--395, 2011.

\bibitem{Cho.2006.PA387}
S.~Cho.
\newblock Two-sided global estimates of the {G}reen's function of parabolic
  equations.
\newblock {\em Potential Anal.}, 25(4):387--398, 2006.

\bibitem{ChoKangKim.2021.JLMS2823}
S.~Cho, J.~Kang, and P.~Kim.
\newblock Estimates of {D}irichlet heat kernels for unimodal {L}\'{e}vy
  processes with low intensity of small jumps.
\newblock {\em J. Lond. Math. Soc. (2)}, 104(2):823--864, 2021.

\bibitem{Chung.1986.63}
K.~L. Chung.
\newblock Doubly-{F}eller process with multiplicative functional.
\newblock In {\em Seminar on stochastic processes, 1985 ({G}ainesville, {F}la.,
  1985)}, volume~12 of {\em Progr. Probab. Statist.}, pages 63--78.
  Birkh\"{a}user Boston, Boston, MA, 1986.

\bibitem{ChungZhao.1995.287}
K.~L. Chung and Z.~X. Zhao.
\newblock {\em From {B}rownian motion to {S}chr\"odinger's equation}, volume
  312 of {\em Grundlehren der Mathematischen Wissenschaften [Fundamental
  Principles of Mathematical Sciences]}.
\newblock Springer-Verlag, Berlin, 1995.

\bibitem{GrigoryanHuHu.2017.JFA3311}
A.~Grigor'yan, E.~Hu, and J.~Hu.
\newblock Lower estimates of heat kernels for non-local {D}irichlet forms on
  metric measure spaces.
\newblock {\em J. Funct. Anal.}, 272(8):3311--3346, 2017.

\bibitem{GrigoryanHuHu.2018.AM433}
A.~Grigor'yan, E.~Hu, and J.~Hu.
\newblock Two-sided estimates of heat kernels of jump type {D}irichlet forms.
\newblock {\em Adv. Math.}, 330:433--515, 2018.

\bibitem{GrzywnyKimKim.2020.SPA431}
T.~Grzywny, K.-Y. Kim, and P.~Kim.
\newblock Estimates of {D}irichlet heat kernel for symmetric {M}arkov
  processes.
\newblock {\em Stochastic Process. Appl.}, 130(1):431--470, 2020.

\bibitem{KassmannKimKumagai.2022.JMPA91}
M.~Kassmann, K.-Y. Kim, and T.~Kumagai.
\newblock Heat kernel bounds for nonlocal operators with singular kernels.
\newblock {\em J. Math. Pures Appl. (9)}, 164:1--26, 2022.

\bibitem{KimWang.2022.SPA165}
K.-Y. Kim and L.~Wang.
\newblock Heat kernel bounds for a large class of {M}arkov process with
  singular jump.
\newblock {\em Stochastic Process. Appl.}, 145:165--203, 2022.

\bibitem{KimMimica.2018.EJP45}
P.~Kim and A.~Mimica.
\newblock Estimates of {D}irichlet heat kernels for subordinate {B}rownian
  motions.
\newblock {\em Electron. J. Probab.}, 23:Paper No. 64, 45, 2018.

\bibitem{Ros-OtonSerra.2014.JMPA}
X.~Ros-Oton and J.~Serra.
\newblock The {D}irichlet problem for the fractional {L}aplacian: regularity up
  to the boundary.
\newblock {\em J. Math. Pures Appl. (9)}, 101(3):275--302, 2014.

\bibitem{Ros-OtonSerra.2016.Duke}
X.~Ros-Oton and J.~Serra.
\newblock Boundary regularity for fully nonlinear integro-differential
  equations.
\newblock {\em Duke Math. J.}, 165(11):2079--2154, 2016.

\bibitem{Ros-OtonSerra.2016.JDE}
X.~Ros-Oton and J.~Serra.
\newblock Regularity theory for general stable operators.
\newblock {\em J. Differential Equations}, 260(12):8675--8715, 2016.

\bibitem{Ros-OtonValdinoci.2016.Adv}
X.~Ros-Oton and E.~Valdinoci.
\newblock The {D}irichlet problem for nonlocal operators with singular kernels:
  convex and nonconvex domains.
\newblock {\em Adv. Math.}, 288:732--790, 2016.

\bibitem{Schaefer.1974.376}
H.~H. Schaefer.
\newblock {\em Banach lattices and positive operators}.
\newblock Springer-Verlag, New York-Heidelberg, 1974.
\newblock Die Grundlehren der mathematischen Wissenschaften, Band 215.

\end{thebibliography}

\vskip 0.2truein

\end{document}